\documentclass[12pt]{amsart}
\usepackage{amsmath,amssymb,amscd,array, mathrsfs }

\def\url#1{\expandafter\string\csname #1\endcsname}

\def\mmat #1,#2,#3,#4,{\text{\small\arraycolsep=3pt $
\begin{pmatrix}#1&#2\\#3&#4\end{pmatrix}$}}

\usepackage{lscape}

\usepackage{DLdef1}

\rmnameii{etr}{etr}
\rmnameii{evv}{ev}

\begin{document}
\title[Double extensions of Lie superalgebras]{Double extensions of Lie superalgebras in characteristic 2 with nondegenerate invariant supersymmetric bilinear form}
\author{Sa\"id Benayadi}
\address{D\'epartement de Math\'ematiques, Institut Elie Cartan de Lorraine, Universit\'e de Lorraine, CNRS-UMR 7502, Ile du Saulcy, F-57012 Metz cedex 01, France.}

\email{said.benayadi@univ-lorraine.fr }

\author{Sofiane Bouarroudj}

\address{Division of Science and Mathematics, New York University Abu Dhabi, P.O. Box 129188, Abu Dhabi, United Arab Emirates.}
\email{sofiane.bouarroudj@nyu.edu}

\thanks{SB was supported by the grant NYUAD-065.}


\date{\today}



\begin{abstract} A Lie (super)algebra with a non-degenerate invariant symmetric bilinear form will be called a NIS-Lie (super)algebra. The double extension of a NIS-Lie (super)algebra is the result of simultaneously adding to it a central element and an outer derivation so that the larger algebra has also a NIS. Affine loop algebras, Lie (super)algebras with symmetrizable Cartan matrix over any field, Manin triples, symplectic reflection (super)algebras are among the Lie (super)algebras suitable to be doubly extended.

We consider double extensions of Lie superalgebras in characteristic 2, and concentrate on peculiarities of these notions related with the possibility for the bilinear form, the center, and the derivation to be odd.  Two Lie superalgebras we discovered by this method are indigenous to the characteristic 2.

\end{abstract}

\maketitle
 
\section {Introduction}
\ssec{NIS-algebras} Hereafter a Lie algebra with a \textbf{Non-degenerate Invariant Symmetric} bilinear form $B$ will be called a \textit{NIS-Lie algebra}\footnote{Initially, the term \textit{quadratic} was used to single out these algebras, see \cite{ABB}, \cite{ABBQ}, \cite{BB}, \cite{BBB}, \cite{B}. The term had been already occupied: V.~Drinfeld \cite{Dr1} introduced \textit{quadratic algebras} --- the ones with quadratic relations --- in a paper published in 1986 (and never translated from Russian, as far as we know); Yu.~Manin used Drinfeld's term in his book, \cite{Ma}. To avoid confusion, we decided from now on to call the algebras we are studying after their main properties. Besides, the form whose properties are vital for us is \textbf{bilinear}, not quadratic.}. If the Lie algebra is simple and finite-dimensional over $\Cee$, then such a form $B$ is induced by the trace in any irreducible module of dimension $>1$ and is a multiple of the Killing form induced by the trace in the adjoint representation. 

The Killing form, however, becomes degenerate if the simple Lie algebra is defined over a field of characteristic $p>0$, see \cite{S}, or on simple Lie superalgebras, see \cite{BKLS}. Some simple Lie algebras over the ground field of characteristic $p>0$ (and some simple Lie superalgebras over any field) have no non-degenerate invariant bilinear form induced by the (super)-trace in any irreducible representation.

\ssec{Double extensions of NIS-algebras} From ancient time people computed nontrivial central extensions and outer derivations of Lie (super)algebras separately, see a recent review \cite{BGLL1}. There were known, however, important examples when both an outer derivation and a nontrivial central extension of a given Lie (super)algebra are present simultaneously; interestingly, this happens often in presence of a non-degenerate invariant (super)symmetric bilinear form.

Medina and Revoy in \cite{MR} introduced the notion of \textit{double extensions} for NIS-Lie algebras in characteristic 0, and showed that any such algebra $\fg$ can be described inductively in terms of another NIS-Lie algebra $\fa$ of dimension $\dim(\fg)-2$, provided the center of $\fg$ is not trivial (see also \cite{FS}). Solvable NIS-Lie algebras can be embraced by this method since the center is not trivial. 

More interesting examples, however, are Manin triples,  Manin--Olshanski triples, and affine Kac--Moody Lie algebras $\fg(A)$ with Cartan matrix $A$. Each of these examples is a double extension. For example, $\mathfrak{g}(A)$ is a double extension\footnote{The central extension in this case is defined by a 2-cocycle not of the form $B_\fa(D(a),b)$ as described by Eq. (\ref{Def}).} of the loop algebra $\fa=\mathfrak{k}\otimes \Cee [t^{-1},t]$ that has no Cartan matrix. The basis of $\fa$ lacks 2 elements of $\fg(A)$: the central element $x$ and the outer derivative $D=t\frac{d}{dt}$ of $\fa$, where $t=e^{i\varphi}$, and $\varphi$ is the angle parameter on the circle. 

Recall the definition from \cite{MR}. Let $\fa$ be a Lie algebra over a ground field $\Kee$ of characteristic $0$, $B_\fa$  a non-degenerate invariant symmetric bilinear form on $\fa$ and $D\in\fder(\fa)$ an \textit{outer} derivation\footnote{If the derivation $D$ is inner, then the double extension $\fg$ is nothing but $\fa\oplus \fc$, where $\fc$ is a 2-dimensional center of $\fg$.}. If $B_\fa$ is $D$-invariant, i.e., 
\[
B_\fa(D(a),b)+B_\fh(a,D(b))=0\text{~~ for any $a,b\in \fa$, }
\]
then there exists a NIS-Lie algebra structure on $\fg:=\mathscr{K} \oplus \fa\oplus \mathscr{K}^*$, hereafter called the \textit{double extension} of $\fa$, where $\mathscr{K}=\Kee x$, and $\mathscr{K}^*=\Kee x^*$ for $x^*=D$, defined as follows (for any $a,b\in\fa$):
\begin{equation}\label{Def}
[x,\fg]=0; \quad [a,b]_\fg=[a,b]_\fa+B_\fa (D(a),b)x; \quad [x^*,a]=D(a).
\end{equation}

Clearly, since $x$ spans an ideal of $\fg$, the dual $x^*$ of $x$ spans an outer derivation of $\fa$; therefore it might be reasonable to abbreviate the term \textit{double extension} and call it more suggestively \textit{D-extension} hinting at the presence of both (central) extension and a derivation ($D$, often given explicitly).

Examples of Lie (super)algebras that naturally have several $k\geq1$ linearly independent central elements $x_i$, where $i=1,\dots, k$,  and corresponding to them outer derivations $x_i^*$, are related to Lie (super)algebras with Cartan matrices of corank $k$, see \cite{BGL2}.

\ssec{Superization} The first superization of results by Medina and Revoy is due to Benamor and Benayadi \cite{BB} provided $\mathrm{dim}(\fg_\od)=2$. They showed that every non-simple NIS-Lie superalgebra can be described by successive $D$-extensions. Further, Benayadi \cite{B} extended the result provided  $\fg_\ev$ is reductive or solvable and $\fg_\od$ is a completely reducible $\fg_\ev$-module. Albuquerque, Barreiro and Benayadi \cite{ABB} generalized the result of \cite{B} to the case where $\fg_\od$ is not completely reducible $\fg_\ev$-module. The inductive description in this case, however, is given by means of another type of double extensions (for details, see \cite{Bor} \cite{B2} and \cite{BBB}). Later on, the same authors introduced the notion of ``odd" NIS-Lie superalgebras, and showed that, in this case, $\fg_\od$ is a completely reducible $\fg_\ev$-module if and only if $\fg_\ev$ is reductive.

Superization of the definition of $D$-extension introduces completely new possibilities: \textbf{the bilinear form $B_\fa$, the derivation $D$, the central element $x$, and/or its dual $x^*$ may be odd}\footnote{The parities are constrained: $p(B_\fa)=p(x)+p(x^*)$.};  the first of these possibilities is related with non-standard Manin-Olshansky triples, see \cite{LSh}, and loops with values in $\fp\fs\fq(n)$ and $\fh^{(1)}(0|2n+1)$. Passage to the ground field of positive characteristic $>2$ brings new examples but nothing conceptually new when it comes to the construction of $D$-extensions, except, perhaps, for inductive descriptions \`a la Medina and Revoy, because Lie's theorem and the Levi decomposition do not hold true anymore in positive characteristic.

The \textbf{case of characteristic 2} is completely different and requires new definitions and methods. Here we investigate all these new features  boldfaced above. The methods described in \cite{B}, \cite{BB} and \cite{ABB} can be applied \textit{mutatis mutentis} to the case of Lie superalgebras in positive characteristic $p>2$ without any difficulties. The only possible difference is that the space of derivations in positive characteristic could be larger than that in the case of characteristic zero, and hence new, up to an isometry, $D$-extended Lie superalgebras might  appear. In characteristic 2, contrariwise, the methods of \cite{B}, \cite{BB} and \cite{ABB} fail completely, mainly because of the squaring on a Lie superalgebra (see \S \ref{SecDef}). Moreover, even though the double extension \`a la Medina and Revoy could be applied to the Lie algebra $F(\fa)$, where $F$ is the functor of forgetting the parity, it will not help to recover the squaring on $\mathscr{K} \oplus F(\fa)\oplus \mathscr{K}^*$ to turn the later into a Lie superalgebra. Hence, we need another ingredient, a peculiarity of characteristic 2: a {\it quadratic form} defined on $\fa_\od$. To summarize: carrying out a $D$-extension in characteristic 2, requires an outer derivation $D\in \fder(\fa)$, a $D$-invariant NIS $B_\fa$ on $\fa$ and a quadratic form $\alpha$ on $\fa_\od$.

\ssec{Semi-trivial cocycles} In \cite{BGLLS}, the authors introduced the notion of {\it semi-trivial} 2-cocycles in the deformation theory of Lie (super)algebras. A 2-cocycle is called semi trivial if it is cohomology class is non-trivial in $\mathrm{H}^2(\fa;\fa$) but the deform in the direction of this 2-cocycle is isomorphic as abstract Lie (super)algebra (not as deformed ones) to the original one.

In this paper we encounter a similar phenomenon: two non-cohomolgous derivations $D$ and $\tilde D$ might produce the same $D$-extension up to an isometry. The description of these derivations, called {\it semi-trivial}, is given in \S \ref{ST1} and \S\ref{ST2}. This phenomenon has been observed first for Lie algebras in characteristic 0 in \cite{FS}. We give a sufficient and necessary condition for two $D$-extensions $(\fg,B_\fg)$ and $(\tilde \fg,B_{\tilde \fg})$ to be isometric.

The Chevalley-Eilenberg differential for Lie superalgebras in characteristic 2 was modified to take the squaring into account, see \cite{BGLL1}. With this definition, we do have $\mathfrak{out}(\fa)=\mathrm{H}^1(\fa;\fa)$, as shown in \cite{BGLL1}.

\ssec{The NIS-Lie superalgebras $\widetilde{\mathfrak{po}}(0|4)$ and $\mathfrak{po}(0|5;m)$} The first derived of the Hamiltonian Lie superalgebra $\fh^{(1)}(0|4)$ is a NIS-Lie superalgebra. Unexpectedly, we show that $\fh^{(1)}(0|4)$ admits {\bf three} non-isometric $D$-extensions: the Poisson Lie superalgebra $\mathfrak{po}(0|4)$, the general Lie superalgebra $\mathfrak{gl}(2|2)$ endowed with the {\it super trace}, and a new one we call $\widetilde{\mathfrak{po}}(0|4)$. This is also a peculiarity of the characteristic 2, since the NIS-Lie superalgebra $\widetilde{\mathfrak{po}}(0|4)$ has no analog in characteristic $p\not=2$ (see \S \ref{Exa}). 

We also construct a parametric family of Poisson NIS-Lie superalgebras $\mathfrak{po}(0|5;\lambda)$ as a result of a $D$-extension of $\fh^{(1)}(0|5)$. We show that, for $\lambda\not=0$, no member of this family is isometric to $\mathfrak{po}(0|5)$ but coincides with it if $\lambda=0$.

Using Grozman's package {\it SuperLie} \cite{Gr}, we computed for $m=6,7$ and $8$ the cohomology $\mathrm{H}^1(\fh^{(1)}(0|m);\fh^{(1)}(0|m))$. The analogues of the derivation that led to $\widetilde{\mathfrak{po}}(0|4)$ seem also to exist for $m=6$ and $8$. Also, analogues of the derivation that led to $\mathfrak{po}(0|5;\lambda)$ seem also to exist for $m=7$.  Therefore, we propose the following conjecture.

{\bf Conjecture.} The NIS-Lie superalgebras $\widetilde{\mathfrak{po}}(0|2m)$ and $\mathfrak{po}(0|2m+1;\lambda)$ exist for every $m\geq2$. 

\ssec{Examples of NIS-algebras} There are several classes of simple Lie (super)algebras on which the existence of the non-degenerate invariant bilinear form is known (for details, see \cite{BKLS}):

$\bullet$ Lie (super)algebras with indecomposable symmetrizable invertible Cartan matrix $A$. For a recipe how to construct (a unique, up to proportionality and non-degenerate if $A$ is invertible) invariant symmetric bilinear form on such Lie algebra, see \cite{K}, Ch.2; superization is immediate if $p\neq 2$. Such Lie (super)algebras are classified in the following cases:

A) Finite dimensional Lie (super)algebras over an algebraically closed field of prime characteristic; see \cite{BGL2}. Sometimes it is possible to extend NIS from these (and several other types of) Lie (super)algebras to their deforms; for a classification of possible cases, see \cite{BKLS}.

B)  Lie (super)algebras with indecomposable Cartan matrices are classified over $\Cee$ if 

(Ba) their growth is finite or polynomial, see \cite{K2}, \cite{Sc}, \cite{Se} and \cite{V},

(Bb) they are \lq\lq almost affine", see \cite{CCLL}.

$\bullet$ Lie superalgebras of Calogero--Moser model, see \cite{KS} and \cite{KT}. It is proven that simple Lie (super)algebra might have several non-degenerate invariant bilinear  forms; this is impossible for finite dimensional Lie algebras.

$\bullet$ There are also other simple Lie superalgebras which are not subquotients of Lie superalgebras with Cartan matrix but still having a NIS, in all the cases but one, this NIS is ODD, see \cite{BKLS}.



\ssec{Notation} The statements proved with the aid of \textit{SuperLie} code are called \textbf{Claims}. To make our text more graphic, we say  \lq\lq $D_\ev$-extensions" or \lq\lq $D_\od$-extensions", depending on the parity of $D$, rather than just \lq\lq $D$-extensions".

The 1-cochain $\hat x\in C^1(\fg)$ denotes the dual of $x\in\fg$; when exterior products of cochains are considered we assume that $p(\hat x)=p(x)+\od$, where $p$ is the parity function.
\section{Lie superalgebras  for $p=2$} \label{SecDef}
For basics on Lie superalgebras over fields $\Kee$ of characteristic 2, see \cite{LeD, LeD2, BGLLS2}. For Lie superalgebras with Cartan matrices, see \cite{BGLL, BGL2}; for descriptions in terms of Cartan-Tanaka-Shchepochkina prolongations, see \cite{BGLLS1, BGLLS2}; for the classification of simple Lie superalgebras, see \cite{BGL2,BLLSq}. 

If $p\neq2$, superization of many notions of Linear Algebra is performed, as is now well-known, with the help of the Sign Rule. For $p=2$, one has to be more subtle.
We recall basic definitions retaining the minus sign from definitions for $p\neq 2$: for clarity. 

A \textit{Lie superalgebra} in characteristic 2 is a
superspace $\fg=\fg_\ev\oplus\fg_\od$ over a field $\Bbb K$ such that the even part
$\fg_\ev$ is a Lie algebra, the odd part $\fg_\od$ is a
$\fg_\ev$-module made two-sided by
\textit{anti}-symmetry, and on the odd part $\fg_\od$ a \textit{squaring} is
defined as a map given by
\begin{equation}\label{squaring}
\begin{array}{c}
s_\fg:\fg_\od \rightarrow \fg_\ev \quad f\mapsto s_\fg(f)\quad \text{such that $s_\fg(\lambda f)=\lambda^2 s_\fg(f)$ for any $f\in
\fg_\od$ and $\lambda \in \Kee$, and}\\
\text{the map } \fg_\od\times \fg_\od \rightarrow \fg_\ev \text{ given by } (f,g)\mapsto       s_\fg(f+g)-s_\fg(f)-s_\fg(g)\\
\text{~is a bilinear map on $\fg_\od$ with values
in $\fg_\ev$.}
\end{array}
\end{equation}
The bracket on $\fg_\ev$ as well as the action of $\fg_\ev$ on $\fg_\od$ is denoted by the same symbol $[\cdot ,\cdot]_\fg$. For any $f, g\in\fg_\od$, their bracket is
\[
[f,g]_\fg:=  s_\fg(f+g)-s_\fg(f)-s_\fg(g).
\]
The Jacobi identity involving the squaring reads as follows:
\begin{equation}\label{JIS}
[s_\fg(f),g]_\fg=[f,[f,g]_\fg]_\fg\;\text{ for any $f\in \fg_\od$ and $g\in\fg$}.
\end{equation}
It is worth noticing that given such a data we get: 
\begin{eqnarray}
\label{SkS}
[f,f]_\fg&=& 0 \text{ for any $f \in \fg$}; \text{ and }\\[2mm]
\label{IJ}
[f,[g,h]_\fg]_\fg + \circlearrowleft (f,g,h) &=& 0\text{ for any $f, g, h \in \fg$}.
\end{eqnarray}
By forgetting the squaring and keeping only the brackets by setting $[f,f]=0$ for $f$ odd, $\fg$ turns into an ordinary Lie algebra, see \cite{BLLSq}.

For any Lie superalgebra $\fg$ in characteristic 2, its \textit{derived algebras} are
defined to be (for $i\geq 0$)
\[
\fg^{(0)}: =\fg, \quad
\fg^{(i+1)}=[\fg^{(i)},\fg^{(i)}]_\fg+\Span\{s_\fg(f)\mid f\in (\fg^{(i)})_\od\}.
\]
A linear map $D:\fg\rightarrow \fg$ is called a \textit{derivation} of the Lie superalgebra $\fg$ if, in addition to
\begin{eqnarray}
\label{Der1} D([f,g]_\fg)&=&[D(f),g]_\fg+[f,D(g)]_\fg\quad \text{for any $f\in \fg_\ev$ and $g\in \fg$ we have}\\[2mm]
\label{Der2} D(s_\fg(f))&=&[D(f),f]_\fg\quad \text{for any $f\in \fg_\od$}.
\end{eqnarray}
It is worth noticing that condition (\ref{Der2}) implies condition (\ref{Der1}) if $f,g\in \fg_\od$. 

We denote the space of all derivations of $\fg$
by $\fder (\fg)$. 

Let $(\fg,[\cdot ,\cdot]_\fg,s_\fg)$ and $(\fh,[\cdot ,\cdot]_\fh,s_\fh)$ be two Lie superalgebras in characteristic 2. An even linear map $\varphi:\fg\rightarrow \fh$ is called a \textit{morphism} (of Lie superalgebras) if, in addition to
\begin{eqnarray*}
\label{Hom1} \varphi([f,g]_\fg)&=&[\varphi(f), \varphi (g)]_\fh\quad \text{for any $f\in \fg_\ev$ and $g\in \fg$ we have}\\[2mm]
\label{Hom2} \varphi(s_\fg(f))&=& s_\fh(\varphi(f))\quad \text{for any $f\in \fg_\od$}.
\end{eqnarray*}

Therefore, morphisms in the category of Lie superalgebras in characteristic 2 preserve not only the bracket but the squaring as well. In particular, subalgebras and ideals have to be stable under the bracket and the squaring. 

An even linear map $\rho: \fg\tto\fgl(V)$ is a \textit{representation
of the Lie superalgebra} $\fg$ in the superspace $V$ called \textit{$\fg$-module} if
\begin{equation}\label{repres}
\begin{array}{l}
\rho([f, g])=[\rho (f), \rho(g)]\quad \text{ for any $f, g\in
\fg$}; \text{ and }
\rho (s_\fg(f))=(\rho (f))^2\text{~for any $f\in\fg_\od$.}
\end{array}
\end{equation}
Let $B_\fg$ be a bilinear form on $\fg:=(\fg, [\cdot ,\cdot]_\fg,s_\fg)$. We say that\begin{enumerate}
\item [(A)] \label{Q1} $B_\fg$ is symmetric if $B_\fg(f,g)=B_\fg(g,f)$ for any $f,g\in\fg$;  and $B_\fg(f,f)=0$ \text{ for any $f\in \fg_\od$; } 
\item [(B)] \label{Q2} $B_\fg$ is invariant if $B_\fg([f,g]_\fg,h)=B_\fg(f,[g,h]_\fg)$ for any $f,g,h\in \fg$.
\end{enumerate}

We call the Lie superalgebra $\fg:=(\fg, [\cdot ,\cdot]_\fg,s_\fg)$ a \textit{NIS-Lie superalgebra} if it admits a homogenous non-degenerate, invariant and symmetric bilinear form $B_\fg$. We denote such a superalgebra by $(\fg, B_\fg)$.

A NIS-Lie superalgebra $(\fg, B_\fg)$is said to be \textit{reducible} if it can be decomposed into direct sums of ideals, namely $\fg=\oplus I_i$, such that all $I_i$ are mutually orthogonal.

\ssec{Manin triples}
Let $(\fh, s_\fh, [\cdot ,\cdot]_\fh)$ be a finite-dimensional Lie superalgebra (not necessarily \lq\lq NIS"), and let the dual space have the structure of an abelian Lie superalgebra.  A NIS-Lie superalgebra structure on $\fg:=\fh\oplus \fh^*$ is naturally defined as follows. Over $\Cee$, the triple $(\fg, \fh, \fh^*)$ is called a \textit{Manin triple} (see \cite{Dr}), although our construction is a particular case of Manin triple. The squaring on $\fg$ is defined as follows (for any $h\in \fh_\od$ and $\pi\in \fh^*_\od $):
\begin{equation}
\label{squaringh}
s_{\fg}(h+\pi):=s_\fh(h)+\pi\circ \text{ad}_h \quad \text{ for any $h+\pi \in \fg_\od$},
\end{equation}
where $\pi\circ \ad_h\in \fh^*_\ev$ should be understood as follows: 
\[
\pi\circ \ad_h(k)=\pi([h,k]_\fh)\; \text {for all $k\in \fh$}.
\]

The bracket of two elements
is defined as follows:
\begin{equation}
\label{bracketh}
[h+\pi, h'+\pi']_{\fg}:=[h,h']_\fh+\pi'\circ \text{ad}_{h}+\pi\circ \text{ad}_{h'}\quad \text{for any $h+\pi,  \ h'+\pi'\in\fg$.}
\end{equation}
It is easy to show that the map $s_{\fg}$ defined by Eq. (\ref{squaringh}) is indeed a squaring, i.e., satisfies Eq. (\ref{squaring}), and the bracket $[\cdot , \cdot]_{\fg}$ defined by Eq. (\ref{bracketh}) satisfies the Jacobi identity.

\underline{$\bullet$ $B_{\fg}$ is even}. We define a bilinear form on $\fg$ as follows:
\begin{equation}\label{Mform}
B_{\fg}(h+\pi, h'+\pi'):=\pi(h')+\pi'(h) \quad \text{for any $h+\pi, \ h'+\pi'\in\fg$.}
\end{equation}
Obviously, $B_\fg$ is even. It is easy to show that the bilinear form $B_\fg$ is non-degenerate, invariant and symmetric. 

\underline{$\bullet$ $B_{\fg}$ is odd}.  A NIS-Lie superalgebra structure on $\fg:=\fh\oplus \Pi(\fh^*)$, where $\Pi$ is the change of parity functor, is naturally defined as follows. Define the squaring on $\fg$ as follows (for any $h\in \fh_\od$ and $\pi\in \fh^*_\ev $):
\begin{equation}
\label{squaringhp}
s_{\fg}(h+\Pi(\pi)):=s_\fh(h)+\Pi(\pi)\circ \text{ad}_h \quad \text{ for any $h+\Pi(\pi) \in \fg_\od$}.
\end{equation}
The bracket of two elements
is defined as follows:
\begin{equation}
\label{brackethp}
[h+\Pi(\pi), h'+\Pi(\pi')]_{\fg}:=[h,h']_\fh+\Pi(\pi')\circ \text{ad}_{h}+\Pi(\pi)\circ \text{ad}_{h'}\quad \text{for any $h+\pi,  \ h'+\pi'\in\fg$.}
\end{equation}
It is easy to show that the map $s_{\fg}$ defined by Eq. (\ref{squaringhp}) is indeed a squaring, and the bracket $[\cdot , \cdot]_{\fg}$ defined by Eq. (\ref{brackethp}) satisfies the Jacobi identity.

We define a bilinear form on $\fg$ as follows:
\begin{equation}\label{OMform}
B_{\fg}(h+\Pi(\pi), h'+\Pi(\pi')):=\pi(h')+\pi'(h) \quad \text{for any $h+\Pi(\pi), \ h'+\Pi(\pi')\in\fg$.}
\end{equation}
This bilinear form is odd, non-degenerate, invariant, and symmetric.
 
\ssec{Quadratic and bilinear forms in characteristic 2} \label{Arf} A given map $\alpha:V\rightarrow \Kee$, where $V$ is a $\Kee$-vector space, is called a \textit{quadratic form} if 
\[
\begin{array}{c}
\alpha(\lambda v)=\lambda^2 \alpha(v)\text{ for any $\lambda \in \Kee$ and for any $v \in V$, and the map }\\
(u,v) \mapsto B_\alpha(u,v):=\alpha(u+v)-\alpha(u)-\alpha(v) \text{ is bilinear.}
\end{array}
\] 
The form $B_\alpha$ is called the \textit{polar form} of $\alpha$. Recall that non-degenerate quadratic forms over a field of characteristic 2 are classified by the 
\textit{Arf} invariant, see \cite{D}. Recently, Lebedev classified non-degenerate bilinear forms over a perfect field $\Kee$ (i.e., such that $\Kee^2=\Kee$), see \cite{LeD, LeD2}. This is a non-trivial result not related with a well-known classification of quadratic forms in any characteristic, because in characteristic 2, the maps 
\[
\alpha \longleftrightarrow B_\alpha
\]
are not one-to-one. On any $2n$-dimensional space, there exists a Darboux basis $(v_i)_{i=1}^{2n}$ in which every non-degenerate quadratic form can be written as:
\[
\alpha \left (\sum_{1\leq i \leq 2n}\lambda_i v_i \right )=\sum_{1\leq i \leq n} \lambda_i \lambda_{i+n}+A(\lambda_n^2+\lambda_{2n}^2), \quad \text{where } A\in \Kee.
\]
Its Arf invariant is $A^2$ (cf. \cite{D}). The polar form associated with it is given by the formula
\[
B_\alpha(x,y)=\sum_{1\leq i \leq n} \left ( \lambda_i \mu_{n+i}-\mu_i\lambda_{n+i}\right ),  \text{ where $x=\sum_{1\leq i\leq 2n}\lambda_i v_i$ and $y=\sum_{1\leq i\leq 2n}\mu_i v_i$.}
\]
\section{The case where the bilinear form $B$ is even}
 
 \subsection{$D_\ev$-extensions}
\sssbegin{Theorem}\label{MainTh} Let $(\fa, B_\fa)$ be a NIS-Lie superalgebra in characteristic $2$ such that $B_\fa$ is even. Let $D\in \fder_\ev(\fa)$ be a derivation satisfying the following conditions:
\begin{eqnarray}
\label{D1} B_\fa(D(a),b)+B_\fa(a,D(b))=0\text{ for any } a,b \in \fa; \; B_\fa(D(a), a)=0\text{ for any } a \in \fa_\ev.
\end{eqnarray}
Let $\alpha: \fa_\od \rightarrow \Bbb K$ be a quadratic form and $B_\alpha$ its polar form which satisfies 
\begin{eqnarray}
\label{D3}
B_\fa(a, D(b))&=&B_\alpha(a,b) \text{ for any  $a,b \in \fa_\od$}. 
\end{eqnarray}

Then there exists a NIS-Lie superalgebra structure on $\fg:=\mathscr{K} \oplus \fa \oplus \mathscr{K} ^*$, where $\mathscr{K} :=\Span\{x\}$ for $x$ even, defined as follows. The squaring is given by 
\[
s_\fg(a):= s_\fa(a)+\alpha (a) x \qquad \text{ for any } a\in \fg_\od \; (=\fa_\od).
\]
The bracket on $\fg$ is defined as follows:
\be\label{*}
[\mathscr{K} ,\fg]_\fg=0, \qquad [a,b]_\fg:=B_\fa(D(a),b)x+[a,b]_\fa, \qquad [x^*,a]_\fg:=D(a)\text{ for any } a,b\in  \fa.
\ee
The non-degenerate symmetric bilinear form $B_\fg$ on $\fg$ is defined as follows: 
\begin{eqnarray*}
&{B_\fg}\vert_{\fa \times \fa}:= B_\fa, \quad B_\fg(\fa,\mathscr{K} ):=0, \quad B_\fg(x,x^*):=1, \quad B_\fg(\fa,\mathscr{K} ^*):=0, \\ 
&B_\fg(x,x):=0,\ B_\fg(x^*,x^*) \; \text{ arbitrary}.
\end{eqnarray*}
Moreover, the form $B_\fg$, obviously even, is invariant on $\fg$. 
\end{Theorem}

We call the Lie superalgebra $(\fg, B_\fg)$ constructed in Theorem \ref{MainTh} a \textit{$D_\ev$-extension} of $(\fa, B_{\fa})$  by means of $D$ and $\alpha$.
\begin{proof}
Let us first show that $s_\fg$ is indeed a squaring on $\fg$. Recall that since $\mathscr{K} $ and $\mathscr{K}^*$ are even vector spaces, then $\fg_\od=\fa_\od$. Now, let $\lambda \in \Bbb K$ and let $a\in\fa_\od$; we have
\[
s_\fg(\lambda a)=\alpha (\lambda a) x+ s_\fa(\lambda a)= \lambda^2 \alpha ( a) x+ \lambda^2 s_\fa(a)=\lambda^2s_\fg(a),
\]
since $\alpha$ is a quadratic form and $s_\fa$ is a squaring on $\fa$. 
%
Let us show that the bracket on $\fg$ is symmetric. Indeed, using condition (\ref{D1}), for any $f=tx+a+sx^*$ and $g=\tilde tx+b+\tilde sx^*$ in $\fg$, where $s,t,\tilde s, \tilde t \in \Kee$, we have
\begin{eqnarray*}
[f, g]_\fg&=&B_\fa(D(a),b)x+[a,b]_\fa+s D(b)+\tilde s D(a)\\
&=& B_\fa(D(b),a)x+[b,a]_\fa+s D(b)+\tilde s D(a)\\
&=& [g,f]_\fg.
\end{eqnarray*}

Let us check the Jacobi identity relative to the squaring. Indeed, for any $a\in \fg_\od$ and for any even element $g=sx+c+tx^*$, where $c\in \fa$, we have
\[
\begin{array}{lcl}
[s_\fg(a),g]_\fg+[a,[a,g]_\fg]_\fg&=&[\alpha(a)x+s_\fa(a),g]_\fg +[a,[a,c]_\fa+tD(a)]_\fg \\[2mm]
&=&[s_\fa(a),g]_\fg + B_\fa(D(a),[a,c]_\fa )x+[a,[a,c]_\fa]_\fa +t[a,D(a)]_\fa \\[2mm]
&&+tB_\fa(D(a),D(a))x\\[2mm]
&=&t D(s_\fa(a))+B_\fa(D(s_\fa(a)),c)x+[s_\fa(a),c]_\fa\\[2mm]
&&+ B_\fa(D(a),[a,c]_\fa )x+[a,[a,c]_\fa]_\fa +t[a,D(a)]_\fa \\[2mm]
&=&[s_\fa(a),c]_\fa+[a,[a,c]_\fa]_\fa =0,
\end{array}
\]
since $s_\fa$ is a squaring on $\fa$, $D\in\fder(\fa)$, condition (\ref{D1}), and the fact that 
\[
B_\fa(D(s_\fa(a)),c)=B_\fa([D(a),a]_\fa,c)=B_\fa(D(a),[a,c]_\fa).
\]
The proof of the Jacobi identity for the Lie bracket
would also follow from the fact that the Medina-Revoy construction
is still valid in characteristic 2, and for the convenience of
the reader we repeat the arguments for that case. 

To check the Jacobi identity, we proceed as follows. 

If $h=x$, the identity
\[
[x,[f,g]_\fg]_\fg+[g,[x,f]_\fg]_\fg+[f,[g,x]_\fg]_\fg=0
\]
is certainly satisfied since $x$ is central in $\fg$. 

If $h=x^*$, the identity 
\[
[x^*,[f,g]_\fg]_\fg+[g,[x^*,f]_\fg]_\fg+[f,[g,x^*]_\fg]_\fg=0
\]
is also satisfied for the following reasons. If either $f$ or $g$ is $x$, then we are done since $x$ is central in $\fg$. Now if $f=x^*$ (or the other way round, $g=x^*$), then
\[
[x^*,[x^*,g]_\fg]_\fg+[g,[x^*,x^*]_\fg]_\fg+[x^*,[g,x^*]_\fg]_\fg=2D(D(g))=0.
\]
Let us assume now that $f, g \in\fa.$ We have (if $f$ and $g$ are not both odd)
\[
\begin{array}{lcl}
[x^*,[f,g]_\fg]_\fg+[g,[x^*,f]_\fg]_\fg+[f,[g,x^*]_\fg]_\fg&=&[x^*, B_\fa(D(f),g)x+[f,g]_\fa]_\fg+ [g, D(f)]_\fg\\[2mm]
&&+ [f,D(g)]_\fg\\[2mm]
&=&D([f,g]_\fa)+B_\fa(D(g),D(f))x+ [g, D(f)]_\fa\\[2mm]
&&+B_\fa(D(f),D(g))x+ [f,D(g)]_\fa\\[2mm]
&=&0,
\end{array}
\]
since $D\in\fder(\fa)$ and $B_\fa$ is (super)symmetric. 

If $f$ and $g$ are both odd, then 
\[
\begin{array}{lcl}
[x^*,[f,g]_\fg]_\fg+[g,[x^*,f]_\fg]_\fg+[f,[g,x^*]_\fg]_\fg&=&[x^*, s_\fg(f+g)+s_\fg(f)+s_\fg(g)]_\fg+ [g, D(f)]_\fg\\[2mm]
&&+ [f,D(g)]_\fg\\[2mm]
&=&[x^*, B_\alpha(f,g)x+[f,g]_\fa]_\fg+B_\alpha(g,D(f))x\\[2mm]
&&+[g,D(f)]_\fa+ B_\alpha(f,D(g))x+[f,D(g)]_\fa=0.\\[2mm]
\end{array}
\]
since $D\in\fder(\fa)$, and $B_\alpha(f,D(g))=B_\alpha(D(f),g)$ from condition (\ref{D3}). 

From now and on we assume that $f,g, h \in \fa$. We distinguish several cases to check the Jacobi identity.

If $f, g$ and $h$ are even, then we have
\[
\begin{array}{lcl}
[f,[g,h]_\fg]_\fg+\circlearrowleft (f,g,h)&=&[f,B_\fa(D(g),h)x+[g,h]_\fa]_\fg+\circlearrowleft (f,g,h)\\[2mm]
&=&[f,[g,h]_\fa]_\fa+B_\fa(D(f),[g,h]_\fa)x +\circlearrowleft (f,g,h)=0,\\[2mm]
\end{array}
\]
because the JI holds on $\fa$ and thanks to Lemma \ref{lemma1}.

\sssbegin{Lemma}
\label{lemma1}
Let $(\fa,B_\fa)$ be a NIS-Lie superalgebra. Let $D$ be in $\fder(\fa)$. Suppose further that $B_\fa$ satisfies condition (\ref{D1}). Then
\[
B_\fa(D(h),[f,g]_\fa)+ \circlearrowleft (f,g,h)=0 \text{ for any $f,g,h\in\fa$}.
\]
\end{Lemma}

\begin{proof} The proof follows from the fact that $B_\fa$ is invariant and $D$ is a derivation. Indeed, 
\[
\begin{array}{lcl}
B_\fa(D(h),[f,g]_\fa)+\ \circlearrowleft (f,g,h) &=&B_\fa(h,D[f,g]_\fa)+B_\fa(h,[D(f),g]_\fa)+ B_\fa(h,[f,D(g)]_\fa)\\[2mm]
&=&B_\fa(h,D[f,g]_\fa+[D(f),g]_\fa+[f,D(g)]_\fa)=0.
\end{array}
\]\end{proof}

If $f$ and $g$ are both odd but $h$ is even, then we have
\[
\begin{array}{lcl}
[f,[g,h]_\fg]_\fg +[h,[f,g]_\fg]_\fg+[g,[h,f]_\fg]_\fg&=&[f,B_\fa(D(g),h)x+[g,h]_\fa]_\fg+[h,B_\alpha(f,g)x+[f,g]_\fa]_\fg \\[2mm]
&& +[g,B_\fa(D(h),f)x+[h,f]_\fa]_\fg\\[2mm]
&=&B_\alpha(f,[g,h]_\fa)x+B_\alpha(g,[h,f]_\fa)x+B_\fa(D(h),[f,g]_\fa)x\\[2mm]
&&+[h,[f,g]_\fa]_\fa+[f,[g,h]_\fa]_\fa+[g,[h,f]_\fa]_\fa=0,\\
\end{array}
\]
because the JI holds on $\fa$, and $B_\alpha(f,[g,h]_\fa)+B_\alpha(g,[h,f]_\fa)+B_\fa(D(h),[f,g]_\fa)=0$. Indeed, using conditions (\ref{D1}) and (\ref{D3}), we have
\[
\begin{array}{lcl}
B_\fa(D(h),[f,g]_\fa)&=& B_\fa(h, D[f,g]_\fa)\\[2mm]
&=&B_\fa(h, [D(f),g]_\fa+[f,D(g)]_\fa)\\[2mm]
&=&B_\fa(D(f),[h,g]_\fa)+B_\fa(D(g),[h,f]_\fa)\\[2mm]
&=&B_\alpha(f,[h,g]_\fa)+B_\alpha(g,[h,f]_\fa).\\[2mm]
\end{array}
\]
If $f$ is odd and both $g$ and $h$ are even, we have
\[
\begin{array}{lcl}
[f,[g,h]_\fg]_\fg+ \circlearrowleft (f,g,h)&=&[f,B_\fa(D(g),h)x+[g,h]_\fa]_\fg+ \circlearrowleft (f,g,h)\\[2mm]
&=&[f,[g,h]_\fa]_\fa+ B_\fa(D(f),[g,h]_\fa)x+\circlearrowleft (f,g,h)=0,
\end{array}
\]
since the JI holds on $\fa$ and $B_\fa$ is even.
 
Let us now show that $B_\fg$ is invariant. We should check that
\[
B_\fg([f,g]_\fg, h)=B_\fg (f, [g,h]_\fg) \quad \text{for any $f,g$ and $h$  $\in \fg$.}
\]
This is true if $f=x$ because $x$ is central and by the very definition of $B_\fg$.  If $f=x^*$, $g=ux+b+vx^*$ and $h=wx+c+zx^*$ (where $g$ and $h$ are not both odd, $u,v,wz\in \Kee$), we have
\[
\begin{array}{lcl}
B_\fg([f,g]_\fg, h)&=&B_\fg(D(b), wx+c+zx^*)\\[2mm]
&=& B_\fa(D(b), c).\\[2mm]
\end{array}
\]
On the other hand, 
\[
\begin{array}{lcl}
B_\fg(f,[g, h]_\fg)&=&B_\fg(x^*, B_\fa(D(b),c)x+[b,c]_\fa+v D(c)+zD(b))\\[2mm]
&=&B_\fa(D(b), c) .\\[2mm]
\end{array}
\]
If $f=x^*$, $g=b$ and $h=c$, where $g, h \in \fa_\od$, we have
\[
\begin{array}{lcl}
B_\fg([f,g]_\fg, h)&=&B_\fa(D(b), c).
\end{array}
\]
On the other hand, using condition (\ref{D3}) we have
\[
\begin{array}{lcl}
B_\fg(f,[g, h]_\fg)&=&B_\fg(x^*, B_\alpha(b,c)x+[b,c]_\fa)\\[2mm]
&=&B_\alpha(b, c) \\[2mm]
&=& B_\fa(D(b),c).
\end{array}
\]
If $f, g, h\in\fa$, the invariance of $B_\fg$ easily follows from that of $B_\fa$.
The fact that the bilinear form $B_\fg$ is non-degenerate, even and symmetric is clear. 
\end{proof}
Now, we need the following definition. The \textit{\lq\lq special center" of $\fg$ relative to $B_\fg$} is the following set
\begin{equation}\label{speCen}
\fz_s(\fg):=\fz(\fg)\cap s_\fg(\fg_\od)^{\perp}.
\end{equation}
Observe that $\fz_s(\fg)_\od=\fz(\fg)_\od$ and $\fz_s(\fg)_\ev=\fz(\fg)_\ev \cap s_\fg(\fg_\od)^{\perp}$. Moreover, $\fz_s(\fg)$ is not necessarily an ideal.

\sssbegin{Proposition}
\label{Rec1}
Let $(\fg,B_\fg)$ be an irreducible NIS-Lie superalgebra. Suppose that $\fz_s(\fg)_\ev\not=\{0\}.$ Then $(\fg,B_\fg)$ is obtained as an $D_\ev$-extension from a NIS-Lie superalgebra $(\fa,B_\fa)$. 
\end{Proposition}
\begin{proof}
Let $x$ be a non-zero element in $\fz_s(\fg)_\ev$. The subspace $\mathscr{K}:=\Span\{x\}$ is an ideal in $(\fg,B_\fg)$ because $x$ is central in $\fg$. Moreover, $\mathscr{K}^\perp$ is also an ideal in $(\fg,B_\fg)$. Indeed, let us show first that $(\mathscr{K}^\perp)_\od=\fg_\od$. This is true because $x$ is orthogonal to any odd element since $B_\fg$ is even. Now let $f\in (\mathscr{K}^\perp)_\od=\fg_\od$. It follows that $B_\fg(x,s_\fg(f))=0$ since $x\in \fz_s(\fg)$. Hence $s_\fg(f)\in \mathscr{K}^\perp$. On the other hand, if $f\in \mathscr{K}^\perp$ and $g\in \fg$, then $[f,g]_\fg \in \mathscr{K}^\perp$ because
\[
B_\fg(x,[f,g]_\fg)=B_\fg([x,f]_\fg,g)=B_\fg(0,g)=0.
\]
Since $\mathscr{K}$ is 1-dimensional, then either $\mathscr{K} \cap \mathscr{K}^\perp=\{0\}$ or $\mathscr{K}\cap \mathscr{K}^\perp=\mathscr{K}$. The first case is to be disregarded because otherwise $\fg=\mathscr{K}\oplus \mathscr{K}^\perp$ and the Lie superalgebra $\fg$ will not be irreducible. Hence, $\mathscr{K}\cap \mathscr{K}^\perp=\mathscr{K}$. It follows that $\mathscr{K}\subset \mathscr{K}^\perp$ and $\dim(\mathscr{K}^\perp)=\dim(\fg)-1$. Therefore, there exists $x^* \in \fg_\ev$ such that
\[
\fg=\mathscr{K}^\perp\oplus \mathscr{K}^*, \quad \text{ where  $\mathscr{K}^*:=\Span\{x^*\}$}.
\]
This $x^*$ can be normalized to have $B_\fg(x,x^*)=1$. Besides, $B_\fg(x,x)=0$ since $\mathscr{K}\cap \mathscr{K}^\perp=\mathscr{K}$.

Let us define $\fa:=(\mathscr{K} +\mathscr{K}^*)^\perp$. We then have a decomposition $\fg=\mathscr{K} \oplus \fa \oplus \mathscr{K}^*$. 

Let us define a bilinear form on $\fa$ by setting: 
\[
B_\fa={B_\fg}\vert_{\fa\times \fa}.
\]  
The form $B_\fa$ is non-degenerate on $\fa$. Indeed, suppose there exists an $a\in\fa$ such that 
\[
B_\fa(a,f)=0\quad \text{ for any $f\in \fa$}.
\] 
But $a$ is also orthogonal to $x$ and $x^*$. It follows that $B_\fg(a,f)=0$ for any $f\in\fg$. Hence, $a=0$, since $B_\fg$ is nondegenerate. 

Let us show now that there exists a NIS-Lie superalgebra structure on the vector space $\fa$ for which $\fg$ is its double extension. Let $a,b \in \fa$. The bracket $[a,b]_\fg$ belongs to $\mathscr{K} \oplus \fa$ because $\fa \subset \mathscr{K}^\perp=\mathscr{K}\oplus \fa$ and the latter is an ideal. It follows that 
\[
[a,b]_\fg=\phi(a,b)x+\psi(a,b), \quad [a,x^*]_\fg=f(a)x+D(a),
\]
where $\psi(a,b), D(a) \in \fa$, and $\phi(a,b), f(a) \in \Bbb K$. Let us show that $\phi(a,b)=B_\fg(a,D(b))$. Indeed, since $B_\fg([a,b]_\fg,x^*)=B_\fg(a,[b,x^*]_\fg)$, it follows that 
\[
B_\fg(\phi(a,b)x+\psi(a,b),x^*)=B_\fg(a,f(b)x+D(b)).
\] 
Since $\psi(a,b)$ is orthogonal to $x^*$, and $a$ is orthogonal to $x$, and $B_\fg(x,x^*)=1$, we get $\phi(a,b)=B_\fg (a,D(b))$. The antisymmetry of the bracket $[\cdot , \cdot]_\fg$ implies the antisymmery of the map $\phi(a,b)$ which, in turn, implies that 
\[
B_\fg(a,D(b))=B_\fg(D(a),b) \quad \text{ for any $a,b\in\fa$.}
\]
Similarly, 
\[
0=B_\fg([x^*,x^*]_\fg,a)=B_\fg(x^*,[x^*,a]_\fg)=B_\fg(x^*, f(a)x+D(a))=f(a) \quad \text{ for any $a\in\fa$.}
\]
This implies that $f(a)=0$ for any $a\in \fa$.

The map $\psi$ is a bilinear on $\fa$ because $[\cdot , \cdot]_\fg$ is bilinear, so for convenience let us re-denote $\psi$ by $[\cdot , \cdot]_\fa$. We still have to show that the Jacobi identity is satisfied for $[\cdot , \cdot]_\fa$.
Let us now describe a squaring on $\fa$. Since $\fa \subset \mathscr{K}^\perp$, then $s_\fg(a)\in (\mathscr{K}^\perp)_\ev=\mathscr{K} \oplus \fa_\ev$, for any $a \in \fa_\od$.  It follows then that
\[
s_\fg(a)=\alpha(a)x+s_\fa(a).
\]
Now, because $s_\fg$ is a squaring on $\fg$,  
it follows that $\alpha$ is a quadratic form on $\fa_\od$ and $s_\fa$ behaves as a squaring on $\fa$. It follows that if $a, b\in \fa_\od$, then
\[
\begin{array}{lcl}
[a,b]_\fg&=&s_\fg(a+b)+s_\fg(a)+s_\fg(b)\\[2mm]
&=&B_\alpha(a,b)x+[a,b]_\fa,
\end{array}
\]
where we have put $[a,b]_\fa:=s_\fa(a+b)+s_\fa(a)+s_\fa(b)$.

Let us show next that $D$ is a derivation on $\fa$. Indeed, the condition 
\[
[s_\fg(a),x^*]_\fg=[a,[a,x^*]_\fg]_\fg
\] 
implies that $D(s_\fa(a))=[a,D(a)]_\fa$ because $B_\fg(D(a),D(a))=0$ since $a$ is odd.  Moreover, for any $a,b$ and $c$ in $\fa$ we have
\[
\begin{array}{lcl}
0&=&[[a,b]_\fg,c]_\fg+ \circlearrowleft (a,b,c)\\[2mm]
&=&[[a,b]_\fa + B_\fa(D(a),b)x,c]_\fg + \circlearrowleft (a,b,c)\\[2mm]
&=&[[a,b]_\fa,c]_\fa+ B_\fg(D[a,b]_\fa,c)x+ \circlearrowleft (a,b,c) \end{array}
\]
This implies that the JI for the bracket $[.,.]_\fa$ is satisfied and 
\[
B_\fg(c,D[a,b]_\fa+[a,D(b)]_\fa+[D(a),b]_\fa)=0\quad \text{ for any $a, b,c\in\fa$}.
\] Since $B_\fg$ is also non-degenerate on $\fa$, it follows that $D([a,b]_\fa)=[D(a),b]_\fa+[a,D(b)]_\fa$.  

Let us assume that $a$ and $b$ are odd and $c$ is even
\[
\begin{array}{lcl}
0&=&[[a,b]_\fg,c]_\fg+ [[c,a]_\fg,b]_\fg+[[b,c]_\fg,a]_\fg\\[2mm]
&=&[[a,b]_\fa,c]_\fg+ [[c,a]_\fa,b]_\fg+[[b,c]_\fa,a]_\fg\\[2mm]
&=&[[a,b]_\fa,c]_\fa+B_\fa(D[a,b]_\fa,c)x+ [[c,a]_\fa,b]_\fa+B_\alpha([c,a]_\fa,b)x+[[b,c]_\fa,a]_\fa\\[2mm]
&&+B_\alpha([b,c]_\fa,a)x.
\end{array}
\]
It follows that the JI is satisfied for the bracket $[\cdot , \cdot]_\fa$, provided  
\begin{equation}
\label{cond}
B_\fa(D[a,b]_\fa,c)+B_\alpha([c,a]_\fa,b)+B_\alpha([b,c]_\fa,a)=0.
\end{equation}
We will show later that this condition follows form another condition. 

We can similarly show that the JI identity on $\fg$ for one element odd and two elements even implies the JI on $\fa$. 

Now, we should prove that the bilinear form $B_\fa:={B_\fg}\vert_{\fa \times \fa}$ is also invariant on $(\fa, [\cdot ,\cdot]_\fa, s_\fa, B_\fa)$. Indeed, for any $a,b,c\in\fa$ with at least two of them even, we have
\[
\begin{array}{lcl}
B_\fa(a,[b,c]_\fa)&=&B_\fg(a, [b,c]_\fg+B_\fa(D(b),c)x)\\[2mm]
&=&B_\fg(a, [b,c]_\fg)\\[2mm]
&=&B_\fg([a,b]_\fg, c)\\[2mm]
&=&B_\fg([a,b]_\fa+B_\fa(D(a),b)x, c)\\[2mm]
&=&B_\fg([a,b]_\fa, c)\\[2mm]
&=&B_\fa([a,b]_\fa, c).\\
\end{array}
\]
Now, for any $a,b\in\fa_\od$ and $c\in\fa_\ev$, we have
\[
\begin{array}{lcl}
B_\fa(a,[b,c]_\fa)&=&B_\fg(a, [b,c]_\fg+B_\fa(D(b),c)x)\\[2mm]
&=&B_\fg(a, [b,c]_\fg)\\[2mm]
&=&B_\fg([a,b]_\fg, c)\\[2mm]
&=&B_\fg([a,b]_\fa+B_\alpha(a,b)x, c)\\[2mm]
&=&B_\fg([a,b]_\fa, c)\\[2mm]
&=&B_\fa([a,b]_\fa, c).\\
\end{array}
\]
Similarly, one can prove  the invariance property when $a, b, c\in \fa_\od$. 

Next, we show that the bilinear form $B_\fa$ satisfies the condition 
\[
B_\fa(D(a),b)=B_\alpha(a,b) \quad \text{ for any $a,b\in\fa_\od$. }
\] 
which, in turn, implies Eq. (\ref{cond}). We have,
\[
\begin{array}{lcl}
B_\fa(D(a),b)&=&B_\fg(D(a),b) \\[2mm]
&=&B_\fg([x^*,a]_\fg,b)\\[2mm]
&=&B_\fg(x^*, [a,b]_\fg)\\[2mm]
&=&B_\fg(x^*, [a,b]_\fa+B_\alpha(a,b)x)\\[2mm]
&=&B_\alpha(a,b).
\end{array}
\]
The proof now is complete.
\end{proof}
\sssbegin{Remark}
Proposition \ref{Rec1} is an attempt to generalize the result of [MR] that any non-simple NIS-Lie algebra for $p=0$ can be obtained as a result of a double extension (or a {\it generalized} double extension) of another NIS-Lie algebra. Several generalizations were  attempted in the papers \cite{BB}, \cite{B} and \cite{BBB} in the case of NIS-Lie superalgebras for $p=0$. The proof of Proposition \ref{Rec1} we provide has been modified a lot to take the squaring into account, and overcome the difficulty that, for Lie superalgebras in characteristic $p=2$, if $I$ is an ideal then $I^\perp$ is not necessarily an ideal.
\end{Remark}

\subsection{$D_\od$-extensions}

\sssbegin{Theorem}\label{MainTh2} Let $(\fa,B_\fa)$ be a NIS-Lie superalgebra in characteristic $2$ such that $B_\fa$ is even. Let $D\in\fder_\od(\fa)$ and $a_0\in \fa_\ev$ satisfy the following conditions:
\begin{eqnarray}
\label{2D1} B_\fa(D(a),b)&=& B_\fa(a,D(b)) \text{ for any } a,b \in\fa;\\[2mm]
\label {2D2} D^2&=&\ad_{a_0};\\[2mm]
\label{2D3} D(a_0)&=&0.
\end{eqnarray}
Then there exists a NIS-Lie superalgebra structure on $\fg:=\mathscr{K} \oplus \fa \oplus \mathscr{K}^*$, where $\mathscr{K}:=~\Span\{x\}$ and $x$ odd, defined as follows. The squaring is given by 
\[
s_\fg(rx+a+t x^*):= s_\fa(a) + t^2 a_0+t D(a) \qquad \text{ for any }\; rx+a+tx^* \in \fg_\od.
\]
The bracket is given by:
\[
\begin{array}{l}
[x,\fg]_\fg:=0; \qquad [a,b]_\fg:=[a,b]_\fa+B_\fa(D(a),b)x\quad\text{~~for any $a,b\in\fa$};\\[2mm]
[x^*,a]_\fg:=D(a) \quad  \text{ for any } a \in \fa.
\end{array}
\]
The bilinear form $B_\fg$ on $\fg$ defined by:
\begin{eqnarray*}
&{B_\fg}\vert_{\fa \times \fa}:= B_\fa, \quad B_\fg(\fa,\mathscr{K}):=0, \quad B_\fg(\fa,\mathscr{K}^*):=0, \\ 
&B_\fg(x,x^*):=1, \quad B_\fg(x,x):= B_\fg(x^*,x^*):=0.
\end{eqnarray*}
is even, non-degenerate, symmetric and invariant on $(\fg, [\cdot ,\cdot]_\fg, s_\fg)$.

Therefore $(\fg,B_\fg)$ is a NIS-Lie superalgebra. 
\end{Theorem}

We call the  NIS-Lie superalgebra $(\fg,B_\fg)$ constructed in Theorem 
\ref{MainTh2} a \textit{$D_\od$-extension} of $(\fa, B_{\fa})$  by means of $D$ and $a_0$.
\begin{proof}
The proof is similar to that of Theorem \ref{MainTh}. Let us first show that $s_\fg$ is indeed a squaring on $\fg$. Recall that since $\mathscr{K}$ and $\mathscr{K}^*$ are odd vector spaces, then $\fg_\od=\mathscr{K} \oplus \fa_\od\oplus \mathscr{K}^*$. Now, let $\lambda\in\Bbb K$ and $f=rx+a+tx^*\in\fg_\od$; we have
\[
s_\fg(\lambda f)= s_\fa(\lambda a) + (\lambda t)^2 a_0+t \lambda D(\lambda a)=\lambda ^2 s_\fa(a) +\lambda^2 t^2 a_0+t \lambda^2 D(a)=\lambda^2s_\fg(f),
\]
since $s_\fa$ is a squaring on $\fa$ and $D$ is a linear map. Besides, for any $f=rx+a+tx^*$ and $g=ux+b+vx^*$ in $\fa_\od$, we see that
\begin{equation}
\label{squaring2}
\begin{array}{lcl}
s_\fg(f+g)+s_\fg(f)+s_\fg(g)&=&s_\fa(a+b)+(t+v)^2 a_0+(t+v)D(a+b)\\[2mm]
&&+s_\fa(a)+t^2 a_0+tD(a)+s_\fa(b)+v^2 a_0+vD(b)\\[2mm]
&=&[a,b]_\fa+vD(a)+tD(b)
\end{array}
\end{equation}
is obviously bilinear since it is expressed in terms of two bilinear maps. 
Let us show that the bracket on $\fg$ defined above is supersymmetric. Indeed, using condition (\ref{2D1}), for any $f=rx+a+tx^*$ and $g=\tilde rx+b+\tilde tx^*$ in $\fg$, we have
\begin{eqnarray*}
[f, g]_\fg&=&B_\fa(D(a),b)x+[a,b]_\fa+t D(b)+\tilde t D(a)\\
&=& B_\fa(D(b),a)x+[b,a]_\fa+t D(b)+\tilde t D(a)\\
&=& [g,f]_\fg.
\end{eqnarray*}

Let us check the Jacobi identity relative to the squaring $s_\fg$. Let $f=r x+a+t x^*\in\fg_\od$ and $g \in \fg$. If $g=x$, then we are done since $x$ is central. If $g=x^*$, we have
\[
\begin{array}{lcl}
[s_\fg(f),g]_\fg+[f,[f,g]_\fg]_\fg&=&[s_\fa(a) +t^2 a_0+t D(a),x^*]_\fg +[f, D(a) ]_\fg\\[2mm]
&=& D(s_\fa(a) +t^2 a_0+t D(a))\\[2mm]
&&+B_\fa(D(a),D(a))x+ [a,D(a)]_\fa+t D(D(a))=0,
\end{array}
\]
since $D\in \fder(\fa)$, condition (\ref{2D3}) holds and $B_\fa(D(a),D(a))=0$. Indeed, using conditions (\ref{2D2}) and (\ref{2D3}) of Theorem \ref{MainTh2}, we have
\[
B_\fa(D(a), D(a))=B_\fa(a, D^2(a))= B_\fa(a, [a_0,a])=B_\fa(a_0,[a,a])=0.
\] 
If $g=b\in\fa_\od$, then we have
\[
\begin{array}{lcl}
[s_\fg(f),g]_\fg+[f,[f,g]_\fg]_\fg&=&[s_\fa(a) +t^2 a_0+t D(a),b]_\fg+[f, [f,g]_\fg]_\fg \\[2mm]
&=&B_\fa(D(s_\fa(a) +t^2 a_0+t D(a)),b)x+[s_\fa(a) +t^2 a_0+t D(a),b]_\fa\\[2mm]
&&+B_\fa(D(a), [f,g]_\fg)x +[a,  [f,g]_\fg]_\fa +t D([f,g]_\fg)\\[2mm]
&=&B_\fa([Da,a] +t D^2(a),b)x+[s_\fa(a) ,b]_\fa+t^2 [a_0,b]_\fa+t [D(a),b]_\fa\\[2mm]
&&+B_\fa(D(a),[a,b]_\fa+ t D(b))x +[a,  [a,b]_\fa+ t D(b)]_\fa\\[2mm]
&&+t D( [a,b]_\fa+ t D(b))\\[2mm]
&=&B_\fa([Da,a] +t D^2(a),b)x+B_\fa(D(a),[a,b]_\fa+ t D(b))x\\[2mm]
&&+t^2 [a_0,b]_\fa+t^2 D^2(b)+[s_\fa(a) ,b]_\fa+[a,  [a,b]_\fa]_\fa\\[2mm]
&&+[a, t D(b)]_\fa+t D( [a,b]_\fa)+t [D(a),b]_\fa=0,
\end{array}
\]
since $s_\fa$ is a squaring, $D$ is a derivation, $B_\fa$ is $\fa$-invariant and $D^2=\ad_{a_0}$.

Now, if $g=b\in\fa_\ev$, we have
\[
\begin{array}{lcl}
[s_\fg(f),g]_\fg+[f,[f,g]_\fg]_\fg&=& [s_\fa(a) +t^2 a_0+t D(a),b]_\fg \\[2mm]
&&+ [f,B_\fa(D(a),b)x+[a,b]_\fa+t D(b)]_\fg \\ [2mm]
&=&B_\fa(D(s_\fa(a) +t^2 a_0+t D(a)),b)x+[s_\fa(a) +t^2 a_0+t D(a),b]_\fa \\ [2mm]
&&+B_\fa(D(a),[a,b]_\fa+tD(b)))x\\[2mm]
&&+[a,B_\fa(D(a),b)x+[a,b]_\fa+tD(b)]_\fa +t D([a,b]_\fa+tD(b))\\[2mm]
&=&B_\fa([Da,a]_\fa +t \ad_{a_0}(a)),b)x+[s_\fa(a) +t^2 a_0+t D(a),b]_\fa \\ [2mm]
&&+B_\fa(D(a),[a,b]_\fa+tD(b)))x\\[2mm]
&&+[a,[a,b]_\fa]_\fa+[a,t D(b)]+t D([a,b]_\fa+t D(b))\\[2mm]
&=&B_\fa([Da,a]_\fa +t \ad_{a_0}(a)),b)x+B_\fa(D(a),[a,b]_\fa+t D(b)))x\\[2mm]
&&+t^2 [a_0,b]_\fa +t^2 D^2(a) +[a,[a,b]_\fa]_\fa+[s_\fa(a),b]_\fa\\[2mm]
&&+t [a,D(b)]+t[D(a),b]+t D([a,b]_\fa)=0,
\end{array}
\]
since $D$ is a derivation, $s_\fa$ is a squaring, $B_\fa$ is $\fa$-invariant and $D^2=\ad_{a_0}$.

To check the Jacobi identity 
\[
[h,[f,g]_\fg]_\fg+[g,[h,f]_\fg]_\fg+[f,[g,h]_\fg]_\fg=0, 
\]
we proceed as follows. 

If $h=x$, the identity
\[
[x,[f,g]_\fg]_\fg+[g,[x,f]_\fg]_\fg+[f,[g,x]_\fg]_\fg=0
\]
is certainly satisfied since $x$ is central. 

If $h=x^*$, the identity 
\[
[x^*,[f,g]_\fg]_\fg+[g,[x^*,f]_\fg]_\fg+[f,[g,x^*]_\fg]_\fg=0
\]
is also satisfied for the following reasons. If either $f$ or $g$ is $x$, then we are done since $x$ is central. Now if $f=x^*$ (or the way around $g=x^*$), then by putting $g=ux+b+vx^*$ we get
\[
[x^*,[x^*,g]_\fg]_\fg+[g,[x^*,x^*]_\fg]_\fg+[x^*,[g,x^*]_\fg]_\fg=2D(D(b))=0.
\]
Let us assume now that $f, g\in\fa$. We see that if $f$ is even and $g$ is odd, then 
\[
\begin{array}{lcl}
[x^*,[f,g]_\fg]_\fg+[g,[x^*,f]_\fg]_\fg+[f,[g,x^*]_\fg]_\fg&=&[x^*, B_\fa(D(f),g)x+[f,g]_\fa]_\fg+ [g, D(f)]_\fg\\[2mm]
&&+ [f,D(g)]_\fg\\[2mm]
&=&D([f,g]_\fa)+ [g, D(f)]_\fa+B_\fa(D(f),D(g))x\\[2mm]
&&+ [f,D(g)]_\fa=0,
\end{array}
\]
since $D\in \fder(\fa)$ and $B_\fa(D(f),D(g))=0$ as $B_\fa$ is even. If $f$ and $g$ are both even, we deduce, using the fact that $B_\fa$ is symmetric, that 
\[
\begin{array}{lcl}
[x^*,[f,g]_\fg]_\fg+[g,[x^*,f]_\fg]_\fg+[f,[g,x^*]_\fg]_\fg&=&[x^*, B_\fa(D(f),g)x+[f,g]_\fa]_\fg+ [g, D(f)]_\fg\\[2mm]
&&+ [f,D(g)]_\fg\\[2mm]
&=&D([f,g]_\fa)+ B_\fa(D(g),D(f))x+ [g, D(f)]_\fa\\[2mm]
&&+ B_\fa(D(f),D(g))x+ [f,D(g)]_\fa=0.
\end{array}
\]
Besides, if $f, g\in\fa_\od$, then
\[
\begin{array}{lcl}
[x^*,[f,g]_\fg]_\fg+[g,[x^*,f]_\fg]_\fg+[f,[g,x^*]_\fg]_\fg&=&[x^*, [f,g]_\fa]_\fg+ [g, D(f)]_\fg+ [f,D(g)]_\fg\\[2mm]
&=&D[f,g]_\fa+[g,D(f)]_\fa+[f,D(g)]_\fa=0,
\end{array}
\]
since $D\in \fder(\fa)$. 

From now and on we will assume that $f,g,h\in\fa$. We distinguish several cases to check the Jacobi identity.

If $f, g$ and $h$ are even, then we have \[
\begin{array}{lcl}
[f,[g,h]_\fg]_\fg+ \circlearrowleft (f,g,h)&=&[f,B_\fa(D(g),h)x+[g,h]_\fa]_\fg + \circlearrowleft (a,b,c)\\[2mm]
&=&[f,[g,h]_\fa]_\fa+B_\fa(D(f),[g,h]_\fa)x + \circlearrowleft (a,b,c) =0,
\end{array}
\]
because the JI holds on $\fa$ and Lemma \ref{lemma1}. 
If $f$ and $g$ are both odd but $h$ is even, then
\[
\begin{array}{lcl}
[f,[g,h]_\fg]_\fg +[h,[f,g]_\fg]_\fg+[g,[h,f]_\fg]_\fg&=&[f,B_\fa(D(g),h)x+[g,h]_\fa]_\fg+[h,[f,g]_\fa]_\fg \\[2mm]
&& +[g,B_\fa(D(h),f)x+[h,f]_\fa]_\fg\\[2mm]
&=& [f,[g,h]_\fa]_\fa+B_\fa(D(h),[f,g]_\fa)x+[h,[f,g]_\fa]_\fa\\[2mm]
&&+B_\fa(D(g),[h,f]_\fa)x+[g,[h,f]_\fa]_\fa=0,
\end{array}
\]
because the Jacobi identity holds on $\fa$ and $B_\fa(D(h),[f,g]_\fa)=B_\fa(D(g),[h,f]_\fa)=0$ as $B_\fa$ is even.

If $f$ is odd and both $g$ and $h$ are even, using Lemma \ref{lemma1} we have
\[
\begin{array}{lcl}
[f,[g,h]_\fg]_\fg+[h,[f,g]_\fg]_\fg+[g,[h,f]_\fg]_\fg&=&[f,B_\fa(D(g),h)x+[g,h]_\fa]_\fg+ [h,B_\fa(D(f),g)x+[f,g]_\fa ]_\fg\\[2mm]
&&+[g,B_\fa(D(h),f)x+[h,f]_\fa ]_\fg\\[2mm]
&=&[f,[g,h]_\fa]_\fa+[h,[f,g]_\fa]_\fa+[g,[h,f]_\fa]_\fa ,\\[2mm]
&&+B_\fa(D(f),[g,h]_\fa)x+B_\fa(D(h),[f,g]_\fa)x\\[2mm]
&&+B_\fa(D(g),[h,f]_\fa)x=0.
\end{array}
\] 
Let us show now that $B_\fg$ is invariant. We should check that
\[
B_\fg([f,g]_\fg, h)=B_\fg (f, [g,h]_\fg) \quad \text{for any $f,g,h\in\fg$}.
\]
This is true if $f=x$ because $x$ is central and by the very definition of $B_\fg$.  If $f=x^*$, $g=ux+b+vx^*$ and $h=wx+c+zx^*$ (where $g$ and $h$ are both odd), using Eq. (\ref{squaring2}) and the fact that $B_\fa$ is even we have
\[
\begin{array}{lcl}
B_\fg([f,g]_\fg, h)&=& B_\fg(D(b), wx+c+zx^*)\\[2mm]
&=& B_\fa(D(b), c)=0.
\end{array}
\]
On the other hand, 
\[
B_\fg(f,[g, h]_\fg)=B_\fg(x^*,[b,c]_\fa+zD(b)+vD(c))=0.
\] 
If $f=x^*$, $g=ux+b+vx^*$ is odd and $h$ is even, we have
\[
\begin{array}{lcl}
B_\fg([f,g]_\fg, h)&=&B_\fa(D(b), h).
\end{array}
\]
On the other hand,
\[
\begin{array}{lcl}
B_\fg(f,[g, h]_\fg)&=&B_\fg(x^*, B_\fa(D(b),h)x+[b,h]_\fa + vD(h))\\[2mm]
&=&B_\fa(D(b), h) .\\[2mm]
\end{array}
\]
If $f=x^*$, but $g$ and $h$ are both even, we have
$ B_\fg([f,g]_\fg, h)=B_\fa(D(g), h)$.
On the other hand, 
\[
\begin{array}{lcl}
B_\fg(f,[g, h]_\fg)&=&B_\fg(x^*, B_\fa(D(g),h)x+[g,h]_\fa )\\[2mm]
&=&B_\fa(D(g), h) .\\[2mm]
\end{array}
\]
If $f,g,h\in\fa$, the invariance property of $B_\fg$ directly follows from that of $B_\fa$.
\end{proof}

We define the cone 
\[
\mathscr{C}(\fg, B_\fg):=\{x\in \fg_\od\; | \; B_\fg(s_\fg(x), s_\fg(t))=0 \quad \text{ for any } t\in \fg_\od \}.
\]
Clearly, if $x\in\mathscr{C}(\fg, B_\fg)$, then $\lambda x\in\mathscr{C}(\fg, B_\fg)$ for any $\lambda\in\Kee$, hence the terminology. 

\sssbegin{Theorem}
\label{Rec2}
Let $(\fg, B_\fg)$ be a NIS-Lie superalgebra in characteristic $2$. Let us suppose $\fz(\fg)_\od \cap~\mathscr{C} (\fg, B)\not=\{0\}$. Then $(\fg,B_\fg)$ is obtained as a $D_\ev$-extension or a $D_\od$-extension from a NIS-Lie superalgebra $(\fa,B_\fa)$ of dimension $\dim(\fg) - 2$.
\end{Theorem}

\begin{proof} Let $x\in \fz_s(\fg)_\od \cap \mathscr{C} (\fg, B_\fg)$. We assume that $s_\fg(x)=0$; otherwise, $s_\fg(x)$ will be a non-zero element of $\fz_s(\fg)$ and hence a $D_\ev$-extension can be constructed by means of this element following the steps of Theorem \ref{Rec1}. 

Now, the subspace $\mathscr{K}:= \Span\{x\}$ is an ideal in $(\fg,[\cdot ,\cdot]_\fg, s_\fg)$  because $x$ is central and $s_\fa(x)=0$. Moreover, the subspace $\mathscr{K}^\perp$ is also an ideal containing $\mathscr{K}$. Indeed, let us first show that $(\mathscr{K}^\perp)_\ev=\fg_\ev$.  This is true because $x$ is orthogonal to any even element since $B_\fg$ is even. Now let $f\in (\mathscr{K}^\perp)_\od$. It follows that $B_\fg(x,s_\fg(f))=0$ since $x\in \fz_s(\fg)_\od$. Hence, $s_\fg(f) \in \mathscr{K}^\perp$. On the other hand, if $f\in \mathscr{K}^\perp$ and $g\in \fg$, then $[f,g]_\fg \in \mathscr{K}^\perp$ because
\[
B_\fg(x,[f,g]_\fg)=B_\fg([x,f]_\fg,g)=B_\fg(0,g)=0.
\]
Same arguments as in Theorem \ref{Rec1} can be used to construct the ideal $\fa=(\mathscr{K} \oplus \mathscr{K}^*)^\perp,$ 
and a decomposition $\fg=\mathscr{K} \oplus \fa\oplus \mathscr{K}^*$, where the generator $x^*$ of $ \mathscr{K}^*$ can be normalized so that $B_\fg(x,x^*)=1$. 

Let us define a bilinear from on $\fa$ by setting:
\[
B_\fa:={B_\fg}\vert_{\fa \times \fa}.
\]
The form $B_\fa$ is non-degenerate on $\fa$. Indeed, suppose there exists an $a\in \fa$ such that $B_\fa(a,f)=0$ for any $f\in \fa$. But $a$ is also orthogonal to $x$ and to $x^*$. It follows that $B_\fa(a,f)=0$ for any $f\in \fg$. Hence, $a=0$.

Let $a, b \in \fa$. Since $\fa$ is an ideal, it follows that the bracket $[a,b]_\fg\in \mathscr{K}\oplus \fa$ because $\fa \subset  \mathscr{K}^\perp \oplus \fa$ and the latter is an ideal. The same arguments as those used in Theorem~\ref{Rec1} yield
\[
\begin{array}{lcl}
[x^*,a]_\fg&=&D(a),\\[2mm]
[a,b]_\fg&=&B_\fg(D(a),b)x+\psi(a,b),
\end{array}
\]
where $D(a)$ and $\psi(a,b)$ are both in $\fa$. Now, let $f\in \fg_\od$. It follows that $s_\fg(f)\in \fa$ since $\fg_\ev=\fa_\ev$. Let us then write $s_\fg(x^*)=a_0$. The fact that $[s_\fg(x^*),a]_\fg=[x^*,[x^*,a]_\fg]_\fg$ for any $a\in \fa$ implies $B_\fg(D(a_0),a)x+[a_0,a]_\fa=D^2(a)$. Therefore, 
\[
D^2=\ad_{a_0}, \quad \text{ and } \quad  B_\fa(D(a_0),a)=0.
\]
Besides, $[s_\fg(x^*),x^*]=0$ implies that $D(a_0)=0$. 

Similarly, $[s_\fg(a),x^*]_\fg=D(s_\fg(a))$ implies that $[a,[a,x^*]_\fg]_\fg=D(s_\fg(a))$ which, in turn, implies 
\[
D(s_\fg(a))=B_\fg(D(a), D(a))x+[a,D(a)]_\fa. 
\]
But 
\[
B_\fg(D(a), D(a))=B_\fg(a,D^2(a))=B_\fg(a,[a_0,a])=0, 
\] 
hence $D(s_\fg(a))=[a,D(a)]_\fa$. The arguments used in Theorem \ref{MainTh2} show that 
\[
D([a,b]_\fa)=[D(a),b]_\fa+[a,D(b)]_\fa,
\] 
and hence $D\in \fder(\fa)$.

Now, we shall prove that the new defined bilinear form $B_\fa:={B_\fg}\vert_{\fa \times \fa}$ is also invariant on $(\fa, [\cdot ,\cdot]_\fa, s_\fa, B_\fa)$. Indeed, for any $a,b$ and $c$ in $\fa$ (at least two of them even), we have
\[
\begin{array}{lcl}
B_\fa(a,[b,c]_\fa)&=&B_\fg(a, [b,c]_\fg+B_\fa(D(b),c)x)\\[2mm]
&=&B_\fg(a, [b,c]_\fg)\\[2mm]
&=&B_\fg([a,b]_\fg, c)\\[2mm]
&=&B_\fg([a,b]_\fa+B_\fa(D(a),b)x, c)\\[2mm]
&=&B_\fg([a,b]_\fa, c)\\[2mm]
&=&B_\fa([a,b]_\fa, c).\\
\end{array}
\]
Now, for any $a,b$ and $c$ in $\fa$ (where $a$ and $b$ are odd and $c$ is even), we have
\[
\begin{array}{lcl}
B_\fa(a,[b,c]_\fa)&=&B_\fg(a, [b,c]_\fg+B_\fa(D(b),c)x)\\[2mm]
&=&B_\fg(a, [b,c]_\fg)\\[2mm]
&=&B_\fg([a,b]_\fg, c)\\[2mm]
&=&B_\fg([a,b]_\fa, c)\\[2mm]
&=&B_\fa([a,b]_\fa, c).\\
\end{array}
\]
Similarly, one can prove that the invariance property is satisfied when $a$, $b$ and $c$ are all odd. 

The proof is complete now. \end{proof}

\subsection{Isometries, and equivalence classes of derivations} \label{ST1}
For a  NIS-Lie superalgebra $\fa$ with a bilinear form $B_\fa$, denote by $\fg$ (resp. $\tilde \fg$) the double extension of $\fa$ by means of a derivation $D$ (resp. $\tilde D$). In the case of $D_\ev$-extensions, $\fg$ (resp. $\tilde \fg$) is also defined by means of a quadratic form $\alpha$ (resp. $\tilde \alpha$).  An \textit{isometry} between $\fg$ and $\tilde \fg$ is an isomorphism $\pi:\fg\rightarrow \tilde \fg$ such that: 
\[
\begin{array}{rcll}
\pi([f,g]_\fg)&=&[\pi(f), \pi(g)]_{\tilde \fg}, & \text{for any $f\in \fg_\ev$ and $f\in \fg$},\\[2mm]
\pi(s_\fg(f))&=&s_{\tilde \fg}(\pi(f)), &  \text{for any $f\in \fg_\od$},\\[2mm]
B_{\tilde \fg}(\pi(f), \pi(g))&=&B_\fg(f,g), & \text{for any $f,g\in \fg$}.
\end{array}
\]

We will investigate how the derivations $D$ and $\tilde D$ are related with each other when $\fg$ and $\tilde \fg$ are isometric. We will assume further that the isometry satisfies $\pi(\mathscr{K} \oplus \fa)=\tilde{\mathscr{K}} \oplus \fa$, and call it an {\it adapted isometry}. Hereafter, $\mathscr{K}=\Span\{x\}$,  $\mathscr{K^*}=\Span\{x^*\}$, $\tilde {\mathscr{K}}=\Span\{\tilde x\}$, and $\tilde {\mathscr{K}}^*=\Span\{\tilde x^*\}$.

\subsection{The $D_\ev$ case} Let $\pr: \tilde{\mathscr{K}} \oplus \fa\rightarrow \fa$ be the projection, and  $\pi_0:=\pr \circ \pi$. The map $\pi_0$ is obviously linear. Let $a\in \fa_\ev$. Since $\pi(a)-\pi_0(a)\in\mathrm{Ker}(\pr)$, it follows that $\pi(a)+\pi_0(a)\in  \tilde{\mathscr{K}} $. Since $B_\fa$ is nondegenerate, there exists a unique $t_\pi\in\fa$ (depending only in $\pi$) such that 
\[
\pi(a)+\pi_0(a)=B_\fa(t_\pi,a) \tilde x \quad \text{for any $a\in\fa_\ev$.}
\]
Let now $a\in \fa_\od$. Since $\pi(a)\in\fa$, it follows that $\pi(a)=\pi_0(a)$. 

Besides, $\pi(x)=\lambda \tilde x$ for some $\lambda\in\Kee$. Indeed, let us write $\pi(x)=\lambda \tilde x+a$, where $a\in \fa$. We have (for any $b$ even)
\[
0=B_\fg(x,b)=B_{\tilde \fg}(\pi(x), \pi(b))=B_{\tilde \fg}(\lambda \tilde x+a, \pi_0(b)+B_\fa(t_\pi, b) \tilde x)=B_\fa(a,\pi_0(b)). 
\]
We have (for any $b$ odd)
\[
0=B_\fg(x,b)=B_{\tilde \fg}(\pi(x), \pi(b))=B_{\tilde \fg}(\lambda \tilde x+a, \pi_0(b))=B_\fa(a,\pi_0(b)). 
\]
Since $\pi_0$ is surjective and $B_\fa$ is nondegenerate, it follows that $a=0$.
Let us show that $\pi_0$ preserves $B_\fa$. Indeed, for any $a,b\in\fa_\ev$, we have
\[
\begin{array}{lcl}
B_\fa(a,b)&=&B_{\tilde \fg}(\pi(a), \pi(b))\\[2mm]
&=& B_{\tilde \fg}(\pi_0(a)+B_\fa(t_\pi, a)\tilde x, \pi_0(b)+B_\fa(t_\pi, b)\tilde x)\\[2mm]
&=&B_{\tilde \fg}(\pi_0(a), \pi_0(b))\\[2mm]
&=&B_{\fa}(\pi_0(a), \pi_0(b)).
\end{array}
\] 
The same arguments can be used if $a$ and $b$ are both odd (or one of them is odd). 

Let us study the squaring. Let $a\in \fg_\od=\fa_\od$. We have 
\[
\pi(s_\fg(a))=\pi(s_\fa(a)+\alpha(a)x)=\pi_0(s_\fa(a))+B_\fa(t_\pi, s_\fa(a))\tilde x+\alpha(a)\lambda \tilde x.
\]
On the other hand, 
\[
s_{\tilde \fg}(\pi(a))=s_{\tilde \fg}(\pi_0(a))=s_{\fa}(\pi_0(a))+\tilde \alpha(\pi_0(a))\tilde x.
\]
It follows that (for any $a$ odd)
\begin{eqnarray}
\nonumber \tilde \alpha(\pi_0(a))+\lambda \alpha(a)+B_\fa(t_\pi, s_\fa(a))&=&0, \; \text{ and }\\[2mm] 
\label{pi01} \pi_0(s_\fa(a))&=&s_\fa(\pi_0(a)).
\end{eqnarray}
We need the following Lemma:

\ssbegin{Lemma}
If 
\begin{equation}
\label{alph}
\tilde \alpha \circ \pi_0+\lambda \alpha+B_\fa(t_\pi, s_\fa(\cdot))=0,
\end{equation}
then $\pi_0^{-1}\tilde D\pi_0+\lambda D+\ad_{t_\pi}=0.$
\end{Lemma}
\label{Deven}

\begin{proof} Let $a, b \in \fa_\od$. Evaluating Eq. (\ref{alph}) at $a+b$, at $a$, and at $b$, and taking the sum of evaluations we get
\[
B_{\tilde \alpha}( \pi_0(a), \pi_0(b))+\lambda B_\alpha(a,b)+B_\fa(t_\pi, [a,b]_\fa)=0.
\]
Using the fact that $B_\alpha(a,b)=B_\fa(D(a),b)$ and $B_{\tilde \alpha}(\pi_0(a),\pi_0(b))=B_\fa(\tilde D(\pi_0(a)),\pi_0(b))$ we get 
\[
B_{\fa}(\tilde D( \pi_0(a)), \pi_0(b))+\lambda B_\fa(D(a),b)+B_\fa([t_\pi,a],b)=0.
\]
On the other hand, $B_{\fa}(\pi_0^{-1}\tilde D( \pi_0(a))+\lambda D(a)+[t_\pi,a],c)=0$, for any even element $c$, since $B_\fa$ is even. The result follows since $B_\fa$ is nondegenerate. 
\end{proof}

Now, the fact that 
\[
\pi([a,b]_\fg)=[\pi(a), \pi(b)]_{\tilde \fg} \text{ for any $a,b\in\fa_\ev$,}
\] implies that
\[
\pi_0([a,b]_\fa)+B_\fa(D(a),b)\pi(\tilde x) + B_\fa(t_\pi, [a,b]_\fa) \tilde x= [\pi_0(a), \pi_0(b)]_\fa+B_\fa(\tilde D(\pi_0(a)), \pi_0(b))\tilde x,
\]
since $x$ is central in $\fg$. It follows that 
\begin{eqnarray}
\nonumber \lambda B_\fa(D(a),b)+B_\fa([t_\pi,a]_\fa,b)+B_\fa(\pi_0^{-1}D\pi_0(a),b)&=&0, \, \text{ and }\\[2mm]
\label{pi02}\pi_0([a,b]_\fa)&=&[\pi_0(a), \pi_0(b)]_\fa
\end{eqnarray}
which implies that $\pi_0^{-1}\tilde D\pi_0+\lambda D+\ad_{t_\pi}=0 $. 

For all $a$ even and $b$ odd, we have $\pi([a,b]_\fg)=\pi([a,b]_\fa)=\pi_0([a,b]_\fa)$. On the other hand, 
\[
[\pi(a), \pi(b)]_{\tilde \fg}=[\pi_0(a)+B_\fa(t_\pi, a)\tilde x, \pi_0(b)]_{\tilde \fg}=[\pi_0(a),\pi_0(b)]_{\tilde \fg}=[\pi_0(a),\pi_0(b)]_{\fa}.
\]
It follows that 
\begin{equation}
\label{pi03}
\pi_0([a,b]_\fa)=[\pi_0(a),\pi_0(b)]_{\fa}.
\end{equation}

For all $a$ and $b$ odd, we have
\[
\begin{array}{lcl}
\pi([a,b]_\fg)=\pi([a,b]_\fa+B_\alpha(a,b) x)&=&\pi([a,b]_\fa)+\lambda B_\alpha(a,b)\tilde x\\[2mm]
&=&\pi_0([a,b]_\fa)+B_\fa(t_\pi,[a,b]_\fa)\tilde x+\lambda B_\alpha(a,b)\tilde  x\\[2mm]
&=&\pi_0([a,b]_\fa)+B_\fa(t_\pi,[a,b]_\fa)\tilde x+\lambda B_\fa(D(a),b)\tilde  x.
\end{array}
\]
On the other hand, 
\[
[\pi(a), \pi(b)]_{\tilde \fg}=[\pi_0(a), \pi_0(b)]_\fa+B_{\tilde \alpha}(\pi_0(a), \pi_0(b))\tilde x=[\pi_0(a), \pi_0(b)]_\fa+B_\fa(\tilde D(\pi_0(a)), \pi_0(b))\tilde x.
\]
It follows that 
\begin{equation}
\label{pi04}
[\pi_0(a), \pi_0(b)]_\fa=\pi([a,b]_\fa),
\end{equation}
and 
\[
B_\fa(\tilde D(\pi_0(a)), \pi_0(b))=B_\fa(t_\pi,[a,b]_\fa)+\lambda B_\fa(D(a),b).
\] 
This  condition results from Lemma \ref{Deven}.

Now, Eqs, (\ref{pi01}), (\ref{pi02}), (\ref{pi03}) and (\ref{pi04}) imply that $\pi_0$ is an automorphism on $\fa$. 

Let us describe $\pi(x^*)$ now. We write $\pi(x^*)$ as $\mu \tilde x^*+a+\nu \tilde x$ for some $a\in\fa_\ev$. We have
\[
1=B_{\tilde \fg}(\pi(x^*), \pi(x))=B_{\tilde \fg}(\mu \tilde x^*+a+\nu \tilde x, \lambda \tilde x)=\lambda \mu.
\]
Therefore $\mu=\lambda^{-1}$. Besides,
\[
\begin{array}{lcl}
B_\fg(x^*, x^*)=B_{\tilde \fg}(\pi(x^*), \pi(x^*))&=&B_{\tilde \fg}(\mu \tilde x^*+a+\nu \tilde x, \mu \tilde x^*+a+\nu \tilde x)\\[2mm]
&=&\mu^2 B_{\tilde \fg}(\tilde x^*, \tilde x^*)+B_\fa(a,a).
\end{array}
\]
Therefore, $B_\fg(x^*,x^*)=B_\fa(a, a)+\mu^2 B_\fg(\tilde x^*, \tilde x^*)$. Besides, for any $b\in \fa_\ev$, we have
\[
0=B_{\tilde \fg}(\pi(x^*), \pi(b))=B_{\tilde \fg}(\mu \tilde x^*+a+\nu \tilde x, \pi_0(b)+B_\fa(t_\pi,b)\tilde x)=\mu B_\fa(t_\pi,b)+B_\fa(a, \pi_0(b)).
\]
Since $B_\fa$ is even and non-degenerate, it follows that $\mu t_\pi+\pi_0^{-1}(a)=0$ which implies that $a=\mu \pi_0(t_\pi)$. Finally, we get
\[
\pi(x^*)=\lambda^{-1}(\tilde x^*+\pi_0(t_\pi))+\nu \tilde x.
\]
We arrive at the following Theorem.

\sssbegin{Theorem} \label{Isom1}
Let $\pi_\ev$ be an isometry of $(\fa, B_\fa)$.  Let $\lambda \in \Kee^{\times}$ and let $t\in \fa_\ev$, satisfying the following conditions:
\begin{eqnarray}
\label{Ca}\tilde \alpha&=&\lambda \alpha \circ \pi_0^{-1}+B_\fa(t, s_\fa \circ \pi_0^{-1}) \quad \text{on $\fa_\od$}; \\[2mm]
\label{Cd}\pi_0^{-1}\tilde D\pi_0&=&\lambda D+\ad_{t} \quad \text{on $\fa_\ev$}; \\[2mm]
B_\fg(x^*,x^*)&=&\lambda^{-2}( B_\fa(t, t)+B_\fg(\tilde x^*, \tilde x^*)).
\end{eqnarray}
Then there exists an adapted isometry $\pi: \mathscr{K} \oplus \fa \oplus \mathscr{K}^*\rightarrow  \tilde{\mathscr{K}} \oplus \fa \oplus \tilde{\mathscr{K}^*}$ given by 
\[
\begin{array}{lcl}
\pi&=& \pi_0+ B_\fa(t,\cdot) \tilde x \quad \text{on $\fa$;}\\[2mm]
\pi(x)&=&\lambda \tilde x;\\[2mm]
\pi(x^*)&=&\lambda^{-1}(\tilde x^*+\pi_0(t))+\nu \tilde x,\; \text{where $\nu$ is arbitrary.}\\
\end{array}
\]
\end{Theorem}

\begin{proof} To check that $\pi$ preserves the Lie bracket, it is enough to check the conditions below.  For every $a$ even, we have
\[
\pi([x^*,a]_\fg)=\pi (D(a))=\pi_0(D(a))+ B_\fa(t, D(a))\tilde x.
\]
On the other hand, 
\[
\begin{array}{lcl}
[\pi(x^*),\pi(a)]_{\tilde \fg}&=&[\lambda^{-1}\tilde x^*+\lambda^{-1}\pi_0(t)+\nu \tilde x , \pi_0(a)+B_\fa(t,a) \tilde x]_{\tilde \fg}\\[2mm]
&=&\lambda^{-1}\tilde D(\pi_0(a)) +[\lambda^{-1}\pi_0(t), \pi_0(a)]_\fa+B_\fa(\tilde D(\lambda^{-1}\pi_0(t)), \pi_0(a))\tilde x\\[2mm]
&=&\lambda^{-1} \pi_0 (\lambda D(a)+ [t,a]) +\lambda^{-1}\pi_0([t, a]_\fa)+\lambda^{-1}B_\fa(\tilde D(\pi_0(t)), \pi_0(a))\tilde x\\[2mm]
&=&\lambda^{-1} \pi_0 (\lambda D(a)+ [t,a]) +\lambda^{-1}\pi_0([t, a]_\fa)+\lambda^{-1}B_\fa(\pi_0(t), \tilde D( \pi_0(a)))\tilde x\\[2mm]
&=& \pi_0 ( D(a))+\lambda^{-1}B_\fa(\pi_0(t), \pi_0(\lambda D(a)+[t, a]))\tilde x\\[2mm]
&=& \pi_0 ( D(a))+B_\fa(t, D(a))\tilde x.\\[2mm]
\end{array}
\]
For every $a$ odd, we have
\[
\pi([x^*,a]_\fg)=\pi (D(a))=\pi_0(D(a)).
\]
On the other hand, 
\[
\begin{array}{lcl}
[\pi(x^*),\pi(a)]_{\tilde \fg}&=&[\lambda^{-1}\tilde x^*+\lambda^{-1}\pi_0(t)+\nu \tilde x , \pi_0(a)]_{\tilde \fg}\\[2mm]
&=&\lambda^{-1}\tilde D(\pi_0(a)) +[\lambda^{-1}\pi_0(t), \pi_0(a)]_\fa\\[2mm]
&=&\lambda^{-1} \pi_0 (\lambda D(a)+ [t,a]) +\lambda^{-1}\pi_0([t_\pi, a]_\fa)\\[2mm]
&=&\lambda^{-1} \pi_0 (\lambda D(a))+ \lambda^{-1}\pi_0([t,a]_\fa) +\lambda^{-1}\pi_0([t, a]_\fa)\\[2mm]
&=&\pi_0 (D(a)).
\end{array}
\]
Besides, $\pi([x^*,x]_\fg)=0=
[\pi(x^*), \pi(x)]_{\tilde \fg}$. Similarly, $\pi([a,x]_{\fg})=
[\pi(a),\pi(x)]_{\tilde \fg}=0$.
Let us show that $\pi$ preserves $B_\fg$. We have $B_\fg(a,x^*)=0$ and 
\[
\begin{array}{lcl}
B_{\tilde \fg}(\pi(x^*),\pi(a))&=&B_{\tilde \fg}(\lambda^{-1}\tilde x^*+\lambda^{-1}\pi_0(t)+\nu \tilde x, \pi_0(a)+B_\fa(t,a )\tilde x)\\[2mm]
&=&\lambda^{-1}B_\fa(t,a )+B_{\fa} (\lambda^{-1}\pi_0(t),\pi_0(a))\\[2mm]
&=&2 \lambda^{-1}B_\fa(t, a)=0.
\end{array}
\]
The remaining conditions follow from the computations preceding Theorem \ref{Isom1}.

\end{proof}



%


\sssbegin{Remark} Condition (\ref{Ca}) implies that condition (\ref{Cd}) also holds on $\fa_\od$.
\end{Remark}

\sssbegin{Corollary}\label{CorAbove} Two even derivations $D$ and $D'$ that satisfy conditions (\ref{D1}) and \eqref{D3}, and are cohomologous, i.e., $[D]=[D']$ in $\mathrm{H}^1_\ev(\fa; \fa)$, define the same even double extension up to an isometry.
\end{Corollary}

\begin{proof}
Since $[D]=[D']$ in $\mathrm{H}^1_\ev(\fa; \fa)$, it follows that $D=\lambda D'+\ad_t$ for some $\lambda\in\Kee$ and some  $t\in\fa_\ev$.  We define $\pi_0=\id$ and 
\[
\begin{array}{rcl}
\alpha'&=&\lambda \alpha +B_\fa(t, s_\fa  ) \quad \text{on $\fa_\od$}; \\[2mm]
B_\fg(x^*,x^*)&=&\lambda^{-2}( B_\fa(t, t)+B_\fg(\tilde x^*, \tilde x^*)).
\end{array}
\]
The proof follows from Theorem \ref{Isom1}.
\end{proof}

\sssbegin{Remark}
The converse of the Corollary \ref{CorAbove} is not necessarily true, see \S~\ref{Exa}. 
\end{Remark}

\subsection{The $D_\od$ case}{}~{}

Using the same arguments as before, we deduce the following. For every $a \in \fa_\ev$, we have $\pi(a)=\pi_0(a)$. 
Now let $a \in \fa_\od$. There exists an element $t_\pi \in \fa_\od$ such that 
\[
\pi(a)=\pi_0(a)+B(t_\pi, a)\tilde x.
\]
Besides, $\pi(x)=\lambda \tilde x$ for some $\lambda$ in $\Kee$. Indeed, let us write $\pi(x)=\lambda \tilde x+a$, where $a\in \fa_\od$. We have (for any $b$ odd)
\[
0=B_\fg(x,b)=B_{\tilde \fg}(\pi(x), \pi(b))=B_{\tilde \fg}(\lambda \tilde x+a, \pi_0(b)+B_\fa(t_\pi, b) \tilde x)=B_\fa(a,\pi_0(b)). 
\]
We have (for any $b$ even)
\[
0=B_\fg(x,b)=B_{\tilde \fg}(\pi(x), \pi(b))=B_{\tilde \fg}(\lambda \tilde x+a, \pi_0(b))=B_\fa(a,\pi_0(b)). 
\]
Since $\pi_0$ is surjective and $B_\fa$ is nondegenerate, it follows that $a=0$.
Let us show that $\pi_0$ preserves $B_\fa$. Indeed, for any $a$ and $b$ odd we have
\[
\begin{array}{lcl}
B_\fa(a,b)&=&B_{\tilde \fg}(\pi(a), \pi(b))\\[2mm]
&=& B_{\tilde \fg}(\pi_0(a)+B_\fa(t_\pi, a)\tilde x, \pi_0(b)+B_\fa(t_\pi, b)\tilde x)\\[2mm]
&=&B_{\tilde \fg}(\pi_0(a), \pi_0(b))\\[2mm]
&=&B_{\fa}(\pi_0(a), \pi_0(b)).
\end{array}
\] 
The same arguments can be used if $a$ and $b$ are both even (or one of them is odd). 

Let us show that $\pi_0$ is an isometry on $(\fa, [\cdot ,\cdot], s_\fa, B_\fa)$. For any $a, b \in \fa_\od$, we have 
\[
\pi([a,b]_\fg)=\pi_0([a,b]_\fg)=\pi_0([a,b]_\fa). 
\]
On the other hand, 
\[
[\pi(a), \pi(b)]_{\tilde \fg}=[\pi_0(a)+B_\fa(t_\pi, a)\tilde x, \pi_0(b)+B_\fa(t_\pi, b)\tilde x]_{\tilde \fg}=[\pi_0(a), \pi_0(b)]_{\tilde \fg}=[\pi_0(a), \pi_0(b)]_{\fa}.
\]
It follows that 
\begin{equation}
\label{O54}
\pi_0([a,b]_\fg)=[\pi_0(a), \pi_0(b)]_{\tilde \fg}.
\end{equation}
For any $a, b \in \fa_\ev$, we have 
\[
\pi([a,b]_\fg)=\pi_0([a,b]_\fg)=\pi_0([a,b]_\fa). 
\]
On the other hand, 
\[
[\pi(a), \pi(b)]_{\tilde \fg}=[\pi_0(a), \pi_0(b)]_{\tilde \fg}=[\pi_0(a), \pi_0(b)]_\fa.\]

Let $a$ be even and $b$ be odd (or vice versa), we get 
\[
\pi([a,b]_\fg)=\pi([a,b]_\fa)+\lambda B_\fa(D(a),b)\tilde x=\pi_0([a,b]_\fa)+B_\fa(t_\pi,[a,b]_\fa)\tilde x+\lambda B_\fa(D(a),b)\tilde x,
\]
and 
\[
[\pi(a), \pi(b)]_{\tilde \fg}=[\pi_0(a), \pi_0(b)+B_\fa(t_\pi,b)x]_{\tilde \fg}=[\pi_0(a), \pi_0(b)]_{\tilde \fg}=[\pi_0(a), \pi_0(b)]_{\fa}+B_\fa(\tilde D(\pi_0(a)),\pi_0(b))\tilde x.
\]
It follows that 
\[
B_\fa(\pi_0^{-1}\tilde D(\pi_0(a)),b)=B_\fa([t_\pi,a],b]_\fa)+\lambda B_\fa(D(a),b).
\] 
Since $B_\fa$ is even, then  
\[
B_\fa(\pi_0^{-1}\tilde D(\pi_0(a)),c)=B_\fa([t_\pi,a],c)+\lambda B_\fa(D(a),c) \text{ for any $c$ in $\fa.$}
\] 
Therefore
\begin{equation}
\label{O55}
\pi_0^{-1}\tilde D(\pi_0(a))=[t_\pi,a]_\fa+\lambda D(a) \text{ for any $a\in\fa_\ev$}.
\end{equation}

Let us study the squaring. Let us write $\pi(x^*)$ as $\mu \tilde x^*+ a+ \nu \tilde x$ for some $a$ in $\fa_\od$. We have $\pi(s_\fg(x^*))=\pi(a_0)=\pi_0(a_0)$. On the other hand, 
\[
s_{\tilde \fg}(\pi(x^*))=s_{\tilde \fg}(\mu \tilde x^*+a+\nu \tilde x)= s_\fa(a)+\mu \tilde D(a)+\mu^{2}\tilde a_0.
\]
It follows that 
\begin{equation} \label{O57}
\tilde a_0= \mu^{-2}(\pi_0(a_0)+s_\fa(a)+\mu \tilde D(a)).
\end{equation} 
Similarly, $\pi(s_\fg(x))=0$, and $s_{\tilde \fg}(\pi(x))=s_{\tilde \fg}(\lambda \tilde x)= 0$. 

Let us compute the most general case.  We write $f=rx+c+\theta x^*$. We have
\[
\pi(s_\fg(f))= \pi(s_\fa(c)+\theta^2 a_0+\theta D(c))=\pi_0(s_\fa(c))+ \theta^2 (\pi_0(a_0))+ \theta \pi_0(D(c)).
\]
On the other hand,
\[
\begin{array}{lcl}
s_{\tilde \fg}(\pi(f))&=& s_{\tilde \fg}(r\pi(x)+\pi(c)+\theta \pi(x^*))=s_{\tilde \fg}(r\lambda \tilde x+\pi_0(c)+B_\fa(t_\pi, c)\tilde x+\theta(\mu \tilde x^*+ a+\nu \tilde x))\\[2mm]
&=&s_{\tilde \fg}((r\lambda+B_\fa(t_\pi, c)+ \theta\nu) \tilde x+\pi_0(c)+\theta a+\theta \mu \tilde x^*)\\[2mm]
&=&s_{\fa}(\pi_0(c)+\theta a)+(\theta \mu)^2\tilde a_0+\theta \mu \tilde D(\pi_0(c)+\theta a)\\[2mm]

\end{array}
\]
Therefore
\[
s_{\fa}(\pi_0(c)+\theta a)+(\theta \mu)^2\tilde a_0+\theta \mu \tilde D(\pi_0(c)+\theta a)=\pi_0(s_\fa(c))+ \theta^2 (\pi_0(a_0))+ \theta \pi_0(D(c)).
\]
The equality above follows from Eqs. (\ref{O57}) and (\ref{O58}). 

Because $\pi$ is an isometry, we have
\[
1=B_{\tilde \fg}(\pi(x^*), \pi(x))=B_{\tilde \fg}(\mu \tilde x^*+a+\nu \tilde x, \lambda \tilde x)=\lambda \mu.
\]
Therefore $\mu=\lambda^{-1}$. Besides, (if $b$ is odd)
\[
0=B_{\tilde \fg}(\pi(x^*), \pi(b))=B_{\tilde \fg}(\mu \tilde x^*+a+\nu \tilde x, \pi_0(b)+B_\fa(t_\pi,b)\tilde x)=\mu B_\fa(t_\pi,b)+B_\fa(a, \pi_0(b)).
\]
It follows that $ B_\fa(\mu t_\pi+\pi_0^{-1}a, b)=0$ and so 
\[
B_\fa(\mu t_\pi+\pi_0^{-1}a, c)=0 \text{~~for any $c\in\fa$.}
\] 
Therefore, $\mu t_\pi+\pi_0^{-1}(a)=0$ which implies that 
\begin{equation}\label{O58}
a=\mu \pi_0(t_\pi).
\end{equation} 
Finally, we get
\[
\pi(x^*)=\lambda^{-1}(\tilde x^*+\pi_0(t_\pi))+\nu \tilde x.
\]

We arrive at the following Theorem.

\sssbegin{Theorem} \label{Isom2}
Let $\pi_0$ be an isometry of $(\fa,  B_\fa)$. Let $ \lambda\in \Kee^{\times}$ and $t\in\fa_\od$ satisfy the following conditions: 
\[
\begin{array}{lcl}
\pi_0^{-1}\tilde D\pi_0&=&\lambda D+\ad_{t} \quad \text{on $\fa$}.
\end{array}
\]
 Then there exists an adapted isometry $\pi: \mathscr{K} \oplus \fa \oplus \mathscr{K}^*\rightarrow  \tilde{\mathscr{K}} \oplus \fa \oplus \tilde{\mathscr{K}^*}$, where  $\tilde{\mathscr{K}}$ and $\tilde{\mathscr{K}^*}$ are spanned by $\tilde x$ and $\tilde x^*$, respectively, given by 
\[
\begin{array}{rcl}
\pi&=& \left \{ \begin{array}{ll}
\pi_0+ B_\fa(t,\cdot) \tilde x & \text{on $\fa_\od$};\\[2mm]
\pi_0&  \text{on $\fa_\ev$};\\[2mm]
\end{array}
\right.
\\[2mm]
\pi(x)&=&\lambda \tilde x;\\[2mm]
\pi(x^*)&=&\lambda^{-1}(\tilde x^*+\pi_0(t))+\nu \tilde x, \; \text{where $\nu$ is arbitrary};\\[2mm]
\tilde a_0&=& \lambda^2 \pi_0(a_0)+s_\fa(\pi_0(t))+ \lambda \pi_0(D(t)).
\end{array}
\]
\end{Theorem}

\begin{proof} To check that $\pi$ preserves the Lie bracket, it is enough to check the conditions below.  For every even element $a$, we have 
\[
\pi([x^*,a]_\fg)=\pi (D(a))=\pi_0(D(a))+ B_\fa(t, D(a)).
\] 
On the other hand, 
\[
\begin{array}{lcl}
[\pi(x^*),\pi(a)]_{\tilde \fg}&=&[\lambda^{-1}\tilde x^*+\lambda^{-1}\pi_0(t)+\nu \tilde x , \pi_0(a)+B_\fa(t,a) \tilde x]_{\tilde \fg}\\[2mm]
&=&\lambda^{-1}\tilde D(\pi_0(a)) +[\lambda^{-1}\pi_0(t), \pi_0(a)]_\fa+B_\fa(\tilde D(\lambda^{-1}\pi_0(t)), \pi_0(a))\tilde x\\[2mm]
&=&\lambda^{-1} \pi_0 (\lambda D(a)+ [t,a]) +\lambda^{-1}\pi_0([t, a]_\fa)+\lambda^{-1}B_\fa(\tilde D(\pi_0(t)), \pi_0(a))\tilde x\\[2mm]
&=&\lambda^{-1} \pi_0 (\lambda D(a)+ [t,a]) +\lambda^{-1}\pi_0([t, a]_\fa)+\lambda^{-1}B_\fa(\pi_0(t), \tilde D( \pi_0(a)))\tilde x\\[2mm]
&=& \pi_0 ( D(a))+\lambda^{-1}B_\fa(\pi_0(t), \pi_0(\lambda D(a)+[t, a]))\tilde x\\[2mm]
&=& \pi_0 ( D(a))+B_\fa(t, D(a))\tilde x.\\[2mm]
\end{array}
\]
For every odd element $a$, we have $\pi([x^*,a]_\fg)=\pi (D(a))=\pi_0(D(a))$. On the other hand, 
\[
\begin{array}{lcl}
[\pi(x^*),\pi(a)]_{\tilde \fg}&=&[\lambda^{-1}\tilde x^*+\lambda^{-1}\pi_0(t_\pi)+\nu \tilde x , \pi_0(a)]_{\tilde \fg}\\[2mm]
&=&\lambda^{-1}\tilde D(\pi_0(a)) +[\lambda^{-1}\pi_0(t_\pi), \pi_0(a)]_\fa\\[2mm]
&=&\lambda^{-1} \pi_0 (\lambda D(a)+ [t_\pi,a]) +\lambda^{-1}\pi_0([t_\pi, a]_\fa)\\[2mm]
&=&\lambda^{-1} \pi_0 (\lambda D(a))+2 \lambda^{-1}\pi_0([t_\pi,a]_\fa) \\[2mm]
&=&\pi_0 (D(a)).
\end{array}
\]
Besides, since $x$ is central in $\fg$ and $\tilde x$ is central in $\tilde \fg$, we have $
\pi([x^*,x]_\fg)=
[\pi(x^*), \pi(x)]_{\tilde \fg}=0.
$
Similarly, $ \pi([a,x]_{\fg})=
[\pi(a),\pi(x)]_{\tilde \fg}=0$. 

Let us show that $\pi$ preserves $B_\fg$. For every element $a$ in $\fa$, we have $B_\fg(x^*,a)=0$. If $a$ is odd, then we have
\[
\begin{array}{lcl}
B_{\tilde \fg}(\pi(x^*),\pi(a))&=&B_{\tilde \fg}(\lambda^{-1}\tilde x^*+\lambda^{-1}\pi_0(t_\pi)+\nu \tilde x, \pi_0(a)+B_\fa(t_\pi,a )\tilde x)\\[2mm]
&=&\lambda^{-1}B_\fa(t_\pi,a )+B_{\fa} (\lambda^{-1}\pi_0(t_\pi),\pi_0(a))\\[2mm]
&=&0.
\end{array}
\]
If $a$ is even, then $B_{\tilde \fg}(\pi(x^*), \pi(a))=0$. 

The other conditions are certainly satisfied as shown by previous computations.

Suppose that $D(a_0)=0$. Let us show that $\tilde D(\tilde a_0)=0$. Indeed, let us apply  $\pi_0^{-1}$ to $\tilde D(\tilde a_0)$:
\[
\begin{array}{lcl}
\pi_0^{-1}\tilde D(\tilde a_0)&=&\pi_0^{-1}\left ( \tilde D (\lambda^2 \pi_0(a_0)+ s_\fa(\pi_0(t))+\lambda \pi_0D(t))\right )\\[2mm]
&=&\lambda^2( D +\ad_t)(a_0)+ (\lambda D+\ad_t) s_\fa(t)+\lambda (\lambda D+\ad_t)D(t) \\[2mm]
&=&\lambda^2 \ad_t(a_0)+ \lambda D (s_\fa(t))+\lambda^2 D^2(t)+\lambda \ad_tD(t) \\[2mm]
&=&2\lambda^2 \ad_t(a_0)+ 2\lambda \ad_t(D(t)) =0.
\end{array}
\]
Let us show that $\tilde D^2=\ad_{\tilde a_0}$. Indeed,
\[
\begin{array}{lcl}
\tilde D^2&=& \pi_0 (\lambda D+\ad_t)^2\pi_0^{-1}\\[2mm]
&=& \pi_0 (\lambda^2 D^2+\lambda D \ad_t+ \lambda \ad_t D+\ad_t^2)\pi_0^{-1}\\[2mm]
&=& \pi_0 (\lambda^2 \ad_{a_0}+\lambda D \ad_t+ \lambda \ad_t D+\ad_t^2)\pi_0^{-1}\\[2mm]
&=& \ad_{\lambda^2\pi_0(a_0)}+\ad_{\lambda \pi_0(D(t))}+2\lambda \pi_0 \circ  \ad_t \circ D \circ \pi_0^{-1}+ \ad_{\pi_0(s_\fa(t))}.\\[2mm]
&=&\ad_{\tilde a_0}.
\end{array}
\]
\end{proof}

\sssbegin{Corollary} Two odd derivations $D$ and $D'$ that satisfy conditions (\ref{2D1}), (\ref{2D2}) and (\ref{2D3}) and are cohomologous, i.e., $[D]=[D']$ in $\mathrm{H}^1_\od(\fa; \fa)$, define the same $D_\od$-extension up to an isometry.
\end{Corollary}

\begin{proof}
Since $[D]=[D']$ in $\mathrm{H}^1_\ev(\fa; \fa)$, it follows that $D=\lambda D'+\ad_t$ for some $\lambda$ in $\Kee$ and some element $t$ in $\fa_\ev$.  We define $\pi_0=\id$. 
The result follows from Theorem \ref{Isom2}.
\end{proof}

\section{The case where $B_\fg$ is odd}

\subsection{$D_\od$-extensions}

\sssbegin{Theorem}\label{MainTh3} Let $(\fa, B_\fa)$ be a NIS-Lie superalgebra in characteristic $2$ such that $B_\fa$ is odd. Let $D\in\fder_\od(\fa)$ and $a_0\in \fa_\ev$ satisfy the following conditions:
\begin{eqnarray}
\label{3D1} B_\fa(D(a),b)+ B_\fa(a,D(b))&=&0 \quad \text{for any } a,b \in \fa; \\[2mm]
\label{3D1p} B_\fa(a,D(a))&=&0 \quad \text{for any $a\in \fa_\ev$} ;\\[2mm]
\label {3D2} D^2&=&\ad_{a_0};\\[2mm]
\label{3D3} D(a_0)&=&0.
\end{eqnarray}
Let $\alpha$ be a quadratic form on $\fa_\od$ such that $B_\alpha(a,b)=B_\fa(D(a),b)$. Then there exists a NIS-Lie superalgebra structure on $\fg:=\mathscr{K} \oplus \fa \oplus \mathscr{E}$ , where $\mathscr{K}:=\Span\{x\}$ and $x$ is even, $\mathscr{E}:=\Span\{e\}$ and $e$ is odd, defined as follows. The squaring is given by (where $\mu \in \Kee$):
\[
s_\fg(a+\mu e):=s_\fa(a)+ (\mu^2m+\alpha(a))x + \mu^2 a_0+\mu D(a) \qquad \text{ for any }\; a+\mu e\in  \fg_\od.
\]
The bracket is given by:
\[
\begin{array}{l}
[x,\fg]_\fg:=0; \qquad [a,b]_\fg:=[a,b]_\fa+B_\fa(D(a),b)x\quad\text{~~for any $a,b\in \fa$};\\[2mm]
[e,a]_\fg:=D(a) \quad  \text{ for any } a \in  \fa.
\end{array}
\]
The bilinear form $B_\fg$ on $\fg$ defined by:
\begin{eqnarray*}
&{B_\fg}\vert_{\fa \times \fa}:= B_\fa, \quad B_\fg(\fa,\mathscr{K}):=0, \quad B_\fg(\fa,\mathscr{E}):=0, \\ 
&B_\fg(x,e):=1, \quad B_\fg(x,x):= B_\fg(e,e):=0,
\end{eqnarray*}
is odd, non-degenerate, invariant, and symmetric.
Therefore $(\fg, B_\fg)$ is a NIS-Lie superalgebra. 
\end{Theorem}

We call the Lie superalgebra $(\fg, B_\fg)$ constructed in Theorem \ref{MainTh3} a \textit{$D_\od$-extension} of $(\fa, B_{\fa})$  by means of $D$ and $a_0$.

\begin{proof}
The proof is similar to that of Theorem \ref{MainTh}. Let us first show that $s_\fg$ is indeed a squaring on $\fg$. Recall that since $\mathscr{E}$ is an odd vector space, then $\fg_\od=\fa_\od\oplus \mathscr{E}$. Now, let $\lambda \in\Bbb K$ and let $f=a+\mu e \in \fg_\od$; we have
\[
s_\fg(\lambda f)= ((\lambda \mu)^2m+\alpha(\lambda a))x+s_\fa(\lambda a) + (\lambda \mu)^2 a_0+\mu \lambda D(\lambda a)=\lambda^2s_\fg(f),
\]
since $s_\fa$ is a squaring on $\fa$ and $D$ is a linear map. Besides, for any $f=a+\mu e$ and $g=b+ve$ in $\fa_\od$, we see that
\begin{equation}
\label{squaring3}
\begin{array}{lcl}
s_\fg(f+g)+s_\fg(f)+s_\fg(g)&=& ((\mu+v)^2m+\alpha(a+b))x+s_\fa(a+b)+(\mu+v)^2 a_0\\[2mm]
&&+(\mu+v)D(a+b)+ (\mu^2m+\alpha(a))x+s_\fa(a)+\mu^2 a_0\\[2mm]
&&+\mu D(a)+ (v^2m+\alpha(b))x+s_\fa(b)+v^2 a_0+vD(b)\\[2mm]
&=&B_\alpha(a,b)x+[a,b]_\fa+vD(a)+\mu D(b)
\end{array}
\end{equation}
is obviously bilinear. Let us show that the bracket on $\fg$ defined above is supersymmetric. Indeed, using condition (\ref{3D1}), for any $f=rx+a+\mu e$ and $g=\tilde rx+b+\tilde \mu e$ in $\fg$, we have (if $a$ and $b$ are not both odd):
\begin{eqnarray*}
[f, g]_\fg&=&B_\fa(D(a),b)x+[a,b]_\fa+\mu D(b)+\tilde \mu D(a)\\
&=& B_\fa(D(b),a)x+[b,a]_\fa+\mu D(b)+\tilde \mu D(a)\\
&=& [g,f]_\fg.
\end{eqnarray*}
We use (\ref{squaring3}) if $a$ and $b$ are both odd.

Let us check the Jacobi identity relative to the squaring $s_\fg$. Let $f=a+\mu e\in\fg_\od$ and $g \in \fg$. If $g=x$, then we are done since $x$ is central. If $g=e$, we have
\[
\begin{array}{lcl}
[s_\fg(f),g]_\fg+[f,[f,g]_\fg]_\fg&=&[(tm+\alpha(a))x+s_\fa(a) +\mu^2 a_0+t D(a),e]_\fg +[f, D(a) ]_\fg\\[2mm]
&=& D(s_\fa(a) +\mu^2 a_0+\mu D(a))\\[2mm]
&&+B_\alpha(a,D(a))x+ [a,D(a)]_\fa+\mu D(D(a))=0,
\end{array}
\]
since $D\in \fder(\fa)$, condition (\ref{3D3}) holds and $B_\fa(D(a),D(a))=0$ since $B_\fa$ is odd. 
If $g=b\in\fa_\od$, then we have
\[
\begin{array}{lcl}
[s_\fg(f),g]_\fg&=&[(\mu m+\alpha(a))x+s_\fa(a) +\mu ^2 a_0+t D(a),b]_\fg \\[2mm]
&=&B_\fa(D(s_\fa(a) +t^2 a_0+\mu D(a)),b)x+[s_\fa(a) +\mu^2 a_0+\mu D(a),b]_\fa\\[2mm]
&=&B_\fa([Da,a]_\fa +\mu D^2(a),b)x+[s_\fa(a) ,b]_\fa+\mu^2 [a_0,b]_\fa+\mu [D(a),b]_\fa\\[2mm]
&=&[a,[a ,b]_\fa]_\fa+\mu^2 [a_0,b]_\fa+\mu [D(a),b]_\fa.
\end{array}
\]
And
\[
\begin{array}{lcl}
[f,[f,g]_\fg]_\fg&=&[f,B_\alpha(a,b)x+[a,b]_\fa+\mu D(b)]_\fg\\[2mm]
&=& \mu D([a,b]_\fa+\mu D(b))+B_\fa(D(a), [a,b]_\fa+\mu D(b))x+[a,[a,b]_\fa+\mu D(b)]_\fa\\[2mm]
&=& \mu D([a,b]_\fa)+\mu^2[a_0,b]+\mu [a,D(b)]_\fa.\\[2mm]
\end{array}
\]
It follows that $[s_\fg(f),g]_\fg+[f,[f,g]_\fg]_\fg=0$ since $s_\fa$ is a squaring, $D$ is a derivation, $B_\fa$ is $\fa$-invariant and $D^2=\ad_{a_0}$.

Now, if $g=b\in\fa_\ev$, we have (since $B_\fa(D(a),b)=0$):
\[
\begin{array}{lcl}
[s_\fg(f),g]_\fg+[f,[f,g]_\fg]_\fg&=& [(tm+\alpha(a))x+s_\fa(a) +\mu^2 a_0+\mu D(a),b]_\fg \\[2mm]
&&+ [f,[a,b]_\fa+\mu D(b)]_\fg \\ [2mm]
&=&B_\fa(D(s_\fa(a) +\mu^2 a_0+t D(a)),b)x+[s_\fa(a) +\mu^2 a_0+\mu D(a),b]_\fa \\ [2mm]
&&+B_\fa(D(a),[a,b]_\fa+\mu D(b))x\\[2mm]
&&+[a,[a,b]_\fa+\mu D(b)]_\fa +t D([a,b]_\fa+\mu D(b))\\[2mm]
&=&B_\fa([Da,a]_\fa +\mu \ad_{a_0}(a)),b)x+[s_\fa(a) +\mu^2 a_0+\mu D(a),b]_\fa \\ [2mm]
&&+B_\fa(D(a),[a,b]_\fa+\mu D(b)))x\\[2mm]
&&+[a,[a,b]_\fa]_\fa+[a,\mu D(b)]_\fa+\mu D([a,b]_\fa+\mu D(b))\\[2mm]
&=&B_\fa([Da,a]_\fa +\mu \ad_{a_0}(a)),b)x+B_\fa(D(a),[a,b]_\fa+\mu D(b))x\\[2mm]
&&+\mu^2 [a_0,b]_\fa +\mu^2 D^2(b) +[a,[a,b]_\fa]_\fa+[s_\fa(a),b]_\fa\\[2mm]
&&+\mu [a,D(b)]_\fa+\mu [D(a),b]_\fa+\mu D([a,b]_\fa)=0,
\end{array}
\]
since $D$ is a derivation, $s_\fa$ is a squaring, $B_\fa$ is $\fa$-invariant and $D^2=\ad_{a_0}$.

To check the Jacobi identity 
\[
[h,[f,g]_\fg]_\fg+[g,[h,f]_\fg]_\fg+[f,[g,h]_\fg]_\fg=0, 
\]
we proceed as follows. 

If $h=x$, the identity
\[
[x,[f,g]_\fg]_\fg+[g,[x,f]_\fg]_\fg+[f,[g,x]_\fg]_\fg=0
\]
is certainly satisfied since $x$ is central. 

If $h=e$, the identity 
\[
[e,[f,g]_\fg]_\fg+[g,[e,f]_\fg]_\fg+[f,[g,e]_\fg]_\fg=0
\]
is also satisfied for the following reasons. If either $f$ or $g$ is equal to $x$, then we are done since $x$ is central. Now if $f=e$ (or, the other way around $g=e$), then by putting $g=ux+b+ve$ we get
\[
[e,[e,g]_\fg]_\fg+[g,[e,e]_\fg]_\fg+[e,[g,e]_\fg]_\fg=2D(D(b))=0.
\]
Let us assume now that $f, g\in\fa$. We see that if $f$ is even and $g$ is odd, then 
\[
\begin{array}{lcl}
[e,[f,g]_\fg]_\fg+[g,[e,f]_\fg]_\fg+[f,[g,e]_\fg]_\fg&=&[e, B_\fa(D(f),g)x+[f,g]_\fa]_\fg+ [g, D(f)]_\fg+ [f,D(g)]_\fg\\[2mm]
&=&D([f,g]_\fa)+ [g, D(f)]_\fa+B_\fa(D(g),D(f))x\\[2mm]
&&+ [f,D(g)]_\fa +B_\fa(D(f),D(g))x=0,
\end{array}
\]
since $D\in \fder(\fa)$ and $B_\fa$ is symmetric. If $f$ and $g$ are both even, we deduce, using the fact that $B_\fa$ is symmetric, that 
\[
\begin{array}{lcl}
[e,[f,g]_\fg]_\fg+[g,[e,f]_\fg]_\fg+[f,[g,e]_\fg]_\fg&=&[e, B_\fa(D(f),g)x+[f,g]_\fa]_\fg+ [g, D(f)]_\fg+ [f,D(g)]_\fg\\[2mm]
&=&D([f,g]_\fa)+ B_\fa(D(g),D(f))x+ [g, D(f)]_\fa\\[2mm]
&&+ B_\fa(D(f),D(g))x+ [f,D(g)]_\fa=0.
\end{array}
\]
Besides, if $f, g\in\fa_\od$, then
\[
\begin{array}{lcl}
[e,[f,g]_\fg]_\fg+[g,[e,f]_\fg]_\fg+[f,[g,e]_\fg]_\fg&=&[e, B_\alpha(a,b)x+[f,g]_\fa]_\fg+ [g, D(f)]_\fg+ [f,D(g)]_\fg\\[2mm]
&=&D([f,g]_\fa)+[g,D(f)]_\fa+[f,D(g)]_\fa=0,
\end{array}
\]
since $D\in \fder(\fa)$. 

From now and on we will assume that $f,g, h \in \fa$. We distinguish several cases to check the Jacobi identity.

If $f, g$ and $h$ are even, then we have \[
\begin{array}{lcl}
[f,[g,h]_\fg]_\fg+ \circlearrowleft (f,g,h)&=&[f,B_\fa(D(g),h)x+[g,h]_\fa]_\fg + \circlearrowleft (f,g,h)\\[2mm]
&=&[f,[g,h]_\fa]_\fa+B_\fa(D(f),[g,h]_\fa)x + \circlearrowleft (f,g,h) =0,
\end{array}
\]
because the JI holds on $\fa$ and thanks to Lemma \ref{lemma1}. 
If $f$ and $g$ are both odd but $h$ is even, then we have
\[
\begin{array}{lcl}
[f,[g,h]_\fg]_\fg +[h,[f,g]_\fg]_\fg+[g,[h,f]_\fg]_\fg&=&[f,B_\fa(D(g),h)x+[g,h]_\fa]_\fg+[h,[f,g]_\fa+B_\alpha(a,b)x]_\fg \\[2mm]
&& +[g,B_\fa(D(h),f)x+[h,f]_\fa]_\fg\\[2mm]
&=& [f,[g,h]_\fa]_\fa+B_\alpha(f,[g,h]_\fa)x+[h,[f,g]_\fa]_\fa\\[2mm]
&&+B_\fa(D(h),[f,g]_\fa)x+[g,[h,f]_\fa]_\fa+B_\fa(D(g), [h,f]_\fa)x.\\[2mm]
&=&0, \\[2mm]
\end{array}
\]
because the JI holds on $\fa$ and $D$ is a derivation.

If $f$ is odd and both $g$ and $h$ are even, and using Lemma \ref{lemma1} we have
\[
\begin{array}{lcl}
[f,[g,h]_\fg]_\fg+[h,[f,g]_\fg]_\fg+[g,[h,f]_\fg]_\fg&=&[f,B_\fa(D(g),h)x+[g,h]_\fa]_\fg+ [h,B_\fa(D(f),g)x+[f,g]_\fa ]_\fg\\[2mm]
&&+[g,B_\fa(D(h),f)x+[h,f]_\fa ]_\fg\\[2mm]
&=&[f,[g,h]_\fa]_\fa+[h,[f,g]_\fa]_\fa+[g,[h,f]_\fa]_\fa ,\\[2mm]
&&+B_\fa(D(f),[g,h]_\fa)x+B_\fa(D(h),[f,g]_\fa)x\\[2mm]
&&+B_\fa(D(g),[h,f]_\fa)x=0.
\end{array}
\] 
Let us show now that $B_\fg$ is invariant. We should check that
\[
B_\fg([f,g]_\fg, h)=B_\fg (f, [g,h]_\fg) \quad \text{for any $f,g,h\in\fg$}.
\]
This is true if $f=x$ because $x$ is central and by the very definition of $B_\fg$.  If $f=e$, $g=b+ve$ and $h=c+ze$ (where $g$ and $h$ are both odd), using Eq. (\ref{squaring3}) and the fact that $B_\fa$ is odd we have
\[
\begin{array}{lcl}
B_\fg([f,g]_\fg, h)&=& B_\fg(D(b), c+ze)\\[2mm]
&=& B_\fa(D(b), c).\\[2mm]
\end{array}
\]
On the other hand, $B_\fg(f,[g, h]_\fg)=B_\fg(e, B_\alpha(b,c)x+ [b,c]_\fa+zD(b)+vD(c))= B_\alpha(b,c)$. 

If $f=e$, $g=b+ve$ is odd and $h=tx+c$ is even, we have
\[
\begin{array}{lcl}
B_\fg([f,g]_\fg, h)&=&B_\fg(D(b), h)=B_\fa(D(b), c).
\end{array}
\]
On the other hand,
\[
\begin{array}{lcl}
B_\fg(f,[g, h]_\fg)&=&B_\fg(e, B_\fa(D(b),c)x+[b,c]_\fa + vD(c))\\[2mm]
&=&B_\fa(D(b), c) .\\[2mm]
\end{array}
\]
If $f=e$, but $g=tx+b$ and $h=ux+c$ are both even, we have
$ B_\fg([f,g]_\fg, h)=B_\fa(D(b), c)$.
On the other hand, 
\[
\begin{array}{lcl}
B_\fg(f,[g, h]_\fg)&=&B_\fg(e, B_\fa(D(b),c)x+[b,c]_\fa )\\[2mm]
&=&B_\fa(D(b), c).\\[2mm]
\end{array}
\]
If $f,g,h\in\fa$, the invariance property of $B_\fg$ directly follows from that of $B_\fa$.
\end{proof}

\sssbegin{Proposition}
\label{Rec3}
Let $(\fg, B_\fg)$ be an irreducible NIS-Lie superalgebra, where $B_\fg$ is odd. Suppose that $\fz(\fg)_\ev\not=\{0\}.$ Then $(\fg, B_\fg)$ is obtained as a $D_\od$-extension of a NIS-Lie superalgebra $\fa$. 
\end{Proposition}

\begin{proof}
Let $x$ be a non-zero element in $\fz(\fg)_\ev$. The subspace $\mathscr{K}:=\Span\{x\}$ is an ideal in $(\fg,B_\fg)$ because $x$ is central in $\fg$. Moreover, $\mathscr{K}^\perp$ is also an ideal in $(\fg,B_\fg)$. Indeed, let us show first that $(\mathscr{K}^\perp)_\ev=\fg_\ev$. This is true because $x$ is orthogonal to any even element since $B_\fg$ is odd. Now let $f\in (\mathscr{K}^\perp)_\od$. It follows that $B_\fg(x,s_\fg(f))=0$ since $B_\fg$ is odd. Hence $s_\fg(f)\in \mathscr{K}^\perp$. On the other hand, if $f\in \mathscr{K}^\perp$ and $g\in \fg$, then $[f,g]_\fg \in \mathscr{K}^\perp$ because
\[
B_\fg(x,[f,g]_\fg)=B_\fg([x,f]_\fg,g)=B_\fg(0,g)=0.
\]
Since $\mathscr{K}$ is 1-dimensional, then either $\mathscr{K} \cap \mathscr{K}^\perp=\{0\}$ or $\mathscr{K}\cap \mathscr{K}^\perp=\mathscr{K}$. The first case is to be disregarded because otherwise $\fg=\mathscr{K}\oplus \mathscr{K}^\perp$ and the Lie superalgebra $\fg$ is not irreducible. Hence, $\mathscr{K}\cap \mathscr{K}^\perp=\mathscr{K}$. It follows that $\mathscr{K}\subset \mathscr{K}^\perp$ and $\dim(\mathscr{K}^\perp)=\dim(\fg)-1$. Therefore, there exists $e \in \fg_\od$ such that
\[
\fg=\mathscr{K}^\perp\oplus \mathscr{E}, \quad \text{ where  $\mathscr{E}:=\Span\{e\}$}.
\]
This $e$ can be normalized to have $B_\fg(x,e)=1$. Besides, $B_\fg(x,x)=0$ since $\mathscr{K}\cap \mathscr{K}^\perp=\mathscr{K}$.

Let us define $\fa:=(\mathscr{K} +\mathscr{E})^\perp$. We then have a decomposition $\fg=\mathscr{K} \oplus \fa \oplus \mathscr{E}$. 

Let us define a bilinear form on $\fa$ by setting: 
\[
B_\fa={B_\fg}|_{\fa\times \fa}.
\]  
The form $B_\fa$ is non-degenerate on $\fa$. Indeed, suppose there exists an $a\in\fa$ such that 
\[
B_\fa(a,f)=0\quad \text{ for any $f\in \fa$}.
\] 
But $a$ is also orthogonal to $x$ and $e$. It follows that $B_\fg(a,f)=0$ for any $f\in\fg$. Hence, $a=0$, since $B_\fg$ is nondegenerate. 

Let us show now that there exists a NIS-Lie superalgebra structure on the vector space $\fa$ for which $\fg$ is its double extension. Let $a,b \in \fa$. The bracket $[a,b]_\fg$ belongs to $\mathscr{K} \oplus \fa$ because $\fa \subset \mathscr{K}^\perp=\mathscr{K}\oplus \fa$ and the latter is an ideal. It follows that
\[
\begin{array}{lcl}
[a,b]_\fg&=&\phi(a,b)x+\psi(a,b), \\[2mm]
[a,e]_\fg&=&f(a)x+D(a), \quad \text{for any $a\in \fa_\ev$},
\end{array}
\]
where $\psi(a,b), D(a) \in \fa$, and $\phi(a,b), f(a) \in \Bbb K$. Let us show that $\phi(a,b)=B_\fg(a,D(b))$. Indeed, since $B_\fg([a,b]_\fg,e)=B_\fg(a,[b,e]_\fg)$, it follows that 
\[
B_\fg(\phi(a,b)x+\psi(a,b),e)=B_\fg(a,f(b)x+D(b)).
\] 
Using the fact that $\psi(a,b)$ is orthogonal to $e$, and $a$ is also orthogonal to $x$, and the fact that $B_\fg(x,e)=1$, we get $\phi(a,b)=B_\fg (a,D(b))$. The antisymmetry of the bracket $[\cdot , \cdot]_\fg$ implies the antisymmery of the map $\phi(a,b)$ which, in turn, implies that 
\[
B_\fg(a,D(b))=B_\fg(D(a),b) \quad \text{ for any $a,b\in\fa$.}
\]
Similarly, 
\[
0=B_\fg([e,e]_\fg,a)=B_\fg(e,[e,a]_\fg)=B_\fg(e, f(a)x+D(a))=f(a) \quad \text{ for any $a\in\fa_\ev$.}
\]
This implies that $f(a)=0$ for any $a\in \fa_\ev$.

Now let $f\in \fg_\od$. It follows that $s_\fg(f)$ is in $\fa\oplus \mathscr{K}$ since $\fg_\ev=\fa_\ev\oplus \mathscr{K}$. Let us then write $s_\fg(e)=a_0+ux$. It follows that $s_\fg(te)=t^2(a_0+ux)$. The fact that $[s_\fg(e),a]_\fg=[e,[e,a]_\fg]_\fg$ for any $a\in \fa_\ev$, implies that $B_\fg(D(a_0),a)x+[a_0,a]_\fa=D^2(a)$. Therefore, 
\[
{{D^2}_|}_{{\fa_\ev}}=\ad_{a_0}, \quad \text{ and } \quad  B_\fg(D(a_0),a)=0.
\]
Besides, $[s_\fg(e),e]_\fg=0$ implies that $D(a_0)=0$. Let us now describe a squaring on $\fa$. Since $\fa \subset \mathscr{K}^\perp$, then $s_\fg(a)\in (\mathscr{K}^\perp)_\ev=\mathscr{K} \oplus \fa_\ev$, for any $a \in \fa_\od$.  It follows then that
\[
s_\fg(a)=\alpha(a)x+s_\fa(a).
\]
Now, because $s_\fg$ is a squaring on $\fg$,  
it follows that $\alpha$ is a quadratic form on $\fa_\od$ and $s_\fa$ behaves as a squaring on $\fa$. Now, the squaring on $\fg$ take the form
\[
s_\fg(a+te)=\alpha(a,te)x+f(a,te),
\]
where $\alpha(a,te) \in \Kee$ and $f(a,te)\in  \fa$. Since 
\[
0=B_\fg([e, te]_\fg,a)=B_\fg(e, [te,a]_\fg),
\] 
it follows that 
\[
0=B_\fg(e,\alpha(a,te)x+f(a,te)+\alpha(a)x+s_\fa(a)+t^2(a_0+ux))=\alpha(a,te)+\alpha(a)+t^2u.
\] 
Therefore, $\alpha(a,te)=\alpha(a)+t^2u$. On the other hand, the fact that $B_\fg(te,[a, b]_\fg)=B_\fg([te,a]_\fg,b)$, it follows that 
\[
B_\alpha(a,b)=B_\fg( f(a,te)+s_\fa(a)+t^2a_0+tD(a), b),\quad \text{for any $a,b\in \fa_\od$}.
\]
Later on, we will show that  $B_\alpha(a,b)=B_\alpha(D(a),b)$ and that would imply that 
\[ 
f(a,te)=s_\fa(a)+t^2a_0+tD(a),
\] 
since $B_\fg$ is non-degenerate and odd. Therefore,
\begin{equation}
\label{squaring7}
s_\fg(a+te)=(\alpha(a)+t^2u)x+tD(a)+s_\fa(a).
\end{equation}

The map $\psi$ is a bilinear on $\fa$ because $[\cdot , \cdot]_\fg$ is bilinear, so for convenience let us re-denote $\psi$ by $[\cdot , \cdot]_\fa$. 

For any $a,b, c\in \fa$, we have
\[
\begin{array}{lcl}
0&=&[[a,b]_\fg,c]_\fg+ \circlearrowleft (a,b,c)\\[2mm]
&=&[[a,b]_\fa + B_\fa(D(a),b)x,c]_\fg + \circlearrowleft (a,b,c)\\[2mm]
&=&[[a,b]_\fa,c]_\fa+ B_\fg(D[a,b]_\fa,c)x+ \circlearrowleft (a,b,c). \end{array}
\]
This implies that the JI for the bracket $[\cdot,\cdot]_\fa$ is satisfied and 
\[
B_\fg(c,D([a,b]_\fa)+[a,D(b)]_\fa+[D(a),b]_\fa)=0\quad \text{ for any $a, b,c\in\fa$}.
\] Since $B_\fg$ is also non-degenerate on $\fa$, it follows that $D([a,b]_\fa)=[D(a),b]_\fa+[a,D(b)]_\fa$.  

Let us assume that $a$ and $b$ are odd and $c$ is even
\[
\begin{array}{lcl}
0&=&[[a,b]_\fg,c]_\fg+ [[c,a]_\fg,b]_\fg+[[b,c]_\fg,a]_\fg\\[2mm]
&=&[[a,b]_\fa,c]_\fg+ [[c,a]_\fa,b]_\fg+[[b,c]_\fa,a]_\fg\\[2mm]
&=&[[a,b]_\fa,c]_\fa+B_\fa(D([a,b]),c)x+ [[c,a]_\fa,b]_\fa+B_\alpha([c,a]_\fa,b)x+[[b,c]_\fa,a]_\fa\\[2mm]
&&+B_\alpha([b,c]_\fa,a)x.
\end{array}
\]
It follows that the JI is satisfied for the bracket $[\cdot , \cdot]_\fa$, provided that 
\begin{equation}
\label{cond2}
B_\fa(D([a,b]),c)+B_\alpha([c,a]_\fa,b)+B_\alpha([b,c]_\fa,a)=0.
\end{equation}
We will show later that this condition follows form another condition. 

We can similarly show that the JI identity on $\fg$ for one element odd and two elements even implies the JI on $\fa$. 

Let us complete the proof of the fact that $D^2=\ad_{a_0}$. Indeed, the fact that $[s_\fg(e),a]_\fg=[e,[e,a]_\fg]_\fg$ for any $a\in \fa_\od$, implies $B_\fg(D(a_0),a)x+[a_0,a]_\fa=[e,[e,a]_\fg]_\fg$. Using Eq. (\ref{squaring7}), we get
\[
[e,a]_\fg=s_\fg(a+e)+s_\fg(a)+s_\fg(e)=D(a).
\] 
Therefore, $D^2=\ad_{a_0}$. 

Now, since $[s_\fg(a),e]_\fg=[s_\fa(a),e]_\fg=D(s_\fg(a))$ it follows that $
[a,[a,e]_\fg]_\fg=D(s_\fa(a))$ which, in turn, implies 
\[
\begin{array}{lcl}
D(s_\fa(a))&=&B_\fg(D(a), D(a))x+[a,D(a)]_\fa\\[2mm]
&=&[a,D(a)]_\fa,
\end{array}
\]
since $B_\fg$ is odd. Hence $D(s_\fg(a))=[a,D(a)]_\fa$. The arguments used in Theorem \ref{MainTh2} show that $D([a,b]_\fa)=[D(a),b]_\fa+[a,D(b)]_\fa$, and hence $D\fder_\od(\fa)$.

Now, we should prove that the bilinear form $B_\fa:={B_\fg}|_{\fa \times \fa}$ is also invariant on $(\fa, [\cdot ,\cdot]_\fa, s_\fa, B_\fa)$. Indeed, let $a,b$ and $c$ in $\fa$ (at least two of them even), we have
\[
\begin{array}{lcl}
B_\fa(a,[b,c]_\fa)&=&B_\fg(a, [b,c]_\fg+B_\fa(D(b),c)x)\\[2mm]
&=&B_\fg(a, [b,c]_\fg)\\[2mm]
&=&B_\fg([a,b]_\fg, c)\\[2mm]
&=&B_\fg([a,b]_\fa+B_\fa(D(a),b)x, c)\\[2mm]
&=&B_\fg([a,b]_\fa, c)\\[2mm]
&=&B_\fa([a,b]_\fa, c).\\
\end{array}
\]
The proof is similar in the other cases. 

Next, we show that the bilinear form $B_\fa$ satisfies the condition 
\[
B_\fa(D(a),b)=B_\alpha(a,b) \quad \text{ for any $a$ and $b\in \fa_\od$,}
\] 
which, in turn, implies Eq. (\ref{cond2}). We have,
\[
\begin{array}{lcl}
B_\fa(D(a),b)&=&B_\fg(D(a),b) \\[2mm]
&=&B_\fg([e,a]_\fg,b)\\[2mm]
&=&B_\fg(e, [a,b]_\fg)\\[2mm]
&=&B_\fg(e, [a,b]_\fa+B_\alpha(a,b)x)\\[2mm]
&=&B_\alpha(a,b).
\end{array}
\]
The proof is complete now.
\end{proof}

\subsection{$D_\ev$-extensions}

\sssbegin{Theorem}\label{MainTh4} Let $(\fa, B_\fa)$ be a NIS-Lie superalgebra in characteristic $2$ such that $B_\fa$ is odd. Let $D\in\fder_\ev(\fa)$ be such that:
\begin{eqnarray}
\label{4D1} B_\fa(D(a),b)+B_\fa(a,D(b))&=& 0 \text{ for any } a,b \in \fa.
\end{eqnarray}
Then there exists a NIS-Lie superalgebra structure on $\fg:=\mathscr{K} \oplus \fa \oplus \mathscr{E}$ , where $\mathscr{K}:=\Span\{x\}$ and $x$ is odd, and where $\mathscr{E}:=\Span\{e\}$ and $e$ is even, defined as follows. The squaring is given by (where $\mu\in \Kee$):
\[
s_\fg(a+\mu x):= s_\fa(a)  \qquad \text{ for any }\; a+\mu x\in  \fg_\od.
\]
The bracket is given by:
\[
\begin{array}{l}
[x,\fg]_\fg:=0; \qquad [a,b]_\fg:=B_\fa(D(a),b)x+[a,b]_\fa\quad\text{~~for any $a,b\in  \fa$};\\[2mm]
[a, e]_\fg:=D(a) \quad  \text{ for any } a \in \fa.
\end{array}
\]
The bilinear form $B_\fg$ on $\fg$ defined by:
\begin{eqnarray*}
&{B_\fg}\vert_{\fa \times \fa}:= B_\fa, \quad B_\fg(\fa,\mathscr{K}):=0, \quad B_\fg(\fa,\mathscr{E}):=0, \\ 
&B_\fg(x,e):=1, \quad B_\fg(x,x):= B_\fg(e,e):=0,
\end{eqnarray*}
is odd, non-degenerate,  invariant, and symmetric.

Therefore $(\fg, B_\fg)$ is a NIS-Lie superalgebra. 
\end{Theorem}

We call the Lie superalgebra constructed in Theorem \ref{MainTh4} the \textit{$D_\ev$-extension} of $(\fa, B_{\fa})$  by means of $D$.

\begin{proof}
The proof is similar to that of Theorem \ref{MainTh3}. Let us first show that $s_\fg$ is indeed a squaring on $\fg$. Recall that since $\mathscr{K}$ is an odd vector spaces, then $\fg_\od=\fa_\od\oplus \mathscr{K}$. Now, let $\lambda\in\Bbb K$ and let $f=a+\mu x\in\fg_\od$; we have
\[
s_\fg(\lambda f)= s_\fa(\lambda a) =\lambda^2s_\fg(a)=\lambda^2s_\fg(f),
\]
since $s_\fa$ is a squaring on $\fa$. Besides, for any $f=a+\mu x$ and $g=b+vx\in\fa_\od$, we see that
\begin{equation}
\label{squaring4}
\begin{array}{lcl}
s_\fg(f+g)+s_\fg(f)+s_\fg(g)&=& s_\fa(a+b)+s_\fa(b)+s_\fa(b)
\end{array}
\end{equation}
is obviously bilinear. Let us show that the bracket on $\fg$ defined above is supersymmetric. Indeed, using condition (\ref{4D1}), for any $f=rx+a+\mu e$ and $g=\tilde rx+b+\tilde \mu e$ in $\fg$, we have (if $a$ and $b$ are not both odd):
\begin{eqnarray*}
[f, g]_\fg&=&B_\fa(D(a),b)x+[a,b]_\fa\\
&=& B_\fa(D(b),a)x+[b,a]_\fa\\
&=& [g,f]_\fg.
\end{eqnarray*}
We use (\ref{squaring4}) if $a$ and $b$ are both odd.

Let us check the Jacobi identity relative to the squaring $s_\fg$. Let $f=a+\mu x\in\fg_\od$ and $g \in \fg$. If $g=x$, then we are done since $x$ is central. If $g=e$, we have
\[
\begin{array}{lcl}
[s_\fg(f),g]_\fg+[f,[f,g]_\fg]_\fg&=&[s_\fa(a),e]_\fg +[a+\mu x, D(a) ]_\fg\\[2mm]
&=& D(s_\fa(a)) + [a,D(a)]_\fa=0,
\end{array}
\]
since $D\in \fder(\fa)$.
If $g=b\in\fa_\od$, then we have
\[
\begin{array}{lcl}
[s_\fg(f),g]_\fg&=&[s_\fa(a),b]_\fg \\[2mm]
&=&[s_\fa(a),b]_\fa +B_\fa(D(s_\fa(a)),b)x\\[2mm]
&=&[a,[a,b]_\fa]_\fa +B_\fa(D(s_\fa(a)),b)x
.\end{array}
\]
On the other hand,
\[
\begin{array}{lcl}
[f,[f,g]_\fg]_\fg&=&[a+\mu x,[a,b]_\fa]_\fg\\[2mm]
&=&[a,[a,b]_\fa]_\fg\\[2mm]
&=&[a,[a,b]_\fa]_\fa+B_\fa(D(a),[a,b]_\fa)x\\[2mm]
&=&[a,[a,b]_\fa]_\fa+B_\fa([D(a),a]_\fa,b)x\\[2mm]
\end{array}
\]
It follows that $[s_\fg(f),g]_\fg+[f,[f,g]_\fg]_\fg=0$ since $s_\fa$ is a squaring and $D$ is a derivation.
Now, if $g=b\in\fa_\ev$, we have:
\[
\begin{array}{lcl}
[s_\fg(f),g]_\fg&=&[s_\fa(a),b]_\fg \\[2mm]
&=&[s_\fa(a),b]_\fa \\[2mm]
&=&[a,[a,b]_\fa]_\fa.
\end{array}
\]
On the other hand,
\[
\begin{array}{lcl}
[f,[f,g]_\fg]_\fg&=&[a+\mu x,[a,b]_\fa+B_\fa(D(a),b)x]_\fg\\[2mm]
&=&[a,[a,b]_\fa]_\fg\\[2mm]
&=&[a,[a,b]_\fa]_\fa.
\end{array}
\]
To check the Jacobi identity 
\[
[h,[f,g]_\fg]_\fg+[g,[h,f]_\fg]_\fg+[f,[g,h]_\fg]_\fg=0, 
\]
we proceed as follows. 

If $h=x$, the identity
\[
[x,[f,g]_\fg]_\fg+[g,[x,f]_\fg]_\fg+[f,[g,x]_\fg]_\fg=0
\]
is certainly satisfied since $x$ is central. 

If $h=e$, the identity 
\[
[e,[f,g]_\fg]_\fg+[g,[e,f]_\fg]_\fg+[f,[g,e]_\fg]_\fg=0
\]
is also satisfied for the following reasons. If either $f$ or $g$ is $x$ then we are done since $x$ is central. Now if $f=e$ (or the way around $g=e$), then by putting $g=ux+b+ve$ we get
\[
[e,[e,g]_\fg]_\fg+[g,[e,e]_\fg]_\fg+[e,[g,e]_\fg]_\fg=2D(D(b))=0.
\]
Let us assume now that $f, g\in\fa$. We see that if $f$ is even and $g$ is odd, then 
\[
\begin{array}{lcl}
[e,[f,g]_\fg]_\fg+[g,[e,f]_\fg]_\fg+[f,[g,e]_\fg]_\fg&=&[e, B_\fa(D(f),g)x+[f,g]_\fa]_\fg+ [g, D(f)]_\fg+ [f,D(g)]_\fg\\[2mm]
&=&D([f,g]_\fa)+ B_\fa(D(g),D(f))x+ [g, D(f)]_\fa\\[2mm]
&&+ B_\fa(D(f),D(g))x+ [f,D(g)]_\fa=0.
\end{array}
\]
since $D$ is a derivation on $(\fa,[\cdot ,\cdot]_\fa,s_\fa)$ and $B_\fa$ is symmetric. If $f$ and $g$ are both even, we deduce, using the fact that $D$ is a derivation, that 

\[
\begin{array}{lcl}
[e,[f,g]_\fg]_\fg+[g,[e,f]_\fg]_\fg+[f,[g,e]_\fg]_\fg&=&[e, [f,g]_\fa]_\fg+ [g, D(f)]_\fg+ [f,D(g)]_\fg\\[2mm]
&=&D([f,g]_\fa)+ [g, D(f)]_\fa+ [f,D(g)]_\fa=0,
\end{array}
\]

Besides, if $f, g\in\fa_\od$, then
\[
\begin{array}{lcl}
[e,[f,g]_\fg]_\fg+[g,[e,f]_\fg]_\fg+[f,[g,e]_\fg]_\fg&=&[e, [f,g]_\fa]_\fg+ [g, D(f)]_\fg+ [f,D(g)]_\fg\\[2mm]
&=&D([f,g]_\fa)+[g,D(f)]_\fa+[f,D(g)]_\fa=0,
\end{array}
\]
since $D$ is a derivation on $(\fa,[\cdot ,\cdot]_\fa,s_\fa)$. 

From now and on we will assume that $f,g$ and $h\in \fa$. We distinguish several cases to check the Jacobi identity.

If $f, g$ and $h$ are even, then we have \[
\begin{array}{lcl}
[f,[g,h]_\fg]_\fg+ \circlearrowleft (f,g,h)&=&[f,[g,h]_\fa]_\fg + \circlearrowleft (f,g,h)\\[2mm]
&=&[f,[g,h]_\fa]_\fa + \circlearrowleft (f,g,h) =0,
\end{array}
\]
because the JI holds on $\fa$.
If $f$ and $g$ are both odd but $h$ is even, then we have
\[
\begin{array}{lcl}
[f,[g,h]_\fg]_\fg +[h,[f,g]_\fg]_\fg+[g,[h,f]_\fg]_\fg&=&[f,B_\fa(D(g),h)x+[g,h]_\fa]_\fg+[h,[f,g]_\fa]_\fg \\[2mm]
&& +[g,B_\fa(D(h),f)x+[h,f]_\fa]_\fg\\[2mm]
&=& [f,[g,h]_\fa]_\fa+[h,[f,g]_\fa]_\fa+[g,[h,f]_\fa]_\fa=0, \\[2mm]
\end{array}
\]
because the Jacobi identity holds on $\fa$.

If $f$ is odd and both $g$ and $h$ are even, we have
\[
\begin{array}{lcl}
[f,[g,h]_\fg]_\fg+[h,[f,g]_\fg]_\fg+[g,[h,f]_\fg]_\fg&=&[f,[g,h]_\fa]_\fg+ [h,B_\fa(D(f),g)x+[f,g]_\fa ]_\fg\\[2mm]
&&+[g,B_\fa(D(h),f)x+[h,f]_\fa ]_\fg\\[2mm]
&=&[f,[g,h]_\fa]_\fa+[h,[f,g]_\fa]_\fa+[g,[h,f]_\fa]_\fa ,\\[2mm]
&&+B_\fa(D(f),[g,h]_\fa)x+B_\fa(D(h),[f,g]_\fa)x\\[2mm]
&&+B_\fa(D(g),[h,f]_\fa)x=0.
\end{array}
\] 
Let us show now that $B_\fg$ is invariant. We should check that
\[
B_\fg([f,g]_\fg, h)=B_\fg (f, [g,h]_\fg) \quad \text{for any $f,g,h\in\fg$}.
\]
This is true if $f=x$ because $x$ is central and by the very definition of $B_\fg$.  If $f=e$, $g=b+vx$ and $h=c+sx$ (where $g$ and $h$ are both odd), using Eq. (\ref{squaring4}) and the fact that $B_\fa$ is odd we have
\[
\begin{array}{lcl}
B_\fg([f,g]_\fg, h)&=& B_\fg(D(b), c+sx)\\[2mm]
&=& B_\fa(D(b), c)=0.\\[2mm]
\end{array}
\]
On the other hand, $B_\fg(f,[g, h]_\fg)=B_\fg(e,  [b,c]_\fa)=0$. 

If $f=e$, $g=b+vx$ is odd and $h=c+te$ is even, we have
\[
\begin{array}{lcl}
B_\fg([f,g]_\fg, h)&=&B_\fg(D(b), h)=B_\fa(D(b), c).
\end{array}
\]
On the other hand,
\[
\begin{array}{lcl}
B_\fg(f,[g, h]_\fg)&=&B_\fg(e, B_\fa(D(b),c)x+[b,c]_\fa+tD(b))\\[2mm]
&=&B_\fa(D(b), c) .\\[2mm]
\end{array}
\]
If $f=e$, but $g=te+b$ and $h=ue+c$ are both even, we have
$ B_\fg([f,g]_\fg, h)=B_\fa(D(b), c)$.
On the other hand, 
\[
\begin{array}{lcl}
B_\fg(f,[g, h]_\fg)&=&B_\fg(e, B_\fa(D(b),c)x+[b,c]_\fa +tD(c)+uD(b))\\[2mm]
&=&B_\fa(D(b), c).\\[2mm]
\end{array}
\]
If $f,g,h\in\fa$, the invariance property of $B_\fg$ directly follows from that of $B_\fa$.
\end{proof}

\sssbegin{Proposition}\
\label{Rec4}
Let $(\fg, B_\fg)$ be an irreducible NIS-Lie superalgebra, where $B_\fg$ is odd. Suppose that $s_\fg(\fz(\fg)_\od) \cap s_\fg(\fg_\od)^\perp\not=\{0\}.$ Then $(\fg, B_\fg)$ is obtained as either a $D_\ev$-extension or a $D_\od$-extension from a NIS-Lie superalgebra $(\fa,B_\fa)$.
\end{Proposition}

\begin{proof}
Let $x$ be a non-zero element in $\fz(\fg)_\od \cap s_\fg(\fg_\od)^\perp$. Let us assume that $s_\fg(x)=0$. If not, we can use Theorem \ref{Rec3} to construct a double extension. The subspace $\mathscr{K}:=\Span\{x\}$ is an ideal in $(\fg,B_\fg)$ because $x$ is central in $\fg$. Moreover, $\mathscr{K}^\perp$ is also an ideal in $(\fg,B_\fg)$. Indeed, let us show first that $(\mathscr{K}^\perp)_\od=\fg_\od$. This is true because $x$ is orthogonal to any odd  element since $B_\fg$ is odd. Now let $f\in (\mathscr{K}^\perp)_\od$. It follows that $B_\fg(x,s_\fg(f))=0$ since $x\in s_\fg(\fg_\od)^\perp$. Hence, $s_\fg(f)\in \mathscr{K}^\perp$. On the other hand, if $f\in \mathscr{K}^\perp$ and $g\in \fg$, then $[f,g]_\fg \in \mathscr{K}^\perp$ because
\[
B_\fg(x,[f,g]_\fg)=B_\fg([x,f]_\fg,g)=B_\fg(0,g)=0.
\]
Since $\mathscr{K}$ is 1-dimensional, then either $\mathscr{K} \cap \mathscr{K}^\perp=\{0\}$ or $\mathscr{K}\cap \mathscr{K}^\perp=\mathscr{K}$. The first case is to be disregarded because otherwise $\fg=\mathscr{K}\oplus \mathscr{K}^\perp$ and the Lie superalgebra $\fg$ will not be irreducible. Hence, $\mathscr{K}\cap \mathscr{K}^\perp=\mathscr{K}$. It follows that $\mathscr{K}\subset \mathscr{K}^\perp$ and $\dim(\mathscr{K}^\perp)=\dim(\fg)-1$. Therefore, there exists $e \in \fg_\od$ such that
\[
\fg=\mathscr{K}^\perp\oplus \mathscr{E}, \quad \text{ where  $\mathscr{E}:=\Span\{e\}$}.
\]
This $e$ can be normalized to have $B_\fg(x,e)=1$. Besides, $B_\fg(x,x)=0$ since $\mathscr{K}\cap \mathscr{K}^\perp=\mathscr{K}$.

Let us define $\fa:=(\mathscr{K} +\mathscr{E})^\perp$. We then have a decomposition $\fg=\mathscr{K} \oplus \fa \oplus \mathscr{E}$. 

Let us define a bilinear form on $\fa$ by setting: 
\[
B_\fa={B_\fg}|_{\fa\times \fa}.
\]  
The form $B_\fa$ is non-degenerate on $\fa$. Indeed, suppose there exists $a\in\fa$ such that 
\[
B_\fa(a,f)=0\quad \text{ for any $f\in \fa$}.
\] 
But $a$ is also orthogonal to $x$ and $e$. It follows that $B_\fg(a,f)=0$ for any $f\in\fg$. Hence, $a=0$, since $B_\fg$ is nondegenerate. 

Let us show now that there exists a NIS-Lie superalgebra structure on the vector space $\fa$ for which $\fg$ is its double extension. Let $a,b \in \fa$. The bracket $[a,b]_\fg$ belongs to $\mathscr{K} \oplus \fa$ because $\fa \subset \mathscr{K}^\perp=\mathscr{K}\oplus \fa$ and the latter is an ideal. It follows that
\[
\begin{array}{lcl}
[a,b]_\fg&=&\phi(a,b)x+\psi(a,b), \\[2mm]
[a,e]_\fg&=&f(a)x+D(a), \quad \text{for any $a\in \fa$}
\end{array}
\]
where $\psi(a,b), D(a) \in \fa$, and $\phi(a,b), f(a) \in \Bbb K$. Let us show that $\phi(a,b)=B_\fg(a,D(b))$. Indeed, since $B_\fg([a,b]_\fg,e)=B_\fg(a,[b,e]_\fg)$ it follows that 
\[
B_\fg(\phi(a,b)x+\psi(a,b),e)=B_\fg(a,f(b)x+D(b)).
\] 
Using the fact that $\psi(a,b)$ is orthogonal to $e$, and $a$ is also orthogonal to $x$, and the fact that $B_\fg(x,e)=1$, we get $\phi(a,b)=B_\fg (a,D(b))$. The antisymmetry of the bracket $[\cdot , \cdot]_\fg$ implies the antisymmery of the map $\phi(a,b)$ which, in turn, implies that 
\[
B_\fg(a,D(b))=B_\fg(D(a),b) \quad \text{ for any $a,b\in\fa$.}
\]
Similarly, 
\[
0=B_\fg([e,e]_\fg,a)=B_\fg(e,[e,a]_\fg)=B_\fg(e, f(a)x+D(a))=f(a) \quad \text{ for any $a\in\fa$.}
\]
This implies that $f(a)=0$ for any $a\in \fa$.

Let us now describe a squaring on $\fa$. Since $\fa \subset \mathscr{K}^\perp$, then $s_\fg(a)\in (\mathscr{K}^\perp)_\ev=\mathscr{E} \oplus \fa_\ev$, for any $a \in \fa_\od$.  It follows then that
\[
s_\fg(a)=\alpha(a)e+s_\fa(a).
\]
Now, because $s_\fg(a)$ is orthogonal to $x$, it follows that $0=B_\fg(x,\alpha(a)e+s_\fa(a))=\alpha(a)$. Therefore, $s_\fg(a)=s_\fa(a).$ It follows that $s_\fa$ is a squaring on $\fa$.

The map $\psi$ is a bilinear on $\fa$ because $[\cdot , \cdot]_\fg$ is bilinear, so for convenience let us re-denote $\psi$ by $[\cdot , \cdot]_\fa$. We still have to show that the Jacobi identity is satisfied for $[\cdot , \cdot]_\fa$.
For any $a,b, c \in \fa_\ev$, we have
\[
\begin{array}{lcl}
0&=&[[a,b]_\fg,c]_\fg+ \circlearrowleft (a,b,c)\\[2mm]
&=&[[a,b]_\fa ,c]_\fg + \circlearrowleft (a,b,c)\\[2mm]
&=&[[a,b]_\fa,c]_\fa+ \circlearrowleft (a,b,c). \end{array}
\]
This implies that the JI for the bracket $[\cdot,\cdot]_\fa$ is satisfied.

Let us assume that $a$ and $b$ are odd and $c$ is even
\[
\begin{array}{lcl}
0&=&[[a,b]_\fg,c]_\fg+ [[c,a]_\fg,b]_\fg+[[b,c]_\fg,a]_\fg\\[2mm]
&=&[[a,b]_\fa,c]_\fg+ [[c,a]_\fa+B_\fa(D(c),a)x,b]_\fg+[[b,c]_\fa+B_\fa(D(b),c)x,a]_\fg\\[2mm]
&=&[[a,b]_\fa,c]_\fa+B_\fa(D([a,b]_\fa),c)x+ [[c,a]_\fa,b]_\fa+B_\fa(D([c,a]_\fa),b)x+[[b,c]_\fa,a]_\fa\\[2mm]
&&+B_\fa(D([b,c]_\fa),a)x.
\end{array}
\]
This implies that the JI for the bracket $[\cdot,\cdot]_\fa$ is satisfied and 
\[
B_\fg(c,D([a,b]_\fa)+[a,D(b)]_\fa+[D(a),b]_\fa)=0\quad \text{ for any $a, b\in\fa_\od$ and $c\in \fa_\ev$}.
\] Since $B_\fg$ is also non-degenerate on $\fa$, it follows that $D([a,b]_\fa)=[D(a),b]_\fa+[a,D(b)]_\fa$.  

Similarly, we can easily show that the JI identity on $\fg$ for one element odd and two elements even implies the JI on $\fa$, and also $D([a,b]_\fa)=[D(a),b]_\fa+[a,D(b)]_\fa$.

Now, $[s_\fg(a),e]_\fg=[s_\fa(a),e]_\fg=D(s_\fg(a))$ implies that $
[a,[a,e]_\fg]_\fg=D(s_\fa(a))$ which, in turn, implies 
\[
\begin{array}{lcl}
D(s_\fa(a))&=&B_\fg(D(a), D(a))x+[a,D(a)]_\fa\\[2mm]
&=&[a,D(a)]_\fa,
\end{array}
\]
since $B_\fg$ is odd. Hence $D(s_\fg(a))=[a,D(a)]_\fa$. Hence $D$ is an even derivation of $(\fa, [\cdot ,\cdot]_\fa, s_\fa)$.

Now, we should prove that the new defined bilinear form $B_\fa:={B_\fg}|_{\fa \times \fa}$ is also invariant on $(\fa, [\cdot ,\cdot]_\fa, s_\fa, B_\fa)$. Indeed, for any $a,b, c \in \fa_\ev$, we have
\[
\begin{array}{lcl}
B_\fa(a,[b,c]_\fa)&=&B_\fg(a, [b,c]_\fg)\\[2mm]
&=&B_\fg([a,b]_\fg, c)\\[2mm]
&=&B_\fg([a,b]_\fa, c)\\[2mm]
&=&B_\fa([a,b]_\fa, c).\\
\end{array}
\]
Now, for any $a,b$ and $c$ in $\fa$ (where $a$ and $b$ are odd and $c$ is even), we have
\[
\begin{array}{lcl}
B_\fa(a,[b,c]_\fa)&=&B_\fg(a, [b,c]_\fg+B_\fa(D(b),c)x)\\[2mm]
&=&B_\fg(a, [b,c]_\fg)\\[2mm]
&=&B_\fg([a,b]_\fg, c)\\[2mm]
&=&B_\fg([a,b]_\fa, c)\\[2mm]
&=&B_\fa([a,b]_\fa, c).
\end{array}
\]
Similarly, one can prove the invariance property when $a$, $b$ and $c$ are all odd. 

The proof is complete now.
\end{proof}

\subsection{Isometries, and equivalence classes of derivations}\label{ST2}
For a NIS-Lie superalgebra $\fa$ with a bilinear form $B_\fa$, denote by $\fg$ (resp. $\tilde \fg$) the double extension of $\fa$ by means of $D, a_0, m$ (resp. $\tilde D, \tilde a_0, \tilde m$). In the case of $D_\od$-extensions, to define $\fg$ (resp. $\tilde \fg$) we also need a quadratic form $\alpha$ (resp. $\tilde \alpha$).  We will investigate how the derivations $D$ and $\tilde D$ are related with each other when $\fg$ and $\tilde \fg$ are isometric. We will assume further that the isometry satisfies $\pi(\mathscr{K} \oplus \fa)=\tilde{\mathscr{K}} \oplus \fa$.

\subsubsection{The $D_\od$ case}

Let $\pr: \tilde{\mathscr{K}} \oplus \fa\rightarrow \fa$ be the projection, and $\pi_0:=\pr \circ \pi$. The map $\pi_0$ is obviously linear. Let $a\in \fa_\ev$. Since $\pi(a)+\pi_0(a)\in \mathrm{Ker}(\pr)$, it follows that $\pi(a)+\pi_0(a)\in  \tilde{\mathscr{K}} $. Since $B_\fa$ is nondegenerate, it follows that there exists a unique $t_\pi \in \fa_\od$ (depending only on $\pi$) such that 
\[
\pi(a)+\pi_0(a)=B_\fa(t_\pi,a) \tilde x \quad \text{for any $a\in\fa_\ev$.}
\]
On the other hand, $\pi_0(a)=\pi(a)$ for every $a\in \fa_\od$. Let us write $\pi(x)=\lambda \tilde x+u$, where $u\in \fa_\ev$. We have (for any $b$ even)
\[
0=B_\fg(x,b)=B_{\tilde \fg}(\pi(x), \pi(b))=B_{\tilde \fg}(\lambda \tilde x+u, \pi_0(b)+B_\fa(t_\pi, b) \tilde x)=B_\fa(u,\pi_0(b)). 
\]
We have (for any $b$ odd)
\[
0=B_\fg(x,b)=B_{\tilde \fg}(\pi(x), \pi(b))=B_{\tilde \fg}(\lambda \tilde x+u, \pi_0(b))=B_\fa(u,\pi_0(b)). 
\]
It follows that $u=0$ since $\pi_0$ is surjective. 
 
Let us show that $\pi_0$ preserves $B_\fa$. Indeed, for any $a, b \in \fa$ (either both even or both odd) we have
\[
B_\fa(a,b)=B_{\fa}(\pi_0(a), \pi_0(b))=0.
\] 
For any $a\in \fa_\ev$ and $b\in\fa_\od$, we have
\[
\begin{array}{lcl}
B_\fa(a,b)&=&B_{\tilde \fg}(\pi(a), \pi(b))\\[2mm]
&=& B_{\tilde \fg}(\pi_0(a)+B_\fa(t_\pi, a)\tilde x, \pi_0(b))\\[2mm]
&=&B_{\tilde \fg}(\pi_0(a), \pi_0(b))\\[2mm]
&=&B_{\fa}(\pi_0(a), \pi_0(b)).
\end{array}
\]

Let us study the squaring. Let $a\in \fa_\od$. We have 
\[
\pi(s_\fg(a))=\pi(s_\fa(a)+\alpha(a)x)=\pi_0(s_\fa(a))+B_\fa(t_\pi, s_\fa(a))\tilde x+\alpha(a)(\lambda \tilde x).
\]
On the other hand, 
\[
s_{\tilde \fg}(\pi(a))=s_{\tilde \fg}(\pi_0(a))=s_{\fa}(\pi_0(a))+\tilde \alpha(\pi_0(a))\tilde x.
\]
It follows that (for any $a$ odd)
\begin{eqnarray}
\nonumber \tilde \alpha(\pi_0(a))+\lambda \alpha(a)+B_\fa(t_\pi, s_\fa(a))&=&0, \; \text{ and }\\[2mm] 
\label{3pi01} \pi_0(s_\fa(a))&=&s_\fa(\pi_0(a)).
\end{eqnarray}

We need the following Lemma:
\sssbegin{Lemma}
If 
\[
\tilde \alpha \circ \pi_0+\lambda \alpha+B_\fa(t_\pi, s_\fa(\cdot))=0,
\]
then $\pi_0^{-1}\tilde D\pi_0+\lambda D+\ad_{t_\pi}=0.$
\end{Lemma}
\label{3Dodd}

\begin{proof} Let $a$ and $b$ in $\fa_\od$. Evaluating $\tilde \alpha \circ \pi_0+\lambda \alpha+B_\fa(t_\pi, s_\fa(\cdot))$ at $a+b$, at $a$, and at $b$, and taking the sum of evaluations we get
\[
B_{\tilde \alpha}( \pi_0(a), \pi_0(b))+\lambda B_\alpha(a,b)+B_\fa(t_\pi, [a,b]_\fa)=0.
\]
Using the fact that $B_\alpha(a,b)=B_\fa(D(a),b)$ and $B_{\tilde \alpha}(\pi_0(a),\pi_0(b))=B_\fa(\tilde D(\pi_0(a)),\pi_0(b))$ we get 
\[
B_{\fa}(\tilde D( \pi_0(a)), \pi_0(b))+\lambda B_\fa(D(a),b)+B_\fa([t_\pi,a]_\fa,b)=0.
\]
On the other hand, $B_{\fa}(\pi_0^{-1}\tilde D( \pi_0(a))+\lambda D(a)+[t_\pi,a],c)=0$, for any even element $c$, since $B_\fa$ is odd. The result follows since $B_\fa$ is nondegenerate. 
\end{proof}

Now, the fact that 
\[
\pi([a,b]_\fg)=[\pi(a), \pi(b)]_{\tilde \fg} \quad \text{ for any $a$ and $b$ in $\fa_\ev$,}
\] implies that
\[
\pi_0([a,b]_\fa)+B_\fa(D(a),b)\pi(x) + B_\fa(t_\pi, [a,b]_\fa) \tilde x= [\pi_0(a), \pi_0(b)]_\fa+B_\fa(\tilde D(\pi_0(a)), \pi_0(b))\tilde x,
\]
since $x$ is central in $\fg$. It follows that 
\begin{eqnarray}
\nonumber \lambda B(D(a),b)+B([t_\pi,a]_\fa,b)+B(\pi_0^{-1}D\pi_0(a),b)&=&0, \, \text{ and }\\[2mm]
\label{3pi02}\pi_0([a,b]_\fa)&=&[\pi_0(a), \pi_0(b)]_\fa.
\end{eqnarray}
which implies that $\pi_0^{-1}\tilde D\pi_0+\lambda D+\ad_{t_\pi}=0 $. 

For any $a$ even and $b$ odd, we have $\pi([a,b]_\fg)=\pi([a,b]_\fa)=\pi_0([a,b]_\fa)$. On the other hand, 
\[
[\pi(a), \pi(b)]_{\tilde \fg}=[\pi_0(a)+B_\fa(t_\pi, a)\tilde x, \pi_0(b)]_{\tilde \fg}=[\pi_0(a),\pi_0(b)]_{\tilde \fg}=[\pi_0(a),\pi_0(b)]_{\fa}.
\]
It follows that 
\begin{equation}
\label{3pi03}
\pi_0([a,b]_\fa)=[\pi_0(a),\pi_0(b)]_{\fa}.
\end{equation}

For any $a$ and $b$ odd, we have
\[
\begin{array}{lcl}
\pi([a,b]_\fg)=\pi([a,b]_\fa+B_\alpha(a,b) x)&=&\pi([a,b]_\fa)+\lambda B_\alpha(a,b)\tilde x\\[2mm]
&=&\pi_0([a,b]_\fa)+B_\fa(t_\pi,[a,b]_\fa)\tilde x+\lambda B_\alpha(a,b)\tilde  x\\[2mm]
&=&\pi_0([a,b]_\fa)+B_\fa(t_\pi,[a,b]_\fa)\tilde x+\lambda B_\fa(D(a),b)\tilde  x.
\end{array}
\]
On the other hand, 
\[
[\pi(a), \pi(b)]_{\tilde \fg}=[\pi_0(a), \pi_0(b)]_\fa+B_{\tilde \alpha}(\pi_0(a), \pi_0(b))\tilde x=[\pi_0(a), \pi_0(b)]_\fa+B_\fa(\tilde D(\pi_0(a)), \pi_0(b))\tilde x
\]
It follows that 
\begin{equation}
\label{3pi04}
[\pi_0(a), \pi_0(b)]_\fa=\pi([a,b]_\fa)
\end{equation}
and 
\[
B_\fa(\tilde D(\pi_0(a)), \pi_0(b))=B_\fa(t_\pi,[a,b]_\fa)+\lambda B_\fa(D(a),b).
\] 
This  condition results from Lemma \ref{3Dodd}.

Now, Eqns, (\ref{3pi01}), (\ref{3pi02}), (\ref{3pi03}) and (\ref{3pi04}) imply that $\alpha_0$ is an automorphism on $\fa$. 

Let us study the squaring. Let us write $\pi(e)$ as $\mu \tilde e+ a$ for some $a$ in $\fa_\od$. We have
\[
1=B_{\tilde \fg}(\pi(e), \pi(x))=B_{\tilde \fg}(\mu \tilde e+a, \lambda \tilde x)=\lambda \mu.
\]
Therefore $\mu=\lambda^{-1}$. Besides,
\[
\begin{array}{lcl}
0=B_\fg(e, e)=B_{\tilde \fg}(\pi(e), \pi(e))&=&B_{\tilde \fg}(\mu \tilde e+a, \mu \tilde e+a)=B_\fa(a,a).
\end{array}
\]
Therefore, $B_\fa(a, a)=0$, which is true since $B_\fa$ is odd. We have $\pi(s_\fg(e))=\pi(mx+a_0)=m\lambda \tilde x +\pi_0(a_0)+B_\fa(t_\pi,a_0)\tilde x$. On the other hand, 
\[
s_{\tilde \fg}(\pi(e))=s_{\tilde \fg}(\mu \tilde e+a)= \mu^{2}\tilde a_0+\mu^2 \tilde m \tilde x+\tilde \alpha(a)\tilde x +\mu \tilde D(a)+s_\fa(a).
\]
It follows that 
\begin{eqnarray} \label{3C57}
\tilde a_0&=& \mu^{-2}(\pi_0(a_0)+\mu \tilde D(a)+s_\fa(a)),\\[2mm]
 \label{3C58} \tilde m &=&\mu ^{-2}(\tilde \alpha(a)+\lambda m+ B_\fa(t_\pi,a_0)).
\end{eqnarray} 
Besides, for any $b\in \fa_\ev$, we have
\[
0=B_{\tilde \fg}(\pi(e), \pi(b))=B_{\tilde \fg}(\mu \tilde e+a, \pi_0(b)+B_\fa(t_\pi,b)\tilde x)=B_\fa(a, \pi_0(b))+\mu B_\fa(t_\pi,b).
\]
If follows that $B_\fa(\mu t_\pi+\pi^{-1}_0(a),b)=0$ for any $b\in \fa_\ev$. Thus, $a=\mu \pi_0(t_\pi)$ since $B_\fa$ is non-degenerate and odd. 

Let us go back to the squaring in the the most general case. We have (for $f=c+\theta e$):
\[
\begin{array}{lcl}
\pi(s_\fg(f))&=& \pi(s_\fa(c)+\theta^2 a_0+\theta D(c)+(\theta^2 m+\alpha(c))x)\\[2mm]
&=&\pi_0(s_\fa(c)+ \theta^2 a_0+\theta D(c))+B_\fa(t_\pi,s_\fa(c)+ \theta^2 a_0+\theta D(c))\tilde x+\lambda (\theta^2 m+\alpha(c))\tilde x.
\end{array}
\]
On the other hand,
\[
\begin{array}{lcl}
s_{\tilde \fg}(\pi(f))&=& s_{\tilde \fg}(\pi(c)+\theta \pi(e))=s_{\tilde \fg}(\pi_0(c)+\theta a_0 + \theta \mu \tilde e)\\[2mm]
&=&s_{\fa}(\pi_0(c)+\theta a_0)+(\theta \mu)^2\tilde a_0+\theta \mu \tilde D(\pi_0(c)+\theta a_0)+(\tilde \alpha(\pi_0(c)+\theta a_0)+\theta^2\tilde m)\tilde x.
\end{array}
\]
The fact that $s_{\tilde \fg}(\pi(f))=\pi(s_\fg(f))$, Lemma \ref{3Dodd}, Eqs. (\ref{3C57}), (\ref{3C58}) imply that $\pi_0(s_\fa(c))=s_\fa(\pi_0(c))$ for any $c\in \fa_\od$. Therefore, $\pi_0$ is an automorphism of $\fa$.

We arrive at the following Theorem.

\sssbegin{Theorem} \label{Isom3}
Let $\pi_0$ be an isometry of $(\fa,  B_\fa)$.  Let $\lambda \in \Kee^\times$ and $t\in \fa_\od$ satisfy the following conditions:
\begin{eqnarray}
\label{3Ca}\tilde \alpha&=&\lambda \alpha \circ \pi_0^{-1}+B_\fa(t, s_\fa \circ \pi_0^{-1}) \quad \text{on $\fa_\od$}; \\[2mm]
\label{3Cd}\pi_0^{-1}\tilde D\pi_0&=&\lambda D+\ad_{t} \quad \text{on $\fa_\ev$}; \\[2mm]
\label{3Ce} \tilde a_0&=&\lambda^2 \pi_0(a_0)+\lambda \pi_0(D(t))+s_\fa(\pi_0(t));\\[2mm]
\label{3Cf} \tilde m&=&\lambda^2 \alpha(t)+\lambda B_\fa(t,s_\fa(t)+\lambda a_0)+\lambda^3m.
\end{eqnarray}
Then there exists an adapted isometry $\pi: \mathscr{K} \oplus \fa \oplus \mathscr{E}\rightarrow  \tilde{\mathscr{K}} \oplus \fa \oplus \tilde{\mathscr{E}}$ given by 
\[
\begin{array}{lcl}
\pi&=& \pi_0+ B_\fa(t,\cdot) \tilde x \quad \text{on $\fa$};\\[2mm]
\pi(x)&=&\lambda \tilde x;\\[2mm]
\pi(e)&=&\lambda^{-1}(\tilde e+\pi_0(t)).\\
\end{array}
\]
\end{Theorem}

\begin{proof} To check that $\pi$ preserves the Lie bracket, it is enough to check the conditions below.  For every $a$ odd, we have
\[
\pi([e,a]_\fg)=\pi (D(a))=\pi_0(D(a))+ B_\fa(t, D(a))\tilde x.
\]
On the other hand, 
\[
\begin{array}{lcl}
[\pi(e),\pi(a)]_{\tilde \fg}&=&[\lambda^{-1}\tilde e+\lambda^{-1}\pi_0(t) , \pi_0(a)]_{\tilde \fg}\\[2mm]
&=&\lambda^{-1}\tilde D(\pi_0(a)) +[\lambda^{-1}\pi_0(t), \pi_0(a)]_\fa+B_\fa(\tilde D(\lambda^{-1}\pi_0(t)), \pi_0(a))\tilde x\\[2mm]
&=&\lambda^{-1} \pi_0 (\lambda D(a)+ [t,a]_\fa) +\lambda^{-1}\pi_0([t, a]_\fa)+\lambda^{-1}B_\fa(\tilde D(\pi_0(t)), \pi_0(a))\tilde x\\[2mm]
&=&\lambda^{-1} \pi_0 (\lambda D(a)+ [t,a]_\fa) +\lambda^{-1}\pi_0([t, a]_\fa)+\lambda^{-1}B_\fa(\pi_0(t), \tilde D( \pi_0(a)))\tilde x\\[2mm]
&=& \pi_0 ( D(a))+\lambda^{-1}B_\fa(\pi_0(t), \pi_0(\lambda D(a)+[t, a]_\fa))\tilde x\\[2mm]
&=& \pi_0 ( D(a))+B_\fa(t, D(a))\tilde x.\\[2mm]
\end{array}
\]
For every $a$ even, we have
\[
\pi([e,a]_\fg)=\pi (D(a))=\pi_0(D(a)).
\]
On the other hand, 
\[
\begin{array}{lcl}
[\pi(e),\pi(a)]_{\tilde \fg}&=&[\lambda^{-1}\tilde e+\lambda^{-1}\pi_0(t) , \pi_0(a)+B_\fa(t,a) \tilde x]_{\tilde \fg}\\[2mm]
&=&\lambda^{-1}\tilde D(\pi_0(a)) +[\lambda^{-1}\pi_0(t), \pi_0(a)]_\fa\\[2mm]
&=&\lambda^{-1} \pi_0 (\lambda D(a)+ [t,a]_\fa) +\lambda^{-1}\pi_0([t, a]_\fa)\\[2mm]
&=&\lambda^{-1} \pi_0 (\lambda D(a))+ \lambda^{-1}\pi_0([t,a]_\fa) +\lambda^{-1}\pi_0([t_\pi, a]_\fa)\\[2mm]
&=&\pi_0 (D(a)).
\end{array}
\]
Besides, $0=\pi([e,x]_\fg)=[\pi(e), \pi(x)]_{\tilde \fg}$. Similarly, $\pi([a,x]_{\fg})=
[\pi(a),\pi(x)]_{\tilde \fg}=0$.
Let us show that $\pi$ preserves $B_\fg$. We have $B_\fg(a,e)=0$ and 
\[
\begin{array}{lcl}
B_{\tilde \fg}(\pi(e),\pi(a))&=&B_{\tilde \fg}(\lambda^{-1}\tilde e+\lambda^{-1}\pi_0(t), \pi_0(a)+B_\fa(t,a )\tilde x)\\[2mm]
&=&\lambda^{-1}B_\fa(t,a )+B_{\fa} (\lambda^{-1}\pi_0(t),\pi_0(a))\\[2mm]
&=&2 \lambda^{-1}B_\fa(t, a)=0.
\end{array}
\]
The other conditions are certainly satisfied because of the previous computations.

Let us show now that $B_\fa(\tilde D(a), a)=0$ for any $a$ odd. Indeed, 
\[
\begin{array}{lcl}
B_\fa(\tilde D(a), a)&=&B_\fa(\pi_0(\lambda D+\ad_t)\pi_0^{-1}(a), a)\\[2mm]
&=&B_\fa((\lambda D+\ad_t)\pi_0^{-1}(a), \pi_0^{-1}(a))\\[2mm]
&=&B_\fa((\ad_t)\pi_0^{-1}(a), \pi_0^{-1}(a))\\[2mm]
&=&B_\fa([t,\pi_0^{-1}(a)]_\fa, \pi_0^{-1}(a))\\[2mm]
&=&B_\fa(t,[\pi_0^{-1}(a),\pi_0^{-1}(a)]_\fa)=0.
\end{array}
\]
Suppose that $D(a_0)=0$. Let us show that $\tilde D(\tilde a_0)=0$. Indeed, let us apply  $\pi_0^{-1}$ to $\tilde D(\tilde a_0)$:
\[
\begin{array}{lcl}
\pi_0^{-1}\tilde D(\tilde a_0)&=&\pi_0^{-1}\left ( \tilde D (\lambda^2 \pi_0(a_0)+\lambda \pi_0(D(t))+s_\fa(\pi_0(t))\right )\\[2mm]
&=&\pi_0^{-1} \tilde D \pi_0\left ( \lambda^2 a_0+\lambda D(t)+s_\fa(t)\right ) \\[2mm]
&=&(\lambda D+\ad_{t})\left ( \lambda^2 a_0+\lambda D(t)+s_\fa(t)\right )  \\[2mm]
&=&2\lambda^2 \ad_t(a_0)+ 2\lambda [t,D(t)]_\fa =0.
\end{array}
\]
Let us show that $\tilde D^2=\ad_{\tilde a_0}$. Indeed,
\[
\begin{array}{lcl}
\tilde D^2&=& \pi_0 (\lambda D+\ad_t)^2\pi_0^{-1}\\[2mm]
&=& \pi_0 (\lambda^2 D^2+\lambda D \circ \ad_t+ \lambda \ad_t \circ D+\ad_t^2)\pi_0^{-1}\\[2mm]
&=& \pi_0 (\lambda^2 \ad_{a_0}+\lambda D \circ \ad_t+ \lambda \ad_t \circ D+\ad_t^2)\pi_0^{-1}\\[2mm]
&=& \ad_{\lambda^2\pi_0(a_0)}+\ad_{\lambda \pi_0(D(t))}+2\lambda \pi_0 \circ  \ad_t \circ D \circ \pi_0^{-1}+ \ad_{\pi_0(s_\fa(t))}.\\[2mm]
&=&\ad_{\tilde a_0}.
\end{array}
\]

\end{proof}

\sssbegin{Remark} Condition (\ref{3Ca}) implies that condition (\ref{3Cd}) also holds on $\fa_\od$.
\end{Remark}

\sssbegin{Corollary}\label{3CorAbove} Two odd derivations $D$ and $D'$ that satisfy conditions (\ref{D1}) and  (\ref{D3}) and are cohomologous, i.e., $[D]=[D']$ in $\mathrm{H}^1_\od(\fa; \fa)$, define the same even double extension up to an isometry.
\end{Corollary}

\begin{proof}
Since $[D]=[D']$ in $\mathrm{H}^1_\od(\fa; \fa)$, it follows that $D=\lambda D'+\ad_t$ for some $\lambda$ in $\Kee$ and some element $t$ in $\fa_\ev$.  We define $\pi_0=\id$. The rest of the proof follows from Theorem \ref{Isom3}.
\end{proof}

\subsubsection{The $D_\ev$ case}

Let $\pr: \tilde{\mathscr{K}} \oplus \fa\rightarrow \fa$ be the projection, and $\pi_0:=\pr \circ \pi$. The map $\pi_0$ is obviously linear. Let $a\in \fa_\od$. Since $\pi(a)+\pi_0(a) \in \mathrm{Ker}(\pr)$, it follows that $\pi(a)+\pi_0(a)\in  \tilde{\mathscr{K}} $. Since $B_\fa$ is nondegenerate, it follows that there exists a unique $t_\pi$ in $\fa_\ev$ (depending only on $\pi$) such that 
\[
\pi(a)+\pi_0(a)=B_\fa(t_\pi,a) \tilde x \quad \text{for any $a$ in $\fa_\od$.}
\]
On the other hand, $\pi_0(a)=\pi(a)$ for every $a\in \fa_\ev$. Let us write $\pi(x)=\lambda \tilde x+u$, where $u\in \fa_\od$. We have (for any $b$ odd)
\[
0=B_\fg(x,b)=B_{\tilde \fg}(\pi(x), \pi(b))=B_{\tilde \fg}(\lambda \tilde x+u, \pi_0(b)+B_\fa(t_\pi, b) \tilde x)=B_\fa(u,\pi_0(b)). 
\]
We have (for any $b$ even)
\[
0=B_\fg(x,b)=B_{\tilde \fg}(\pi(x), \pi(b))=B_{\tilde \fg}(\lambda \tilde x+u, \pi_0(b))=B_\fa(u,\pi_0(b)). 
\]
It follows that $u=0$ since $\pi_0$ is surjective. 
 
Let us show that $\pi_0$ preserves $B_\fa$. Indeed, for any $a, b \in\fa$ (either both even or both odd) we have
\[
B_\fa(a,b)=B_{\fa}(\pi_0(a), \pi_0(b))=0.
\] 
For any $a\in \fa_\od$ and $b\in\fa_\ev$, we have
\[
\begin{array}{lcl}
B_\fa(a,b)&=&B_{\tilde \fg}(\pi(a), \pi(b))\\[2mm]
&=& B_{\tilde \fg}(\pi_0(a)+B_\fa(t_\pi, a)\tilde x, \pi_0(b))\\[2mm]
&=&B_{\tilde \fg}(\pi_0(a), \pi_0(b))\\[2mm]
&=&B_{\fa}(\pi_0(a), \pi_0(b)).
\end{array}
\]

Let us study the squaring. Let $a\in \fa_\od$. We have 
\[
\pi(s_\fg(a))=\pi(s_\fa(a))=\pi_0(s_\fa(a)).
\]
On the other hand, 
\[
s_{\tilde \fg}(\pi(a))=s_{\tilde \fg}(\pi_0(a)+B_\fa(t_\pi, a)\tilde x)=s_{\fa}(\pi_0(a)).
\]
It follows that (for any $a$ odd) 
\begin{equation}
\label{4pi01} \pi_0(s_\fa(a))=s_\fa(\pi_0(a)).
\end{equation}

Now, the fact that 
\[
\pi([a,b]_\fg)=[\pi(a), \pi(b)]_{\tilde \fg} \quad \text{ for any $a, b\in\fa_\ev$,}
\] implies that
\[
\pi_0([a,b]_\fa)= [\pi_0(a), \pi_0(b)]_\fa.
\]
 
For any $a$ even and $b$ odd, we have 
\[
\pi([a,b]_\fg)=\pi([a,b]_\fa+B_\fa(D(a),b)\tilde x)=\pi_0([a,b]_\fa)+B_\fa(t_\pi,[a,b]_\fa) \tilde x+\lambda B_\fa(D(a),b)\tilde x.
\] 
On the other hand, 
\[
[\pi(a), \pi(b)]_{\tilde \fg}=[\pi_0(a), \pi_0(b)+B_\fa(t_\pi, b)\tilde x]_{\tilde \fg}=[\pi_0(a),\pi_0(b)]_{\tilde \fg}=[\pi_0(a),\pi_0(b)]_{\fa}+B_\fa(\tilde D(\pi_0(a)),\pi_0(b))\tilde x.
\]
It follows that 
\begin{equation}
\label{4pi03}
\pi_0([a,b]_\fa)=[\pi_0(a),\pi_0(b)]_{\fa},
\end{equation}
and 
\begin{equation}
\label{4pi04}
B_\fa(\tilde D(\pi_0(a)), \pi_0(b))=B_\fa(t_\pi,[a,b]_\fa)+\lambda B_\fa(D(a),b).
\end{equation}
For any $a$ and $b$ odd, we have
\[
\begin{array}{lcl}
\pi([a,b]_\fg)=\pi([a,b]_\fa)&=&\pi_0([a,b]_\fa).
\end{array}
\]
On the other hand, 
\[
[\pi(a), \pi(b)]_{\tilde \fg}=[\pi_0(a), \pi_0(b)]_\fa.
\]
It follows that 
\begin{equation}
\label{4pi05}
[\pi_0(a), \pi_0(b)]_\fa=\pi_0([a,b]_\fa)
\end{equation}
Now, Eqs, (\ref{4pi01}),  (\ref{4pi03}) and (\ref{4pi05}) imply that $\alpha_0$ is an automorphism of $\fa$. 

Let us write $\pi(e)$ as $\mu \tilde e+ a$ for some $a\in\fa_\ev$. We have
\[
1=B_{\tilde \fg}(\pi(e), \pi(x))=B_{\tilde \fg}(\mu \tilde e+a, \lambda \tilde x)=\lambda \mu.
\]
Therefore $\mu=\lambda^{-1}$. Besides,
\[
\begin{array}{lcl}
0=B_\fg(e, e)=B_{\tilde \fg}(\pi(e), \pi(e))&=&B_{\tilde \fg}(\mu \tilde e+a, \mu \tilde e+a)=B_\fa(a,a).
\end{array}
\]
Therefore, $B_\fa(a, a)=0$, which is true since $B_\fa$ is odd. 
Besides, for any $b\in \fa_\od$, we have
\[
0=B_{\tilde \fg}(\pi(e), \pi(b))=B_{\tilde \fg}(\mu \tilde e+a, \pi_0(b)+B_\fa(t_\pi,b)\tilde x)=B_\fa(a, \pi_0(b))+\mu B_\fa(t_\pi,b).
\]
If follows that $B_\fa(\mu t_\pi+\pi^{-1}_0(a),b)=0$ for any $b\in \fa_\od$. Thus, $a=\mu \pi_0(t_\pi)$ since $B_\fa$ is non-degenerate and odd. 

We arrive at the following Theorem.

\sssbegin{Theorem} \label{Isom4}
Let $\pi_0$ be an isometry of $(\fa, B_\fa)$.  Let $\lambda\in \Kee$ and $t\in \fa_\ev$ satisfy the following conditions:
\begin{eqnarray}
\label{4Cd}\pi_0^{-1}\tilde D\pi_0&=&\lambda D+\ad_{t} \quad \text{on $\fa$}.
\end{eqnarray}
Then there exists an adapted isometry $\pi: \mathscr{K} \oplus \fa \oplus \mathscr{E}\rightarrow  \tilde{\mathscr{K}} \oplus \fa \oplus \tilde{\mathscr{E}}$ given by 
\[
\begin{array}{lcl}
\pi&=& \pi_0+ B_\fa(t,\cdot) \tilde x \quad \text{on $\fa$};\\[2mm]
\pi(x)&=&\lambda \tilde x;\\[2mm]
\pi(e)&=&\lambda^{-1}(\tilde e+\pi_0(t)).\\
\end{array}
\]
\end{Theorem}

\begin{proof} To check that $\pi$ preserves the Lie bracket, it is enough to check the conditions below.  For every $a$ odd, we have
\[
\pi([e,a]_\fg)=\pi (D(a))=\pi_0(D(a))+ B_\fa(t, D(a))\tilde x.
\]
On the other hand, 
\[
\begin{array}{lcl}
[\pi(e),\pi(a)]_{\tilde \fg}&=&[\lambda^{-1}\tilde e+\lambda^{-1}\pi_0(t) , \pi_0(a)+B(t,a)\tilde x]_{\tilde \fg}\\[2mm]
&=&\lambda^{-1}\tilde D(\pi_0(a)) +[\lambda^{-1}\pi_0(t), \pi_0(a)]_\fa+B_\fa(\tilde D(\lambda^{-1}\pi_0(t)), \pi_0(a))\tilde x\\[2mm]
&=&\lambda^{-1} \pi_0 (\lambda D(a)+ [t,a]_\fa) +\lambda^{-1}\pi_0([t, a]_\fa)+\lambda^{-1}B_\fa(\tilde D(\pi_0(t)), \pi_0(a))\tilde x\\[2mm]
&=&\lambda^{-1} \pi_0 (\lambda D(a)+ [t,a]_\fa) +\lambda^{-1}\pi_0([t, a]_\fa)+\lambda^{-1}B_\fa(\pi_0(t), \tilde D( \pi_0(a)))\tilde x\\[2mm]
&=& \pi_0 ( D(a))+\lambda^{-1}B_\fa(\pi_0(t), \pi_0(\lambda D(a)+[t, a]_\fa))\tilde x\\[2mm]
&=& \pi_0 ( D(a))+B_\fa(t, D(a))\tilde x.\\[2mm]
\end{array}
\]
For every $a$ even, we have
\[
\pi([e,a]_\fg)=\pi (D(a))=\pi_0(D(a)).
\]
On the other hand, 
\[
\begin{array}{lcl}
[\pi(e),\pi(a)]_{\tilde \fg}&=&[\lambda^{-1}\tilde e+\lambda^{-1}\pi_0(t) , \pi_0(a)+B_\fa(t,a) \tilde x]_{\tilde \fg}\\[2mm]
&=&\lambda^{-1}\tilde D(\pi_0(a)) +[\lambda^{-1}\pi_0(t), \pi_0(a)]_\fa\\[2mm]
&=&\lambda^{-1} \pi_0 (\lambda D(a)+ [t,a]) +\lambda^{-1}\pi_0([t, a]_\fa)\\[2mm]
&=&\lambda^{-1} \pi_0 (\lambda D(a))+ \lambda^{-1}\pi_0([t,a]_\fa) +\lambda^{-1}\pi_0([t, a]_\fa)\\[2mm]
&=&\pi_0 (D(a)).
\end{array}
\]
Besides, $0=\pi([e,x]_\fg)=[\pi(e), \pi(x)]_{\tilde \fg}$. Similarly, $\pi([a,x]_{\fg})=
[\pi(a),\pi(x)]_{\tilde \fg}=0$.
Let us show that $\pi$ preserves $B_\fg$. We have $B_\fg(e,a)=0$ and 
\[
\begin{array}{lcl}
B_{\tilde \fg}(\pi(e),\pi(a))&=&B_{\tilde \fg}(\lambda^{-1}\tilde e+\lambda^{-1}\pi_0(t), \pi_0(a)+B_\fa(t,a )\tilde x)\\[2mm]
&=&\lambda^{-1}B_\fa(t,a )+B_{\fa} (\lambda^{-1}\pi_0(t),\pi_0(a))\\[2mm]
&=&2 \lambda^{-1}B_\fa(t, a)=0.
\end{array}
\]
The other conditions are certainly satisfied because of the previous computations.
\end{proof}

\sssbegin{Corollary}\label{4CorAbove} Two even derivations $D$ and $D'$ that are cohomologous, i.e., $[D]=[D']$ in $\mathrm{H}^1_\ev(\fa; \fa)$, define the same double extension up to an isometry.
\end{Corollary}

\begin{proof}
Since $[D]=[D']$ in $\mathrm{H}^1_\ev(\fa; \fa)$, it follows that $D=\lambda D'+\ad_t$ for some $\lambda\in\Kee$ and some element $t$ in $\fa_\ev$.  We define $\pi_0=\id$. The rest of the proof follows from Theorem \ref{Isom4}.
\end{proof}

\section{Examples: the case where the bilinear form $B$ is even} \label{Exa}
We give here a few examples of $D$-extensions in the case where the bilinear form $B$ is even and the derivation $D$ is either even and odd. We use non-equivalent derivations to obtain non-isometric NIS-Lie superalgebras. 

\subsection{The Manin triple constructed from the Heisenberg superalgebra} Consider the Heisenberg superalgebra $\fhei(0|2)$ spanned by  $p,q$ (odd) and $z$ (even), with the only nonzero bracket: $[p,q]=z$. We consider the NIS-Lie superalgebra $\fa:=\fhei(0|2)\oplus \fhei(0|2)^*$ constructed as in \S~\ref{SecDef}.  In this case, $\fa_\od=\text{Span } \{p,q,p^*,q^*\}.$ A direct computation using Eq. (\ref{squaringh}) shows that (for any $s,w,u,v\in \Kee$)
\[
\begin{array}{lcl}
s_\fa(s p+wq+up^*+v q^*)&=&sw z+s(up^*+v q^*)\circ \ad_{p}+w(uq^*+v q^*)\circ \ad_{q}=swz.\\[2mm]
\end{array}
\]
A direct computation using Eqs. (\ref{squaringh}) and (\ref{bracketh}) shows that the only nonzero brackets are 
\[
[p,q]_\fa = z,\quad [p,z^*]_\fa = q^*,\quad [q, z^*]_\fa=p^*.
\]

\sssbegin{Claim}
\label{der1}
The space $\mathfrak{out}(\fa)$ is spanned by the (classes of the) following cocycles 
(odd cocycles are underlined):
\[
\begin{array}{lcllcllcllcl}
D_1&=&q^*\otimes  \widehat{p}, &D_2&=&q^*\otimes  \widehat{q}, &\underline{D_3}&=&q^*\otimes  \widehat{z^*},\\[2mm] 
D_4&=&p^*\otimes  \widehat{p},&\underline{D_5}&=&p^*\otimes  \widehat{z^*}, &\underline{D_6}&=&q\otimes  \widehat{z^*}+z\otimes  \widehat{q^*},\\[2mm]
\underline{D_7}&=&p\otimes  \widehat{z^*}+z\otimes  \widehat{p^*},& D_8&=&z\otimes  \widehat{z^*},&D_9&=&p\otimes  \widehat{p}+q^*\otimes  \widehat{q^*}+z\otimes  \widehat{z}, \\[2mm]
D_{10}&=&q\otimes  \widehat{q}+p^*\otimes  \widehat{p^*}+z\otimes  \widehat{z},\\[2mm]
 D_{11}&=&q^*\otimes  \widehat{q^*}+p^*\otimes  \widehat{p^*}+z^*\otimes  \widehat{z^*}.
\end{array}
\]
\end{Claim}


Let us fix an ordered basis as follows: $p,q,z, p^*, q^*, z^*$. In this basis, the Gram matrix of the bilinear form $B_\fa$ in (\ref{Mform}) is given by (where $I_n$ denotes the identity $n\times n$-matrix)
\[
\left (
\begin{array}{c|c}
0&I_3\\
\hline 
I_3& 0
\end{array}
\right ).
\]
Any derivation $D$ has the following supermatrix representation:
\[
D=\left (
\begin{array}{c|c}
A&B\\
\hline
C& F
\end{array}
\right ).
\] It follows that $D$ is compatible with the bilinear form $B_\fa$ if and only if $F=A^t, B^t=B$ and $C^t=C$. Let us consider the most general derivation $D=\mathop{\sum}\limits_{1 \leq i\leq 11}\alpha_i D_i$, where $D_i$ are the cocycles given in Claim~\ref{der1}. In the same basis $p,q,z,p^*,q^*,z^*$, we have
\[
D=\left (
\begin{array}{c|c}
\alpha_9 (E^{1,1}+E^{3,3})+\alpha_{10}(E^{2,2}+E^{3,3})&\alpha_{7} (E^{1,3}+E^{3,1})+\alpha_6(E^{2,3}+E^{3,2})+\alpha_8 E^{3,3}\\
\hline
\alpha_4 E^{1,1}+\alpha_1 E^{2,1}+\alpha_2 E^{2,2}&\alpha_{11} I_3+\alpha_{10}E^{1,1}+\alpha_{9} E^{2,2}+\alpha_5 E^{1,3}+\alpha_3 E^{2,3}
\end{array}
\right ),
\]
where $E^{i,j}$ is the $3\times 3$ matrix $(i,j)$th unit. 

It follows that $D$ is compatible with the bilinear form $B_\fa$ if and only if $\alpha_1=\alpha_3=\alpha_5=0$, and $\alpha_{11}=\alpha_9+\alpha_{10}$.

\subsubsection{The $D_\ev$-extension}



In this case, the most general even derivation is of the form $\alpha_2D_2+\alpha_4 D_4+\alpha_8 D_8+\alpha_9D_9+\alpha_{10}D_{10}+(\alpha_9+\alpha_{10})D_{11}$. 

Now, we define a quadratic form on $\fa_\od$ as follows: 
\[
\alpha(\mu p+ \nu q+\lambda p^*+\beta q^*)=\mu \lambda+\nu \beta+A(\nu^2+\beta^2),
\] 
where $A\in \Kee$, see \S \ref{Arf}. 

Let us check the condition (\ref{D3}).  Let $a=sz+tz^*$ be an even element. The fact that $B_\fa(a,D(a))=0$ implies that $\alpha_8=0$. Now, let 
\[
a=sp+tp^*+wq+uq^*,\; b=\tilde sp+\tilde tp^*+\tilde wq+\tilde uq^* \in \fa_\od.
\]
We have 

\[
\begin{array}{lcl}
B_\fa(a, D(b))&=&B_\fa(sp+tp^*+wq+uq^*,\alpha_9 \tilde sp+(\alpha_{10}\tilde t+ \alpha_4 \tilde s)p^*+\alpha_{10}\tilde wq+(\alpha_2\tilde w+\alpha_{10}\tilde u) q^*)\\[2mm]
&=&s(\alpha_{10}\tilde t+ \alpha_4 \tilde s)+t\alpha_9 \tilde s+w (\alpha_2\tilde w+\alpha_{10}\tilde u) +u \alpha_{10}\tilde w.
\end{array}
\]
On the other hand,  
\[
\begin{array}{lcl}
B_\alpha(a,b)&=&(s+\tilde s)(t+\tilde t)+(w+\tilde w)(u +\tilde u)+A((w+\tilde w)^2+(u+\tilde u)^2)\\
&&+st +wu+A (w^2+u^2)+\tilde s \tilde t +\tilde w \tilde u+A (\tilde w^2+\tilde u^2)\\[2mm]
&=& s \tilde t+ \tilde s t+w \tilde u+ \tilde w u.
\end{array}
\]
It follows that $\alpha_{10}=\alpha_9=1$ and $\alpha_4=\alpha_2=0 $. 


The even double extension of $\fa$ is then a NIS-Lie superalgebra of $\sdim=4|4$, by means of the derivation $D=p\otimes \hat p+q^*\otimes \hat q^*+q\otimes \hat q+p^*\otimes \hat p^*$ and the quadratic form $\alpha$. 
\subsubsection{The $D_\od$-extension} The only odd derivations compatible with the bilinear form $B_\fa$ are $D_6$ and $D_7$. Put $D:=\alpha_6 D_6+\alpha_7 D_7$. Now, the only even element $a_0$ for which $D(a_0)=0$ is of the form $a_0=k z$, where $k\in \Kee$. A direct computation shows that $D^2=\ad_{k z}=0$. The $D_\od$-extension of $\fa$ is then a NIS-Lie superalgebra of $\sdim=2|6$, by means of the derivation $D$ and $a_0=k z$, where $k\in\Kee$. Let us show that there is an isometry between the extensions corresponding to pairs $(\tilde D=\alpha_6D_6+\alpha_7D_7, \tilde a_0=k z)$ and $(D_6, a_0=0)$. Indeed, the isometry is given by (for notation, see Theorem \ref{Isom2})
\[
\begin{array}{rclrclrclrcl}
\pi_0(z)&=& z, & \pi_0(z^*)&=&z^*, &  \pi_0(p)&=&s_1p+s_3q, \\[2mm]
 \pi_0(q)&=&\alpha_7 p+\alpha_6q, & \pi_0(q^*)&=&s_1q^*+s_3p^*,& \pi_0(p^*)&=&\alpha_7 q^*+\alpha_6 p^*,\\[2mm]
 t&=&kq^*, & \lambda&=&1,& \nu&=&0,
\end{array}
\]
where $s_1 \alpha_6+s_3 \alpha_7=1$.

On the other hand, let us show that the $D_\od$-extension by means of $D_6$ is not a trivial one; namely, it is not isometric to the one by means of $\ad_T$ for some $T\in \fa$. Suppose there is an isometry, say $\pi$. Let us write 
\[
\pi_0(z)=mz, \quad \pi_0(z^*)=m_1z+m^{-1}z^*.
\]
Now, because $q^*=[z^*,p]$, it follows that 
\[\pi_0(q^*)=[m^{-1}z^*, \pi_0(p)]=c_1p^*+c_2q^*\quad  \text{for some $c_1, c_2 \in\Kee$.}
\] 
Similarly, since $p^*=[z^*,q]$, it follows that 
\[
\pi_0(p^*)=[m^{-1}z^*, \pi_0(q)]=\tilde c_1p^*+\tilde c_2q^*\quad  \text{for some $\tilde c_1, \tilde c_2\in \Kee$}.
\] 
We have (here $T=W_1p+W_2p^*+W_3q+W_4q^*$, where $W_i \in \Kee$):
\[
\begin{array}{lcl}
(D_6\circ \pi_0+\pi_0\circ \ad_T)(z^*)&=&D_6(m_1z+m^{-1}z^*)+\pi_0 [T, z^*]\\[2mm]
&=&m^{-1}q+ W_1\pi_0(q^*)+W_3\pi_0(p^*)\\[2mm]
&=&m^{-1}q+ W_1(c_1p^*+c_2q^*)+W_3 (\tilde c_1p^*+\tilde c_2q^*).
\end{array}
\] 
But this is never zero, hence a contradiction.
\subsection{The Manin triple constructed from the Lie superalgebra  $\mathfrak{ba}(1)$} The Lie superalgebra
$\mathfrak{ba}(n)$ is an analog of the Heisenberg Lie superalgebra $\fhei(2n|m)$: the latter is the negative part (in the standard $\Zee$-grading) of the Poisson superalgebra $\fpo(2n|m)$, the former is the negative part (in the standard $\Zee$-grading)  of the antibracket Lie superalgebra $\fb(n)$, see \cite{BGLLS2}. 

Consider the superalgebra $\mathfrak{ba}(1)$ spanned by  $\theta,z$ (odd) and $q$ (even), with the only nonzero bracket: $[q,\theta]=z$. We consider the NIS-Lie superalgebra $\fa:=\mathfrak{ba}(1)\oplus \mathfrak{ba}(1)^*$ constructed in \S~\ref{SecDef}.  In this case, $\fa_\od=\text{Span } \{\theta,z,\theta^*, z^*\}.$ A direct computation using Eq. (\ref{squaringh}) shows that (for any $s,w,u,v\in \Kee$)
\[
\begin{array}{lcl}
s_\fa(s \theta+wz+u\theta^*+v z^*)&=& s(u\theta^*+v z^*)\circ \ad_{\theta}+w(u\theta^*+v z^*)\circ \ad_{z}=svq^*.\\[2mm]
\end{array}
\]
A direct computation using Eqs. (\ref{squaringh}) and (\ref{bracketh}) shows that the only nonzero brackets:
\[
[\theta,q]_\fa = z,\quad [q,z^*]_\fa = \theta^*,\quad [\theta, z^*]_\fa=q^*.
\]

\sssbegin{Claim}
\label{der1ba}
The space $\mathfrak{out}(\fa)$ is spanned by the (classes of the) following cocycles 
(odd cocycles are underlined):
\[\footnotesize
\begin{array}{lcllcllcllcl}
D_1&=&q^*\otimes \widehat{q}, &\underline{D_2}&=&q^*\otimes \widehat{\theta}, &\underline{D_3}&=&q^*\otimes \widehat{z^*},\\[2mm] 
\underline{D_4}&=&q^*\otimes \widehat{\theta^*}  +\theta\otimes \widehat{q},&D_5&=&\theta^* \otimes \widehat{\theta}, &D_6&=&\theta^*\otimes \widehat{z^*},\\[2mm]
D_7&=&z\otimes \widehat{z^*},& \underline{D_8}&=&q^* \otimes \widehat{z} +z^* \otimes \widehat{q},&D_9&=&q\otimes \widehat{q} +
   \theta^*\otimes \widehat{\theta^*} + z \otimes \widehat{z}, \\[2mm]
D_{10}&=&q^* \otimes \widehat{q^*}+ \theta\otimes \widehat{\theta}  +
   z\otimes \widehat{z},&&&&
 D_{11}&=&q^*\otimes \widehat{q^*}+
 \theta^*\otimes \widehat{\theta^*} +z^*\otimes \widehat{z^*}.\end{array}
\]
\end{Claim}


Let us fix an ordered basis in $\fa$ as follows: $\theta,q,z,\theta^*,q^*,z^*$. In this basis, the Gram matrix of the bilinear form $B_\fa$ in (\ref{Mform}) is given by 
\[
\left (
\begin{array}{c|c}
0&I_3\\
\hline 
I_3& 0
\end{array}
\right ).
\]
Any derivation $D$ has the following supermatrix representation:
\[
D=\left (
\begin{array}{c|c}
A&B\\
\hline
C& F
\end{array}
\right ).
\] 
It follows that $D$ is compatible with the bilinear form $B_\fa$ if and only if $F=A^t, B^t=B$ and $C^t=C$. Let us consider the most general derivation $D=\mathop{\sum}\limits_{1\leq i\leq 11}\alpha_i D_i$, where $D_i$ are the cocycles given in Claim~\ref{der1ba}. In the same basis $\theta,q,z,\theta^*,q^*,z^*$, we have
\tiny
\[
D=\left (
\begin{array}{c|c}
\alpha_{10} E^{1,1}+\alpha_4 E^{1,2}+\alpha_9E^{2,2}+(\alpha_9+\alpha_{10})E^{3,3}&\alpha_{7} E^{3,3}\\
\hline
\alpha_5 E^{1,1}+\alpha_2 E^{2,1}+\alpha_1 E^{2,2}+\alpha_8(E^{2,3}+E^{3,2})&\alpha_{9}E^{1,1}+\alpha_{6}E^{1,3}+\alpha_{4} E^{2,1}+\alpha_{10} E^{2,2}+\alpha_3 E^{2,3}+\alpha_{11}I
\end{array}
\right ).
\]
\normalsize

It follows that $D$ is compatible with the bilinear form $B_\fa$ if and only if $\alpha_2=\alpha_3=\alpha_6=0$, and $\alpha_{11}=\alpha_9+\alpha_{10}$.

\subsubsection{The $D_\ev$-extension}



In this case, the most general even derivation is of the form 
\[
\alpha_1D_1+\alpha_5 D_5+\alpha_7 D_7+\alpha_9D_9+\alpha_{10}D_{10}+(\alpha_9+\alpha_{10})D_{11}. 
\]
Now, we define a quadratic form on $\fa_\od$ as follows: 
\[
\alpha(\mu \theta+ \nu z+\lambda \theta^*+\beta z^*)=\mu \lambda+\nu \beta+A(\nu^2+\beta^2),
\] 
where $A\in \Kee$, see \S \ref{Arf}.

Let $a=sq+tq^*$ be an even element. The fact that $B_\fa(a,D(a))=0$ implies that $\alpha_1=0$. Now, let 
\[
a=s\theta+t\theta^*+wz+uz^*\text{~~ and $b=\tilde s\theta+\tilde t\theta^*+\tilde wz+\tilde uz^*$}
\]
 be two odd elements. We have 
\small
\[
\begin{array}{lcl}
B_\fa(a, D(b))&=&B_\fa(s\theta+t\theta^*+wz+uz^*,( \alpha_5 \tilde s +\alpha_{10}\tilde t) \theta^*+ (\alpha_7 \tilde u +\alpha_9\tilde w+\alpha_{10}\tilde w) z+\alpha_{10}\tilde s \theta+\alpha_{11}\tilde u z^*)\\[2mm]
&=&s( \alpha_5 \tilde s+\alpha_{10}\tilde t)+u (\alpha_7 \tilde u +\alpha_9\tilde w+\alpha_{10}\tilde w)+t\alpha_{10}\tilde s+w (\alpha_{9}+\alpha_{10})\tilde u.
\end{array}
\]
\normalsize
On the other hand,  
\[
\begin{array}{lcl}
B_\alpha(a,b)&=&(s+\tilde s)(t+\tilde t)+(w+\tilde w)(u +\tilde u)+A((w+\tilde w)^2+(u+\tilde u)^2)\\
&&+st +wu+A (w^2+u^2)+\tilde s \tilde t +\tilde w \tilde u+A (\tilde w^2+\tilde u^2)\\[2mm]
&=& s \tilde t+ \tilde s t+w \tilde u+ \tilde w u.
\end{array}
\]
It follows that $\alpha_{10}=1$ and $\alpha_5=\alpha_7=\alpha_9=0 $. 


The even double extension of $\fa$ is then a NIS-Lie superalgebra of $\sdim=4|4$, by means of the derivation $D= \theta\otimes \widehat{\theta}  +
   z\otimes \widehat{z}+
 \theta^*\otimes \widehat{\theta^*} +z^*\otimes \widehat{z^*}$ and the quadratic form $\alpha$. 
\subsubsection{The $D_\od$-extension} In this case, the only odd derivations compatible with the bilinear form $B_\fa$ are $D_4$ and $D_8$. Put $D:=\alpha_4 D_4+\alpha_8 D_8$. Now, the only even element $a_0$ for which $D(a_0)=0$ is of the form $a_0=k q^*$, where $k\in \Kee$. A direct computation shows that $D^2=\ad_{k q^*}=0$. The odd double extension of $\fa$ is then a NIS-Lie superalgebra of $\sdim=2|6$, by means of the derivation $D$ and $a_0=k q^*$, where $k\in\Kee$. 

Let us show that there is an isometry between the extensions corresponding to pairs $(D, a_0=k q^*)$ and $(\tilde D=D_4, \tilde a_0=0)$. Indeed, the isometry is given by the formulas (for notation, see Theorem \ref{Isom2})
\[
\begin{array}{rclrclrclrcl}
\pi_0(z)&=& \alpha_4z+\alpha_{8}\theta^* & \pi_0(z^*)&=&s_3 \theta+\alpha_8\theta^*+\alpha_4z+s_6z^*, &  \pi_0(q)&=&q, \\[2mm]
 \pi_0(q^*)&=&q^*, & \pi_0(\theta)&=&\alpha_4 \theta+\alpha_{8}z^*,& \pi_0(\theta^*)&=&s_3 z+s_6 \theta^*,\\[2mm]
 t&=&k\theta^*, & \lambda&=&1,& \nu&=&0,
\end{array}
\]
where $s_6 \alpha_4+s_3 \alpha_8=1$.

On the other hand, let us show that the double extension by means of $D_4$ is not a trivial one; namely, it is not isometric to the one by means of $\ad_T$ for some $T\in \fa$. Suppose there is an isometry, say $\pi$. Let us write 
\[
\pi_0(q^*)=mq^*, \quad \pi_0(q)=m_1q^*+m^{-1}q.
\]
Now, because $z=[\theta,q]$, it follows that 
\[\pi_0(z)=[\pi_0(\theta),m_1q^*+m^{-1}q]=c_1z+c_2\theta^*\quad  \text{for some $c_1, c_2 \in\Kee$.}
\] 
Similarly, since $\theta^*=[q,z^*]$, it follows that 
\[
\pi_0(\theta^*)=[m_1q^*+m^{-1}q,\pi_0(z^*)]=\tilde c_1z+\tilde c_2\theta^*\quad  \text{for some $\tilde c_1, \tilde c_2\in \Kee$}.
\] 
We have (here $T=W_1z+W_2z^*+W_3\theta+W_4\theta^*$, where  $W_i\in \Kee$):
\[
\begin{array}{lcl}
(D_4\circ \pi_0+\pi_0\circ \ad_T)(q)&=&D_4(m_1q^*+m^{-1}q)+\pi_0 [T, q]\\[2mm]
&=&m^{-1}\theta+ W_3\pi_0(z)+W_2\pi_0(\theta^*)\\[2mm]
&=&m^{-1}q+ W_3(c_1z+c_2\theta^*)+W_2 (\tilde c_1z+\tilde c_2\theta^*).
\end{array}
\] 
But this is never zero, hence a contradiction.
\subsection{Purely odd superalgebra} Consider $\fa:=\Span\{ a,b\}$, where $a$ and $b$ are both odd, with the bracket $[a,b]_\fa=0$, and $s_\fa(g)=0$ for any $g$ odd. Define a bilinear form as follows: 
\[
B_\fa(a,a):=B_\fa(b,b):=0\text{~~and $B_\fa(a,b):=B_\fa(b,a):=1$. }
\]
Obviously, $B_\fa$ is even, symmetric and non-degenerate. The superalgebra $(\fa, B_\fa)$ is then a NIS-Lie superalgebra. 

\sssbegin{Claim}
The space $\mathfrak{out}(\fa)=\mathrm{H}^1(\fa; \fa)$ is spanned by the even cocycles: 
\[
\begin{array}{lcllcllcllcl}
D_1=&a \otimes \widehat{a}, &D_2=a \otimes \widehat{b}, & D_3&=&b\otimes \widehat{a}, & D_4&=b\otimes \widehat{b}.
\end{array}
\]
\end{Claim}


In the ordered basis $a, b$ of $\fa$, the matrix representation of the derivation $D=\mathop{\sum}\limits_{1\leq i\leq 4} \alpha_i D_i$ is
\[
\left (
\begin{array}{cc}
\alpha_1& \alpha_2\\
\alpha_3 & \alpha_4
\end{array}
\right ).
\]
It follows that $D$ is compatible with the bilinear form $B_\fa$, whose Gram matrix is $\mathrm{antidiag}\{1,1\}$, if and only if $\alpha_1=\alpha_4$.
 

We define a quadratic form on $\fa_\od$ by setting
\[
\alpha(\lambda a+ \mu b)=\lambda \mu, \quad \text{where $\lambda, \mu \in \Kee$}.
\]  
The condition (\ref{D3}) of Theorem \ref{MainTh} is satisfied since $B_\alpha(\lambda a+\beta b, \mu a+\nu b)=\lambda \nu+\beta \mu$; on the other hand,
\[
\begin{array}{lcl}
B_\fa(\lambda a+\beta b,D(\mu a+\nu b))&=&B_\fa(\lambda a+\beta b,\alpha_1(\mu a)+\alpha_4(\nu b)+\alpha_2(\nu a)+\alpha_3(\mu b))\\[2mm]
&=&(\beta \alpha_2+\lambda\alpha_4)\nu+(\beta \alpha_1+\lambda \alpha_3)\mu.
\end{array}
\]
It follows that $\alpha_1=\alpha_4=1$ and $\alpha_2=\alpha_3=0$. On the other hand, 
\[
B_\fa(D(\lambda a+\beta b), D(\lambda a+\beta b))=B(\lambda a+\beta b,\lambda a+\beta b)=2\lambda \beta=0.
\] 

The $D_\ev$-extension of $(\fa, B_\fa)$ is a NIS-Lie superalgebra of $\sdim=2|2$.

\subsection{An exceptional example: $\fpsl(2|2)\simeq\fh^{(1)}(0|4)$} Consider the Hamiltonian superalgebra $\fh(0|4)$, see \cite{BGLLS2}. As a vector space (here $\xi$'s and $\eta$'s are odd indeterminates) it can be considered as follows:
\[
\fh(0|4)\simeq\Span\{H_f\; |\; f\in \Kee[\xi,\eta]\}\simeq \Kee[\xi,\eta]/\Kee\cdot 1,\]
where  \[
H_f= \frac{\partial f}{\partial {\xi_1}} \frac{\partial  }{\partial {\eta_1}} +\frac{\partial f}{\partial {\eta_1}} \frac{\partial  }{\partial {\xi_1}}+\frac{\partial f}{\partial {\xi_2}} \frac{\partial  }{\partial {\eta_2}} +\frac{\partial f}{\partial {\eta_2}} \frac{\partial  }{\partial {\xi_2}}.\]
The Lie bracket $[H_f, H_g]=H_{\{f,g\}}$ is given by the Poisson bracket:
\[
\{f,g\}:= \frac{\partial f}{\partial {\xi_1}} \frac{\partial  g}{\partial {\eta_1}} +\frac{\partial f}{\partial {\eta_1}} \frac{\partial g }{\partial {\xi_1}}+\frac{\partial f}{\partial {\xi_2}} \frac{\partial  g}{\partial {\eta_2}} +\frac{\partial f}{\partial {\eta_2}} \frac{\partial g }{\partial {\xi_2}}.
\]
The derived Lie superalgebra $\fh^{(1)}(0|4)$ admits an invariant non-degenerate (super)symmetric bilinear form given by the Berezin integral
\[
B_{\fh^{(1)}(0|4)}(f,g):=\int fg \, \text{vol}(\xi,\eta) \quad \text{(=the coefficient of the monomial $\xi_1\xi_2\eta_1\eta_2$).}
\]

\sssbegin{Claim} The  space $\mathfrak{out}(\fh^{(1)}(0|4))=\mathrm{H}^1(\fh^{(1)}(0|4); \fh^{(1)}(0|4))$ is spanned by the cocycles: 
\be\label{cocy}\tiny
\begin{array}{ll}
\deg=-2:&D_1=\xi _1  \otimes  \widehat{\left(\xi _1\, \xi _2\, \eta _2\right)}+\xi
   _2  \otimes\widehat{\left(\xi _1\, \xi _2\, \eta _1\right)}+\eta _1 \otimes\widehat{ \left(\xi _2\,
   \eta _1\, \eta _2\right)}+\eta _2 \otimes\widehat{ \left(\xi _1\, \eta _1\, \eta _2\right)};\\[3mm]

\deg=0:  &D_2=\xi _1\otimes \widehat{ \eta _1}+( \xi _1\, \xi
   _2)\otimes  \widehat{ \xi _2\, \eta _1}+  (\xi _1\,
   \eta _2) \otimes \widehat{\eta _1\, \eta _2} +  (\xi
   _1\, \xi _2\, \eta _2) \otimes \widehat{\left(\xi _2\, \eta _1\, \eta _2\right)};\\[2mm]
   
    &D_3= \xi _2 \otimes \widehat{ \eta _2} +  (\xi _1\, \xi
   _2) \otimes \widehat{\xi _1\, \eta _2} + ( \xi _2\,
   \eta _1 )\otimes \widehat{\eta _1\, \eta _2}+ ( \xi
   _1\, \xi _2\, \eta _1 )\otimes \widehat{\left (\xi _1\, \eta _1\, \eta _2\right)};\\[2mm]
   
    &D_4= \eta _1 \otimes\widehat{\xi _1}+  (\xi _2\, \eta
   _1) \otimes\widehat{\xi _1\, \xi _2}+(\eta _1\,
   \eta _2) \otimes \widehat{ \xi _1\, \eta _2}+ ( \xi
   _2\, \eta _1\, \eta _2)  \otimes\widehat{\left(\xi _1\, \xi _2\, \eta _2\right)};\\[2mm]
   
   &D_5=\eta _2 \otimes \widehat{\xi _2}+ (\xi _1\, \eta
   _2)  \otimes\widehat{\xi _1\, \xi _2}+ (\eta _1\,
   \eta _2) \otimes \widehat{\xi _2\, \eta _1}+ (\xi
   _1\, \eta _1\, \eta _2)  \otimes\widehat{\left(\xi _1\, \xi _2\, \eta _1\right)};\\[2mm]

&D_6=\xi _2 \otimes \widehat{\xi _2}+ \eta
   _1 \otimes\widehat{\eta _1}+ (\xi _1\, \eta
   _2) \otimes \widehat{\xi _1\, \eta _2}+ (\xi _2\,
   \eta _1) \otimes \widehat{\xi _2\, \eta _1}+ (\xi
   _1\, \xi _2\, \eta _2) \otimes \widehat{\left(\xi _1\, \xi _2\, \eta _2\right)} + (\xi _1\, \eta _1\, \eta _2) \otimes \widehat{\left(\xi _1\, \eta
   _1\, \eta _2\right)};\\[3mm]
   
\deg=2:&D_7= (\xi _1\, \xi _2\, \eta _1)\otimes  \widehat{\xi
   _2}+(\xi _1\, \xi _2\, \eta _2) \otimes \widehat{\xi
   _1}+(\xi _1\, \eta _1\, \eta _2)  \otimes\widehat{\eta
   _2}+(\xi _2\, \eta _1\, \eta _2)\otimes  \widehat{\eta _1}.
\end{array}
\ee
\end{Claim}


Fix a lexicographically ordered basis on $\fh^{(1)}(0|4)$. In this basis, we identify the bilinear form $B$ with its Gram matrix $\antidiag(1,...,1)$. 
All derivations~\eqref{cocy} are compatible with $B$, i.e., they satisfy 
\[
\begin{array}{lcl}
B_{\fh^{(1)}(0|4)}(D(f),g)&=&B_{\fh^{(1)}(0|4)}(f,D(g)) \quad \text{for any $f,g \in \fh^{(1)}(0|4)$ and}\\[2mm]
B_{\fh^{(1)}(0|4)}(D(f),f)&=&0 \quad  \text{for any $f \in \fh^{(1)}(0|4)_\ev.$} 
\end{array}
\]

Let us give the proof only for the cocycle $D_1$. The proof is identical for the other derivations.  The matrix representation of $D_1$ in the same basis is $D_1\simeq E^{1,12}+E^{2,11}+ E^{3,14}+E^{4,13}$. Now, the condition $B^tD_1=D_1B$ is easily seen. Besides, since $D_1$ acts by zero on the even part, it follows that 
\[
B_{\fh^{(1)}(0|4)}(D_1(f),f)=0\text{~~ for any $f \in \fh^{(1)}(0|4)_\ev$. }
\]

Let $a=\lambda_1 \xi_1+\ldots+\lambda_8 \xi_2\eta_1\eta_2$. The quadratic forms associated with $D_1$ and $D_7$ are given, respectively, by:
\[
\alpha_1(a)=\lambda_6 \lambda_8+\lambda_5\lambda_7, \quad \alpha_7(a)=\lambda_2 \lambda_4+\lambda_1\lambda_3.
\]
Let us show that, up to an isometry, the derivations $(D_1,\alpha_1)$ and $(D_7,\alpha_7)$ give the same Lie superalgebra. Indeed, the isometry is given by (other generators are fixed):
\[
\xi_1\longleftrightarrow \xi_1\xi_2\eta_2, \; \xi_2\longleftrightarrow \xi_1\xi_2\eta_1, \;\eta_1\longleftrightarrow \xi_2\eta_1\eta_2, \; \eta_2\longleftrightarrow \xi_1\eta_1\eta_2, \;\lambda=1,\; \nu=0,\; t=0.
\]
On the other hand, the even double extension of the Lie superalgebra $\fh^{(1)}(0|4)$ by means of $(D_7,\alpha_7)$ is isomorphic to $\mathfrak{po}(0|4)$. Indeed the isomorphism is given by 
\[
\varphi(x)=1, \quad \varphi(x^*)= \xi_1\xi_2\eta_1\eta_2, \quad  {{\varphi}_|}_{\fh^{(1)}(0|4)}=\id.\]
All extensions by means of the derivations $D_2, D_3, D_4$ and $D_5$ are isometric. Here is the list of isometries  relating $D_2$ with $D_3$, $D_4$ and $D_5$:
\[
\begin{array}{ll}
D_3: & \xi_1\longleftrightarrow \xi_2, \quad \xi_2\longleftrightarrow \xi_1, \quad \eta_1\longleftrightarrow \eta_2, \quad \eta_2\longleftrightarrow \eta_1, \quad \lambda=1, \quad \nu=0; \quad t=0,\\
D_4: & \xi_1\longleftrightarrow \eta_1, \quad \xi_2\longleftrightarrow \xi_2, \quad \eta_2\longleftrightarrow \eta_2, \quad \eta_1\longleftrightarrow \xi_1, \quad \lambda=1, \quad \nu=0; \quad t=0,\\ 
D_5: &  \xi_1\longleftrightarrow \eta_2, \quad \eta_1\longleftrightarrow \xi_2, \quad \eta_2\longleftrightarrow \xi_1, \quad \xi_2\longleftrightarrow \eta_1, \quad \lambda=1, \quad \nu=0; \quad t=0.
\end{array}
\]
then extended to monomials in $\xi$'s and $\eta$'s.

Moreover, the double extension of $\fh^{(1)}(0|4)$ by means of $(D_6, \alpha_6)$, see Table~\eqref{h04}, is isometric to $\mathfrak{gl}(2|2)$ with the standard NIS given by the \textit{supertrace}. The isometry is explicitly given on generators by the following correspondences (other elements are obtained by bracketing)
\[
\begin{array}{lcllcllcllcllcl}
\eta_1&\longleftrightarrow &E^{3,2}, &\xi_1\eta_2 &\longleftrightarrow &E^{2,1}, & \xi_1\xi_2& \longleftrightarrow &E^{4,3}, &\xi_1\xi_2 \eta_2 & \longleftrightarrow & E^{2,3},\\[2mm]
\xi_2\eta_1 &\longleftrightarrow & E^{1,2},& \eta_1\eta_2 & \longleftrightarrow & E^{3,4}, &
x &\longleftrightarrow & I, & x^* & \longleftrightarrow & E^{2,2}.
\end{array}
\]
Table (\ref{h04}) gives the quadratic form associated with each derivation, and the respective even double extension, up to an isometry

\begin{equation}\label{h04}
\footnotesize
\renewcommand{\arraystretch}{1.4}
\begin{tabular}{|c|c|c|} \hline
Derivation & $\alpha(a)$&
Double extension\\
\hline
$D_2$ & $\lambda_3 \lambda_8$ & $\widetilde{\mathfrak{po}}(0|4)$ \\ \hline
$D_6$ & $\lambda_2 \lambda_7+\lambda_3 \lambda_6$ & $\fgl(2|2)$\\ \hline
$D_7$ & $\lambda_2 \lambda_4+\lambda_1 \lambda_3$ & $\mathfrak{po}(0|4)$\\ \hline
\end{tabular}
\end{equation}

\sssbegin{Claim}  $\dim\mathfrak{out}(\mathfrak{po}(0|4))=3$, $\dim\mathfrak{out}(\mathfrak{gl}(2|2))=1$ and $\dim\mathfrak{out}(\widetilde{\mathfrak{po}}(0|4))=5$, hence $\widetilde{\mathfrak{po}}(0|4)$, $\mathfrak{gl}(2|2)$ and $\mathfrak{po}(0|4)$ are pairwise not isomorphic. 
\end{Claim}
The Lie superalgebra $\widetilde{{\mathfrak{po}}}(0|4)$ has no analog for $p\not =2$. To show that, we need the following result:  
\sssbegin{Claim} For $p\not =2$, the  space $\mathrm{H}^1(\fh^{(1)}(0|4); \fh^{(1)}(0|4))$ is spanned by three cocycles whose degrees are $-2,0 $ and $2$.
\end{Claim}
The $D$-extension of $\fh^{(1)}(0|4)$ corresponding to the derivatives of degree $2$ and $-2$ are isometric to the Poisson Lie superalgebra $\fpo(0|4)$. The $D$-extension corresponding to the derivative of degree $0$ is isometric to $\mathfrak{gl}(2|2)$. 

\section{Examples: the case where the bilinear form $B$ is odd} \label{Exa2}
We give here an example of $D$-extensions in the case where the bilinear form $B$ is odd.
\subsection{The Lie superalgebra $\fh^{(1)}(0|5)$}

Consider the Hamiltonian superalgebra $\fh(0|5)$, see \cite{BGLLS2}. As a vector space (here $\xi$'s, $\eta$'s and $\theta$ are odd indeterminates) it can be considered as follows:
\[
\fh(0|5)\simeq\Span\{H_f\; |\; f\in \Kee[\xi,\eta,\theta]\}\simeq \Kee[\xi,\eta,\theta]/\Kee\cdot 1,\]
where  
\[
H_f= \frac{\partial f}{\partial {\xi_1}} \frac{\partial  }{\partial {\eta_1}} +\frac{\partial f}{\partial {\eta_1}} \frac{\partial  }{\partial {\xi_1}}+\frac{\partial f}{\partial {\xi_2}} \frac{\partial  }{\partial {\eta_2}} +\frac{\partial f}{\partial {\eta_2}} \frac{\partial  }{\partial {\xi_2}}+\frac{\partial f}{\partial {\theta}} \frac{\partial  }{\partial {\theta}}.\]
The Lie bracket $[H_f, H_g]=H_{\{f,g\}}$ is given by the Poisson bracket:
\[
\{f,g\}:= \frac{\partial f}{\partial {\xi_1}} \frac{\partial  g}{\partial {\eta_1}} +\frac{\partial f}{\partial {\eta_1}} \frac{\partial g }{\partial {\xi_1}}+\frac{\partial f}{\partial {\xi_2}} \frac{\partial  g}{\partial {\eta_2}} +\frac{\partial f}{\partial {\eta_2}} \frac{\partial g }{\partial {\xi_2}}+\frac{\partial f}{\partial {\theta}} \frac{\partial g }{\partial {\theta}}.
\]
The derived Lie superalgebra $\fh^{(1)}(0|5)$ admits an invariant non-degenerate odd (super)symmetric bilinear form given by the Berezin integral
\[
B_{\fh^{(1)}(0|5)}(f,g):=\int fg \, \text{vol} \quad \text{(=the coefficient of the monomial $\xi_1\xi_2\eta_1\eta_2\theta$).}
\]

\sssbegin{Claim} The  space $\mathfrak{out}(\fh^{(1)}(0|5))=\mathrm{H}^1(\fh^{(1)}(0|5); \fh^{(1)}(0|5))$ is spanned by the cocycles: (the odd cocycle is underlined)
\be\label{cocy2}\tiny
\begin{array}{lll}
\deg=0:&D_1=&\xi _1\otimes\widehat {\eta_1}+ (\xi_1\, \theta)\otimes \widehat {
   \eta _1 \theta }+
   (\xi _1\, \zeta_2 )\otimes \widehat {\xi_2\, \eta_1}+(\xi_1\, \eta _2 )\widehat{ (\eta _1
   \eta _2)}+ (\xi_1\, \xi_2\, \theta) \otimes\widehat {\xi_2\, \eta_1\, \theta} + (\xi_1\, \eta _2\, \theta)\otimes \widehat {(\eta _1\, \eta _2\,\theta)}\\
 &&  +(\xi_1\, \xi_2\, \eta_2)\otimes \widehat {(\xi_2\, \eta
   _1\, \eta _2)}+ (\zeta
   _1\, \xi_2\, \eta _2\, \theta) \otimes\widehat {(\xi_2\, \eta_1\, \eta _2\, \theta)};\\[3mm]

&D_2=&\xi_2 \otimes\widehat {\eta_2}+ (\xi_2\, \theta)\otimes\widehat {
  ( \eta _2\, \theta)} +(\xi _1\, \xi_2) \otimes\widehat {(\xi _1\, \eta_2)} + (\xi_2\, \eta _1)\otimes\widehat {(\eta _1\,\eta _2)}+
  (\xi _1\, \xi _2\, \theta) \otimes\widehat {(\xi _1\, \eta_2\, \theta)}+(\xi _2\, \eta _1\, \theta)\otimes\widehat {(\eta _1\, \eta _2\,\theta)} \\ &&+
   (\xi _1\, \xi _2\, \eta)\otimes\widehat { (\xi _1\, \eta_1\, \eta _2)}+ (\xi_1\, \xi _2\, \eta _1\, \theta)\otimes\widehat {( \xi _1\, \eta_1\, \eta _2\, \theta)};\\[3mm]
   
    &D_3=& \eta _1\otimes\widehat {\xi_1}+ (\eta_1\, \theta) \otimes\widehat {\xi _1\, \theta}+(\xi _2\, \eta_1)\otimes \widehat {(\xi _1\, \xi_2)}+(\eta_1\, \eta _2)\otimes \widehat {(\xi_1\, \eta _2)}+(\xi _2\, \eta _1\, \theta)\otimes\widehat {(\xi _1\, \xi_2\, \theta) }+(\eta _1\, \eta _2\, \theta) \otimes\widehat {(\xi _1\, \eta_2\, \theta) }\\    &&+(\xi _2\, \eta _1\, \eta_2) \otimes\widehat { (\xi _1\, \xi_2\, \eta _2)}+(\xi_2\, \eta _1\, \eta _2\, \theta)\otimes\widehat {( \xi _1\,\xi _2\, \eta _2\, \theta)};\\[3mm]
   
    &D_4=& \eta _2 \otimes\widehat { \xi_2}+ (\eta_2\, \theta)\otimes \widehat {
  ( \xi _2\, \theta)}+
(\xi _1\, \eta_2) \otimes\widehat { \xi _1\, \xi
   _2} + (\eta_1\, \eta _2) \otimes\widehat { (\xi_2\, \eta _1)}+
(\xi _1\, \eta _2\, \theta) \otimes\widehat {(\xi _1\, \xi_2\, \theta) }+(\eta _1\, \eta _2\, \theta)\otimes \widehat {(\xi _2\, \eta_1\, \theta) }\\
&&+ (\xi _1\, \eta _1\, \eta_2) \otimes\widehat { (\xi _1\, \xi_2\, \eta _1)}+(\xi_1\, \eta _1\, \eta _2\, \theta) \otimes\widehat {( \xi _1\,
   \xi _2\, \eta _1\, \theta)};\\[3mm]
   
   &D_5=&\theta\otimes \widehat {\theta} + \xi_1\otimes \widehat {\xi _1}+ \xi _2\otimes \widehat {\xi_2 } + \eta_1 \otimes\widehat {\eta _1}+ \eta _2\otimes \widehat {\eta_2} + (\xi_1\, \xi _2\, \theta ) \otimes\widehat { (\xi _1\, \xi_2\, \theta )}  + (\xi _1\, \eta _1\, \theta)\otimes\widehat { (\xi _1\, \eta_1\, \theta) }\\  && +
 (\xi _1\, \eta _2\, \theta) \otimes\widehat { (\xi _1\, \eta_2\, \theta) }+ (\xi _2\, \eta _1\, \theta) \otimes\widehat { (\xi _2\, \eta_1\, \theta) }
 +(\xi _2\, \eta _2\, \theta)\otimes\widehat { (\xi _2\, \eta_2\, \theta)}+(\eta _1\, \eta _2\, \theta)\otimes\widehat { (\eta _1\, \eta _2\,\theta) }\\
 &&
  + (\xi _1\, \xi _2\, \eta_1)\otimes\widehat { (\xi _1\, \xi_2\, \eta _1 )}+ (\xi _1\, \xi _2\, \eta_2)\otimes \widehat { (\xi _1\, \xi_2\, \eta _2)} +
   (\xi _1\, \eta _1\, \eta_2 ) \otimes\widehat { (\xi _1\, \eta_1\, \eta _2 )} +(\xi _2\, \eta _1\, \eta_2 ) \otimes\widehat { (\xi _2\, \eta_1\, \eta _2)};\\[3mm]

\deg=3:&\underline{D_6}=&(\xi _1\, \xi _2\, \eta_1\, \theta)\otimes \widehat {\xi_2}+(\xi _1\, \xi _2\, \eta_2\, \theta)\otimes\widehat {\xi_1}+(\xi _1\, \eta _1\, \eta_2\, \theta)\otimes \widehat {\eta_2}+(\xi _2\, \eta _1\, \eta_2\, \theta)\otimes\widehat {\eta_1}+(\xi _1\, \xi _2\, \eta_1\, \eta _2)\otimes\widehat {\theta}\end{array}
\ee
\end{Claim}


\subsubsection{The $D_\ev$-extension}
Let us first show that the derivation $D_5$ is not compatible with the bilinear form. Indeed, 
\[
B_{\fh^{(1)}(0|5)}(D_5(\theta),\xi_1\xi_2\eta_1\eta_2)=1 \text{ while } B_{\fh^{(1)}(0|5)}(\theta,D_5(\xi_1\xi_2\eta_1\eta_2))=0.
\]
Now, the derivations $D_1, D_2, D_3$ and $D_4$ are compatible with the bilinear form $B_{\fh^{(1)}(0|5)}$, and the proof is similar to that of $\fh^{(1)}(0|4)$ in \S~\ref{Exa}.

All extensions by means of the derivations $D_1, D_2, D_3$ and $D_4$ are isometric. Here is the list of isometries  relating $D_1$ with $D_2$, $D_3$ and $D_4$:
\[
\begin{array}{ll}
D_2: & \xi_1\longleftrightarrow \xi_2, \quad \xi_2\longleftrightarrow \xi_1, \quad \eta_1\longleftrightarrow \eta_2, \quad \eta_2\longleftrightarrow \eta_1,\quad \theta\longleftrightarrow \theta, \quad \lambda=1, \quad \nu=0; \quad t=0,\\
D_3: & \xi_1\longleftrightarrow \eta_1, \quad \xi_2\longleftrightarrow \xi_2, \quad \eta_2\longleftrightarrow \eta_2, \quad \eta_1\longleftrightarrow \xi_1, \quad \theta\longleftrightarrow \theta, \quad \lambda=1, \quad \nu=0; \quad t=0,\\ 
D_4: &  \xi_1\longleftrightarrow \eta_2, \quad \eta_1\longleftrightarrow \xi_2, \quad \eta_2\longleftrightarrow \eta_1, \quad \xi_2\longleftrightarrow \xi_1,\quad \theta\longleftrightarrow \theta,  \quad \lambda=1, \quad \nu=0; \quad t=0,
\end{array}
\]
then extended to monomials in $\xi$'s and $\eta$'s. The double extension by means of the derivation $D_1$ is a Lie superalgebra that we denote by $\widetilde {\mathfrak{po}}(0|5)$. 

Table\eqref{h05p} summarizes these results:

\begin{equation}\label{h05p}
\footnotesize
\renewcommand{\arraystretch}{1.4}
\begin{tabular}{|c|c|c|c|} \hline
Derivation & Compatibility with $B_{\fh^{(1)}(0|5)}$& $s_\fg(x)$ & $D_\ev$-extension\\
\hline
$D_1$ & Yes & 0 & $\widetilde {\mathfrak{po}}(0|5)$\\ \hline
$D_5$ & No & -- & --\\ \hline
\end{tabular}
\end{equation}

\subsubsection{The $D_\od$-extension}
Fix a lexicographically ordered basis on $\fh^{(1)}(0|5)$. In this basis, we identify the bilinear form $B_{\fh^{(1)}(0|5)}$ with its Gram matrix $\antidiag(1,...,1)$. 
The derivation $D_6$ is compatible with $B_{\fh^{(1)}(0|5)}$, i.e., satisfies
\[
\begin{array}{lcl}
B_{\fh^{(1)}(0|5)}(D_6(f),g)&=&B_{\fh^{(1)}(0|5)}(f,D_6(g)) \quad \text{for any $f,g \in \fh^{(1)}(0|5)$ and} \\[3mm]
B_{\fh^{(1)}(0|5)}(D_6(f),f)&=&0 \quad \text{for any $f\in \fh^{(1)}(0|5)_\ev$}.
\end{array}
\]
%
The matrix representation of $D_6$ in the same basis is 
\[
D_6\simeq E^{26,3}+E^{27,2}+ E^{28,5}+E^{29,4}+E^{30,1}. 
\]
Now, the condition $B^tD_6=D_6B$ is easily seen. Besides, since $D_6$ acts by zero on the even part, it follows that 
\[
B_{\fh^{(1)}(0|5)}(D_6(f),f)=0\text{~~ for any $f \in \fh^{(1)}(0|5)_\ev$. }
\]
Let $a=\lambda_1\theta+\lambda_2 \xi_1+\ldots+\lambda_{15} \xi_2\eta_1\eta_2$. The quadratic form associated with $D_6$ is given by:
\[
\alpha_6(a)=\lambda_3 \lambda_5+\lambda_2\lambda_4.
\]
Now, since $0=D^2=\ad_{a_0}$, it follows that $a_0=0$ because $\fh^{(1)}(0|5)$ has no center. We have, therefore, a parametric family of double extensions by means of $D_6$, $\alpha_6$, $a_0=0$ and $m\in \Kee$ (see Theorem \ref{MainTh3}) that we denote by $\mathfrak{po}(0|5;m)$. We will show that this family is not isometric to $\mathfrak{po}(0|5;0)$ (the usual Poisson algebra $\mathfrak{po}(0|5)$). Indeed, suppose there is such an isometry, say $\pi$, between $\mathfrak{po}(0|5;m)$ and $\mathfrak{po}(0|5)$. It follows that (see Theorem \ref{Isom3})
\[
D_6\circ \pi_0=\lambda \pi_0 \circ D_6+\ad_t \quad \text{ for some $t\in \mathfrak{h}^{(1)}(0|5)_\od$.}
\] 
Let us evaluate the equation above at any arbitrary  $a \in \mathfrak{h}^{(1)}(0|5)_\ev$. Since $D_6(a)=0$ and $\pi_0$ is even, it follows that $\ad_t(a)=0$ for any $a \in \mathfrak{h}^{(1)}(0|5)_\ev$. A direct computation shows that $t=0$. On the other hand, Eq. (\ref{3Cf}) of Theorem \ref{Isom3} implies that $m =0$. A contradiction. 

On the other hand, the even double extension $\mathfrak{po}(0|5;0)$ is isomorphic to $\mathfrak{po}(0|5)$. Indeed the isomorphism is given by 
\[
\varphi(x)=1, \quad \varphi(e)= \xi_1\xi_2\eta_1\eta_2\theta, \quad  {{\varphi}_|}_{\fh^{(1)}(0|5)}=\id.
\]

The table below summarizes these results:
\begin{equation}\label{h05}
\footnotesize
\renewcommand{\arraystretch}{1.4}
\begin{tabular}{|c|c|c|c|} \hline
Derivation & $\alpha(a)$ & $s_\fg(e)$ & Double extension\\
\hline
$D_6$ & $\lambda_3 \lambda_5+\lambda_2\lambda_4$&$0$ & $\mathfrak{po}(0|5)$\\ \hline
$D_6$ & $\lambda_3 \lambda_5+\lambda_2\lambda_4$&$mx$ & $\mathfrak{po}(0|5;m), m\not =0$\\ \hline
\end{tabular}
\end{equation}

\section{Nilpotent NIS-Lie superalgebras}
The goal of this section is to prove that every 2-step nilpotent NIS-Lie superalgebra in characteristic 2 can be obtained by an inductive process of $D$-extensions. 

Let $(\fg,B_\fg)$ be a NIS-Lie superalgebra such that $B_\fg$ is even. Let $V$ be a subspace of $\fg$. We define the following subspace
\[
V^{\sharp}:=\{x \in \fg\mid B_\fg(x,V)=0, \text{ and } B_\fg(s_\fg(x),V)=0 \text{ whenever $x\in\fg_\od$}  \}.
\]

\ssbegin{Proposition}\label{ide}
Let $I$ be an ideal of $\fg$. The subspace $I^\sharp$ is an ideal of $\fg$. Moreover, $I^\perp$ is an ideal of $\fg$ if and only if $I^\perp=I^\sharp$.
\end{Proposition}

\begin{proof}
Let $x\in I^\sharp$ and $y\in\fg$. We have 
\[
B_\fg([x,y],z)=B_\fg(x,[y,z])=0 \quad \text{for any $z\in I$}.
\] 
It follows that $[x,y]\in I^\sharp$. Suppose now that $x \in (I^\sharp)_\od$. Since $B(s_\fg(x),I)=0$, it follows that $s_\fg(x) \in I^\sharp$. Therefore, $I^\sharp$ is an ideal of $\fg$. Suppose now that $I^\perp$ is an ideal. Let $x\in I^\perp$. If $x$ is even then $x\in I^\sharp$ and we are done. If $x$ is odd, then $s_\fg(x)\in I^\perp$ because $I^\perp$ is an ideal. Therefore, $B_\fg(s_\fg(x),I)=0$. Thus, $x\in I^\sharp$. Conversely, suppose that $I^\perp=I^\sharp$. It follows that $I^\perp$ is an ideal because $I^\sharp$ is an ideal.
\end{proof}

\ssbegin{Lemma}
If $(\fg, B_\fg)$ is a NIS-Lie superalgebra, then
\begin{itemize}
\item[(i)]  $([\fg,\fg])^\perp=\fz(\fg)$;
\item[(ii)]  $([\fg,\fg])^\sharp=\fz_s(\fg)$;
\item[(iii)] $(\fg^{(1)})^\perp=\fz_s(\fg)$.
\end{itemize}
\end{Lemma}

\begin{proof}
(i) Let $z\in\fz(\fg)$. Let $x, y \in \fg$. We have $
B_\fg(z,[x,y])=B_\fg([z,x],y)=0$. It follows that $z\in ([\fg, \fg])^\perp$. The other way around, let $z \in ([\fg,\fg])^\perp$. Let $x, y \in \fg$.  We have 
\[
B_\fg([z,x],y)=B(z,[x,y])=0. 
\]
It follows that $[z,x]=0$ since $B_\fg$ is non-degenerate. Therefore, $z \in \fz(\fg)$. 

(ii) Since $([\fg,\fg])^\sharp=([\fg,\fg])^\perp \cap \{x\in \fg_\od\mid B_\fg(s_\fg(x),V)=0\}$, the result follows from (i).

(iii) Using Part (i) we have
\[
\begin{array}{lcl}
(\fz_s(\fg))^\perp&=&\left (\fz(\fg)\cap (s_\fg(\fg_\od))^\perp \right )^\perp\\[2mm]
&=&[\fg,\fg] + s_\fg(\fg_\od)\\[2mm]
&=& \fg^{(1)}.\qed
\end{array}
\]
\noqed\end{proof}

\ssbegin{Proposition} 
If $\Kee$ is a perfect field and $\dim(\fz(\fg)_\od)>\dim(\fz(\fg)_\ev)$, then $(\fg,B_\fg)$ is obtained as an $D_\od$-extension of a NIS-Lie superalgebra. 
\end{Proposition}

\begin{proof}
Let us show first that $s_\fg$ becomes additive when restricted on $\fz(\fg)_\od$. Indeed, let $x,y\in \fz(\fg)_\od$. We see that 
\[
0=[x,y]_\fg=s_\fg(x+y)-s_\fg(x)-s_\fg(y).
\] 
The result follows. 

Let us show that if the squaring ${{s_\fg}}|_{\fz(\fg)_\od}$ is injective, then 
\[
\dim(\fz(\fg)_\od)\leq \dim(\fz(\fg)_\ev). 
\]
Suppose that $\fz(\fg)_\od=\Span\{e_i \, |\, i=1,\ldots,n\}.$ If $\mathop{\sum}\limits_{i}\alpha_i s_\fg(e_i)=0$, then 
$s_\fg(\mathop{\sum}\limits_{i} \sqrt{\alpha_i} e_i)=0$. Since $s_\fg$ in injective, it follows that $\mathop{\sum}\limits_{i} \sqrt{\alpha_i} e_i=0$, and hence $\alpha_i=0$ for any $i=1,\ldots,n$. Thus,  $\dim(\fz(\fg)_\od)\leq \dim(\fz(\fg)_\ev)$. By our assumption, the squaring ${{s_\fg}}|_{\fz(\fg)_\od}$ must be non-injective. It follows that there exists a non-zero element $x\in \fz(\fg)_\od$ such that $s_\fg(x)=0$. Therefore, $\fz(\fg)_\od \cap \mathscr{C}(\fg,B_\fg)\not =\{0\}$ and Proposition \ref{Rec2} can be applied. 
\end{proof}

\ssbegin{Proposition}
Suppose that $(\fg, B_\fg)$ is a non-abelian irreducible 2-step nilpotent NIS-Lie superalgebra of $\dim(\fg)>2$. Then $(\fg,B_\fg)$ is obtained as an $D$-extension of a NIS-Lie superalgebra.
\end{Proposition}

\begin{proof}
Since $\fg$ is 2-step nilpotent, it follows that
\[
[\fg,[\fg,\fg]_\fg]_\fg=0, \quad [\fg,s_\fg(\fg_\od)]_\fg=0, \quad s_\fg([\fg_\ev,\fg_\od]_\fg)=0.
\]
If $[\fg_\ev,\fg_\od]_\fg\not =0$, then there exists $x\in [\fg_\ev,\fg_\od]_\fg$ such that $s_\fg(x)=0$. In addition, $x\in \fz(\fg)_\od$. Therefore, Proposition \ref{Rec2} can be applied, and hence $\fg$ can be obtained as a $D_\od$-extension from a NIS-Lie superalgebra. 

If $[\fg_\ev,\fg_\od]_\fg=0$, then $B_\fg(\fg_\ev,[\fg_\od,\fg_\od]_\fg)=0$ since $B_\fg$ is invariant. Therefore, $[\fg_\od,\fg_\od]_\fg=0$ since $B_\fg$ is even. It follows that $\fg_\od \subseteq \fz(\fg)_\od$ and, hence, $\fg_\od =\fz(\fg)_\od $. It follows that $s_\fg(\fg_\od)\subseteq \fz(\fg)_\ev$. 

Let us consider $I=s_\fg(\fg_\od)\oplus \fg_\od$. This is obviously an ideal.  If $\fg=I$, then $\fg$ is abelian, ruled out by hypothesis. If $I\not=\fg$, then we write $\fg=I\oplus I^\perp$. Moreover, $I^\perp=I^\sharp$ and using Proposition \ref{ide} it follows that $I^\perp$ is an ideal. Since $(\fg,B_\fg)$ is irreducible it follows that $I$ is degenerate. Hence, there exists a non-zero $x\in I$ such that $B_\fg(x,s_\fg(\fg_\od))=0$. Therefore, $x\in \fz_s(\fg)$, and hence $\fg$ can be obtained as a $D_\ev$-extension from a NIS-Lie superalgebra following the steps in Proposition \ref{Rec1}.
\end{proof}



\begin{thebibliography}{999}


\bibitem[ABB]{ABB} Albuquerque H., Barreiro E. and Benayadi S., Quadratic Lie superalgebras with a reductive even part. J. Pure Appl. Algebra, {\bf 213} (2009), 724-731.

\bibitem[ABBQ]{ABBQ} Albuquerque H., Barreiro E. and Benayadi S., Odd quadratic Lie superalgebras. J. of Geometry and Physics, {\bf 60} (2010), 230-250.

\bibitem[BBB]{BBB} 
Bajo I., Benayadi S. and Bordemann M., Generalized double extension and descriptions of quadratic Lie superalgebras, \texttt{arXiv:0712.0228}.

\bibitem[BB]{BB} 
Benamor H. and Benayadi S., Double extension of quadratic Lie superalgebras. Comm. Algebra, {\bf 27}, No. 1 (1999), 67-88.

\bibitem[B]{B} 
Benayadi S., Quadratic Lie superalgebras with completely reductive action of the even part on the odd part. J. of Algebra, {\bf 223} (2000), 344-366.

\bibitem[B2]{B2} 
Benayadi S., Socle and some invariants of quadratic Lie superalgebras. J. of Algebra, {\bf 261} (2003), 245-291.




\bibitem[Bor]{Bor} Bordemann M., Nondegenerate invariant bilinear forms on nonassociative algebras. Acta Math. Univ. Comenianae, Vol. LXVI, 2 (1997), 151--201.

\bibitem[BGLL]{BGLL}
Bouarroudj S., Grozman P., Lebedev A., Leites D., Divided power
(co)homology. Presentations of simple finite dimensional modular Lie
superalgebras with Cartan matrix. Homology, Homotopy and
Applications, Vol. 12 (2010), No. 1, 237--278;
\texttt{arXiv:0911.0243}

\bibitem[BGLL1]{BGLL1}
Bouarroudj S., Grozman P., Lebedev A., Leites D., Derivations and
central extensions of simple modular Lie algebras and superalgebras;
\texttt{arXiv:1307.1858}


\bibitem[BGLLS]{BGLLS}
Bouarroudj S., Grozman P., Lebedev A., Leites D., Shchepochkina I.,
Lie algebra deformations in characteristic 2.  Math. Research
Letters, v.~22 (2015) no.~2, 353--402; \texttt{arXiv:1301.2781}

\bibitem[BGLLS1]{BGLLS1}
Bouarroudj S., Grozman P., Lebedev A., Leites D., Shchepochkina I.,
New simple Lie algebras in characteristic $2$. IMRN, No. 18, (2016), 5695-5726; \texttt{arXiv:1307.1551}

\bibitem[BGLLS2]{BGLLS2}
Bouarroudj S., Grozman P., Lebedev A., Leites D., Shchepochkina I.,
Simple vectorial Lie algebras in characteristic $2$ and their
superizations; \texttt{arXiv:1510.07255}

%
\bibitem[BGL2]{BGL2}
Bouarroudj S., Grozman P., Leites D., Classification of finite
dimensional modular Lie superalgebras with indecomposable Cartan
matrix. Symmetry, Integrability and Geometry: Methods and
Applications (SIGMA), 5 (2009), 060, 63 pages;
\texttt{arXiv:math.RT/0710.5149}


\bibitem[BKLS]{BKLS}
Bouarroudj S., Krutov A., Leites D. Shchepochkina I., Non-degenerate invariant (super)symmetric bilinear forms on simple Lie (super)algebras. To appear in Algebras and Representation Theory.

\bibitem[BLLSq]{BLLSq}
Bouarroudj S., Lebedev A., Leites D., Shchepochkina I.,
Classifications of simple Lie superalgebras in characteristic $2$;
\texttt{arXiv:1407.1695}







\bibitem[CCLL]{CCLL}
Chapovalov D., Chapovalov M., Lebedev A., Leites D. The classification of almost affine (hyperbolic) Lie
superalgebras. J.~
Nonlinear Math. Phys., vol. 17 (2010), suppl. 1, Special issue in
memory of F.~Berezin, 103--161; \texttt{arXiv:0906.1860}

\bibitem[D]{D} 
Dye R. H., On the Arf invariant. J. Algebra 53 (1978), no. 1, 36--39.

\bibitem[Dr1]{Dr1} Drinfeld V.G., On quadratic commutation relations in the quasiclassical case. In: Mathematical Physics, Functional Analysis, ``Naukova Dumka'', Kiev (1986), 25--34 (in Russian).

\bibitem[Dr]{Dr} 
Drinfeld V. G., Quantum groups. J. Sov. Math. {\bf 41} (1988), 898-915.



\bibitem[FS]{FS} 
Favre G. and Santharoubane L. J., Symmetric, invariant, non-degenerate bilinear form on a Lie algebra. J. of Algebra, {\bf 105} (1987), 451--464. 

\bibitem[Gr]{Gr}
Grozman P., \textbf{SuperLie},
\url{http://www.equaonline.com/math/SuperLie}


\bibitem[K]{K}
Kac V., \textit{Infinite dimensional Lie algebras}. 3rd edition. Cambrodge U. Press (1995), xxi+400

\bibitem[K2]{K2} Kac V., Lie superalgebras, Adv. Math. {\bf 26} (1977), 8--96.

\bibitem[KS]{KS} 
Konstein S. E. and Stekocshchik R., Klein operator and the number of independent traces and supertraces on the superalgebra of observables of rational Calogero model based on the root system. J. Nonlinear Math. Phys., 20 (2013), 295--308.

\bibitem[KT]{KT} 
Konstein S. E. and Tyutin I. V., The number of independent traces and supertraces on symplectic reflection algebras. J. Nonlinear Math. Phys., 20 (2013), 295--308.

\bibitem[LeD]{LeD}
Lebedev A., Simple modular Lie~superalgebras. Ph.D. thesis. Leipzig
University, July, 2008.

\bibitem[LeD2]{LeD2} 
Lebedev A., Analogs of the orthogonal, Hamiltonian, Poisson, and contact Lie superalgebras in characteristic 2. J. Nonlinear Math. Phys. 17 (2010), suppl. 1, 217--251.

\bibitem[LSh]{LSh}
Leites D. and Shapovalov A., Manin--Olshansky triples and Lie superalgebras. J.~Nonlinear Math. Phys., v. 7, 2000, no.~2,
120--125; \texttt{arXiv:math.QA/0004186}

\bibitem[Ma]{Ma} 
Manin Yu. I., \textit{Quantum groups and non-commutative geometry}. CRM, Montr\'eal, 1988, 92pp.

\bibitem[MR]{MR} 
Medina A. and Revoy P. , Alg\`ebres de Lie et produit scalaire invariant. Ann. Scient. \'Ec. Norm. Sup., 4 s\'erie, {\bf 18} (1985), 553--561.


\bibitem[Sc]{Sc} Scheunert M., {\it The Theory of Lie Superalgebras}. Lectures Notes in Mathematics,
Vol. 716, Berlin, Heidelberg: Springer-Verlag, (1976). 

\bibitem[Se]{Se} Serganova V., Automorphisms of simple Lie superalgebras, Izv. Akad. Nauk SSSR Ser. Mat. {\bf 48} (1984),
585--598 (English transl.: Math. USSR-Izv. {\bf 24} (1985), 539--551).


\bibitem[S]{S}
Strade H. {\em Simple Lie algebras over fields of positive
characteristic. I. Structure theory.} de Gruyter Expositions in
Mathematics, 38. Walter de Gruyter \& Co., Berlin,  (2004) viii+540
pp; (2009) vi+385pp; (2012) x+239pp.

\bibitem[V]{V} 
van de Leur J. , A classification of contragredient Lie superalgebras of finite growth, Comm. Algebra, {\bf 17} (1989), 1815--1841.

\end{thebibliography}
\end{document}